\documentclass{article}[12pt]  
\usepackage[utf8x]{inputenc}
\usepackage{amsfonts}
\usepackage{graphics}
\usepackage{graphicx} 
\usepackage{color}
\usepackage[margin=3.cm]{geometry}
\usepackage{amsmath,amssymb}
\usepackage{color}
\usepackage[dvipsnames]{xcolor}
\usepackage{appendix}
\usepackage{amsthm}
\usepackage{array}
\usepackage{mathtools}
\renewcommand{\epsilon}{\varepsilon}

\graphicspath{{figures/}}

\DeclareGraphicsExtensions{.eps,.pdf,.png,.jpg,.JPG,.jpeg,.ps,.pdf_tex}
\numberwithin{equation}{section}

\usepackage{setspace,color}
\usepackage{amsmath,amssymb,mathrsfs}
\usepackage{graphicx,epstopdf}
\usepackage{xcolor}

\numberwithin{equation}{section}

\newcommand{\ihat}{\hat{\i}}
\newcommand{\jhat}{\hat{\j}}
\newcommand{\khat}{\hat{{\rm k}}}
\newcommand{\uhat}{\hat{{\rm u}}}
\newcommand{\vhat}{\hat{{\rm v}}}
\newcommand{\what}{\hat{{\rm w}}}

\begin{document}

 \title{Asymptotic approximations for the close evaluation of
  double-layer potentials}
\author{Camille Carvalho\footnote{Applied Mathematics Unit, School of Natural Sciences,
    University of California, Merced, 5200 North Lake Road, Merced, CA
    95343}  \and  Shilpa Khatri ${}^{\ast}$  \and  Arnold D. Kim ${}^{\ast}$}
   \maketitle

  \begin{abstract}
    When using the boundary integral equation method to solve a boundary value problem, the evaluation
    of the solution near the boundary is challenging to compute
    because the layer potentials that represent the solution
    are nearly-singular integrals. To address this close evaluation
    problem, we apply an asymptotic analysis of these nearly singular
    integrals and obtain an asymptotic approximation. 
	We derive the asymptotic approximation for the case of the double-layer potential in two and three dimensions,
    representing the solution of the interior Dirichlet problem for Laplace's equation.
    By doing so, we obtain an asymptotic approximation given
  by the Dirichlet data at the boundary point nearest to the interior
  evaluation point plus a nonlocal correction.
  We present numerical methods to compute this asymptotic approximation, and we demonstrate the efficiency and accuracy of the asymptotic
  approximation through several examples. These
  examples show that the asymptotic approximation is useful as it
  accurately approximates the close evaluation of the double-layer
  potential while requiring only modest computational resources.

\end{abstract}


\section{Introduction}

The close evaluation problem refers to the non-uniform error produced
by high-order quadrature rules used in boundary integral equation
methods. In particular, the close evaluation problem occurs when evaluating layer potentials at evaluation points close to the boundary. 
These high-order quadrature rules attain spectral accuracy
when computing the solution, represented by layer potentials, far from the boundary whereas they
incur a very large error when computing the solution close to the
boundary.  It is well understood that this growth in error is due to
the fact that the integrand of the layer potentials become
increasingly peaked as the point of evaluation approaches the
boundary. In fact, when the distance between the evaluation point and
its closest boundary point is smaller than the distance between
quadrature points on the boundary for a fixed-order quadrature rule,
the quadrature points do not adequately resolve the peak of the
integrand and therefore produce an $O(1)$ error.

Accurate evaluations of layer potentials close to the boundary of the domain are needed for a wide range of applications, including
the modeling of swimming micro-organisms, droplet suspensions, and blood cells in Stokes flow~\cite{smith2009boundary, barnett2015spectrally, marple2016fast,
  keaveny2011applying}, and to predict accurate measurements of the electromagnetic near-field in the field of
plasmonics~\cite{Maier07} for 
nano-antennas~\cite{akselrod2014probing,novotny2011antennas} and
sensors~\cite{mayer2008label,sannomiya2008situ}.   

Several computational methods have been developed to address this
close evaluation problem.  Schwab and
Wendland~\cite{schwab1999extraction} have developed a boundary
extraction method based on a Taylor series expansion of the layer
potentials.  Beale and Lai~\cite{beale2001method} have developed a
method that first regularizes the nearly singular kernel of the layer potential and then adds
corrections for both the discretization and the regularization. Beale
{\it et al.}~\cite{beale2016simple} have extended the regularization
method to three-dimensional problems. Helsing and
Ojala~\cite{helsing2008evaluation} developed a method that combines a
globally compensated quadrature rule and interpolation to achieve very
accurate results over all regions of the
domain. Barnett~\cite{barnett2014evaluation} has used surrogate local
expansions with centers placed near, but not on, the
boundary. Kl\"{o}ckner {\it et al.}~\cite{klockner2013quadrature}
introduced Quadrature By Expansion (QBX), which uses expansions about
accurate evaluation points far away from the boundary to compute
accurate evaluations close to it. There have been several subsequent
studies of QBX~\cite{epstein2013convergence, af2016fast,
  rachh2017fast, wala20183DQBX, af2017error} that have extended its
use and characterized its behavior.

Recently, the authors have applied asymptotic analysis to study the
close evaluation problem. For two-dimensional problems, the authors
developed a method that used matched asymptotic expansions for the
kernel of the layer potential~\cite{carvalho2018asymptotic}. In that
method, the asymptotic expansion that captures the peaked behavior of
the kernel (namely, the peaked behavior of the integrand of the layer potential) can be integrated exactly and the relatively smooth
remainder is integrated numerically, resulting in a highly accurate
method. For three-dimensional problems, the authors have developed a
simple, three-step method for computing layer
potentials~\cite{carvalho2018_3D}. This method involves first rotating
the spherical coordinate system used to compute the layer potential so
that the boundary point at which the integrand becomes singular is
aligned with the north pole. By studying the asymptotic behavior of
the integral, they found that integration with respect to the
azimuthal angle is a natural averaging operation that regularizes the
integral thereby allowing for a high-order quadrature rule to be used
for the integral with respect to the polar angle. This numerical
method was shown to achieve an error that decays quadratically with
the distance to the boundary provided that the underlying boundary
integral equation for the density is sufficiently resolved.

In this work, we carry out a complete asymptotic analysis of the
double-layer potential for the interior Dirichlet problem for
Laplace's equation in two and three dimensions. By doing so, we derive asymptotic approximations
for the close evaluation of the double-layer potential. These asymptotic approximations provide valuable
insight into the inherent challenges of the close evaluation problem
and an explicit method to address it. We find that the leading-order asymptotic
behavior of the double-layer potential in the close evaluation limit
is given by the Dirichlet data at the boundary point closest to the
evaluation point plus a nonlocal correction. It is the nonlocal
correction that has made the close evaluation problem challenging to
address. Since this asymptotic analysis explicitly finds this nonlocal
correction, we are able to develop a simple and accurate numerical
method to compute the double-layer potential and thus, address the close evaluation problem
systematically. 
We compute several numerical
examples using the asymptotic approximations to evaluate their
efficacy and accuracy.

The asymptotic analysis used here to study the close evaluation
problem also provides valuable insight (and useful asymptotic approximations) for other problems. In particular there is an interesting connection with forward-peaked scattering in radiative transfer, which is used
to describe the multiple scattering of light~\cite{chandrasekhar1960,
  ishimaru}. Forward-peaked scattering is an important problem for
several applications, and is challenging to study. We draw this
connection and apply the asymptotic analysis developed in this paper to forward-peaked scattering
in radiative transfer.

The remainder of this paper is as follows. We precisely define the
close evaluation problem for the double-layer potential in Section
\ref{sec:close-eval}. We compute the leading-order asymptotic behavior
of the double-layer potential in two and three dimensions in Sections
\ref{sec:asymptotics2D} and \ref{sec:asymptotics3D}, respectively.
We describe numerical methods to evaluate the asymptotic
approximations for the close evaluation of the double-layer potential
in Section \ref{sec:numerics}. We give several examples demonstrating
the accuracy of this numerical method in Section
\ref{sec:results}. Section \ref{sec:RTE} describes the connection between the close
evaluation problem and forward-peaked scattering in radiative
transfer. Section \ref{sec:conclusion} gives our
conclusions. The Appendix provides details of the computations for the three-dimensional case: Appendix \ref{sec:rotation} gives details of how we
rotate spherical integrals and Appendix \ref{sec:spherical-laplacian} gives a
useful derivation of the spherical Laplacian.

\section{Close evaluation of the double-layer potential}
\label{sec:close-eval}

Consider a simply connected, open set, denoted by
$D \subset \mathbb{R}^{n}$ with $n = 2, 3$, with an analytic close
boundary, $B$, and let $\bar{D} = D \cup B$. Given some smooth data
$f$, we write the function $u \in C^{2}(D) \cap C^{1}(\bar{D})$
satisfying the interior Dirichlet problem,
\begin{subequations}
  \begin{gather}
    \varDelta u = 0 \quad \text{in $D$},\\
    u = f \quad \text{on $B$},
  \end{gather}
  \label{eq:dirichlet}
\end{subequations}
as the double-layer potential,
\begin{equation}
  u(x) = \frac{1}{2^{n-1} \pi} \int_{B} \frac{\nu_{y} \cdot ( x
    - y )}{|x - y|^{n}} \mu(y) \mathrm{d}\sigma_{y}, \quad x \in D,
  \quad n = 2, 3.
  \label{eq:DLP}
\end{equation}
Here, $\nu_{y}$ denotes the unit outward normal at $y \in B$,
$\mathrm{d}\sigma_{y}$ denotes the boundary element, and $\mu$, the density, is a
continuous function. This double-layer potential
satisfies the following jump relation~\cite{guenther1996partial},
\begin{equation}
  \lim_{\substack{x \to y^{\star} \in B\\x \in D}} u(x) = u(y^{\star}) -
  \frac{1}{2} \mu(y^{\star}).
  \label{eq:jump}
\end{equation}
By requiring that $u$ satisfies \eqref{eq:dirichlet}, we find that, in
light of jump relation \eqref{eq:jump}, $\mu$ must satisfy
\begin{equation}
  \frac{1}{2^{n-1} \pi} \int_{B} \frac{\nu_y \cdot (
    y^{\star} - y )}{|y^{\star} - y|^{n}} \mu(y) \mathrm{d}\sigma_{y}
  -\frac{1}{2} \mu(y^{\star}) = f(y^{\star}), \quad y^{\star}
  \in B,
  \label{eq:bie}
\end{equation}
the boundary integral equation for $\mu$.

\begin{figure}[h!]
  \centering
  \def\svgwidth{0.8\columnwidth} 
  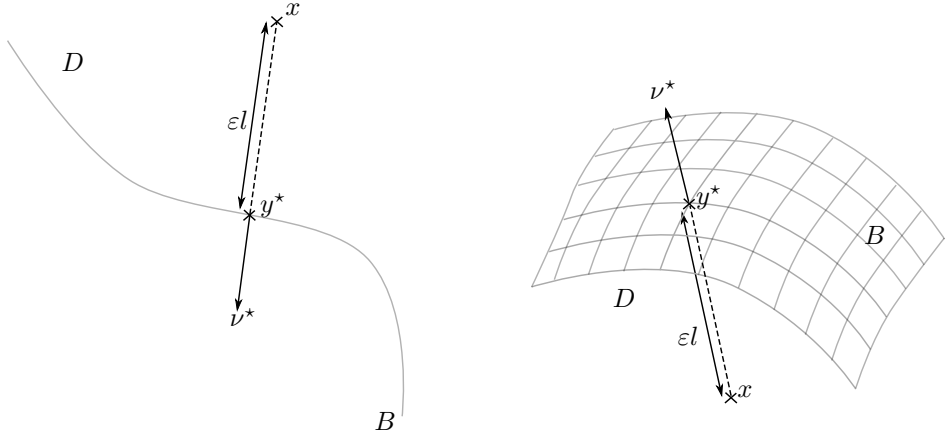
  \caption{Sketch of the quantities introduced in
    \eqref{eq:closept} to study evaluation points close to the
    boundary in 2D (left) and in 3D (right).}
  \label{fig:sketch}
\end{figure}

Here, we seek to evaluate \eqref{eq:DLP} at points close to the
boundary.  To define a close evaluation point precisely, let
$0 < \epsilon \ll 1$ denote a small, dimensionless parameter, and
consider
\begin{equation}
  x = y^{\star} - \epsilon \ell \nu^{\star},
  \label{eq:closept}
\end{equation}
with $y^{\star} \in B$ denoting the closest point to $x$ on the
boundary, $\nu^{\star}$ denoting the unit, outward
normal at $y^{\star}$, and $\ell$ denoting a characteristic length of
the problem such as the signed (2D) or mean (3D) curvature at
$y^{\star}$ (see Fig. \ref{fig:sketch}).  

Because the solution of \eqref{eq:dirichlet}
continuously approaches its boundary data from within $D$, we write
\begin{equation}
u (x)  = u(y^{\star} - \epsilon \ell \nu^{\star}) = f(y^{\star})
  + \epsilon U(y^{\star};\epsilon). 
  \label{eq:constantapprox}
\end{equation}
To determine an expression for $U$, we substitute \eqref{eq:DLP}
evaluated at \eqref{eq:closept} for
$u(y^{\star} - \epsilon \ell \nu^{\star})$ and \eqref{eq:bie} for
$f(y^{\star})$ into \eqref{eq:constantapprox}, and find that
\begin{equation}
  U(y^{\star};\epsilon) = \epsilon^{-1} \left[ \frac{1}{2} \mu(y^{\star}) +
    \frac{1}{2^{n-1} \pi} \int_{B} \left[
      \frac{\nu_{y} \cdot ( y_{d} - \epsilon \ell \nu^{\star})}{|y_{d} -
        \epsilon \ell \nu^{\star} |^{n}} - \frac{\nu_{y} \cdot
        y_{d}}{|y_{d}|^{n}} \right] \mu(y) \mathrm{d}\sigma_{y} \right],
  \label{eq:constantapproxwithbie}
\end{equation}
where we have introduced the notation, $y_{d} = y^{\star} - y$.  Next,
we make use of Gauss' theorem~\cite{guenther1996partial}
\begin{equation}
  \frac{1}{2^{n-1} \pi} \int_{B} \frac{\nu_{y} \cdot ( x - y
    )}{|x - y|^{n}} \mathrm{d}\sigma_{y} = \begin{cases}
    -1 & x \in D,\\
    -\frac{1}{2} & x \in B,\\
    \,\,\,\, 0 & x \not\in \bar{D}, \end{cases}
  \label{eq:gauss}
\end{equation}
to write
\begin{equation}
  \frac{1}{2} \mu(y^{\star}) = - \frac{1}{2^{n-1} \pi} \int_{B}
  \frac{\nu_{y} \cdot ( y_{d} - \epsilon \ell
        \nu^{\star})}{|y_{d} - \epsilon \ell
    \nu^{\star} |^{n}} \mu(y^{\star}) \mathrm{d}\sigma_{y}
  + \frac{1}{2^{n-1} \pi} \int_{B} \frac{\nu_{y} \cdot 
    y_{d}}{|y_{d} |^{n}} \mu(y^{\star})
  \mathrm{d}\sigma_{y}.
  \label{eq:gaussidentity}
\end{equation}
Substituting \eqref{eq:gaussidentity} into
\eqref{eq:constantapproxwithbie} yields
\begin{equation}
  U(y^{\star};\epsilon) = \frac{1}{2^{n-1} \pi} \int_{B}
  \epsilon^{-1} \left[ \frac{\nu_{y} \cdot ( y_{d} - \epsilon
      \ell \nu^{\star})}{|y_{d} - \epsilon \ell \nu^{\star} |^{n}} -
    \frac{\nu_{y} \cdot y_{d}}{|y_{d}|^{n}} \right] [ \mu(y) -
  \mu(y^{\star}) ]\mathrm{d}\sigma_{y}.
  \label{eq:correction}
\end{equation}

We seek to determine the asymptotic expansion of $U$ given in
\eqref{eq:correction} in the limit as $\epsilon \to 0^{+}$.  To
determine this asymptotic expansion, we make use of explicit
parametrizations for $B$. Therefore, we consider the two and three-dimensional problems separately.

\section{Asymptotic analysis in two dimensions}
\label{sec:asymptotics2D}

Suppose $B$ is an analytic, closed curve on the plane. For that case,
we introduce the parameter $t \in [-\pi,\pi]$ such that $y = y(t)$ and
$y^{\star} = y(0)$.  In terms of this parameterization,
\eqref{eq:correction} is given by
\begin{equation}
  U(y^{\star};\epsilon) = \frac{1}{2 \pi} \int_{-\pi}^{\pi}
  \epsilon^{-1} K(t;\epsilon) [ \tilde{\mu}(t) -
  \tilde{\mu}(0) ] \mathrm{d}t,
\end{equation}
with $\tilde{\mu}(t) = \mu(y(t))$,
$\tilde{\mu}(0) = \mu(y(0)) = \mu(y^{\star})$, and
\begin{equation}
  K(t;\epsilon) = \left[ \frac{\tilde{\nu}(t) \cdot ( y_{d}(t)
      - \epsilon \ell \nu^{\star})}{|y_{d}(t) -
      \epsilon \ell \nu^{\star} |^{2}} -
    \frac{\tilde{\nu}(t) \cdot y_{d}(t)}{|y_{d}(t)|^{2}}
  \right] J(t),
  \label{eq:kernel2D}
\end{equation}
with $\tilde{\nu}(t) = \nu(y(t))$, $y_{d}(t) = y(0) - y(t)$, and
$J(t) = | y'(t) |$. Note that $\nu^\star = \tilde{\nu}(0)$.

To determine the asymptotic expansion for $U$, we introduce the small
parameter $\delta$ satisfying $0 < \epsilon \ll \delta \ll 1$, and
write
\begin{equation}
  U(y^{\star};\epsilon) = U^{\text{in}}(y^{\star};\epsilon,\delta) +
  U^{\text{out}}(y^{\star};\epsilon,\delta).
  \label{eq:3.3}
\end{equation}
Here, the inner expansion, $U^{\text{in}}$, is given by
\begin{equation}
  U^{\text{in}}(y^{\star};\epsilon,\delta) = \frac{1}{2\pi}
  \int_{-\delta/2}^{\delta/2} \epsilon^{-1} K(t;\epsilon) \left[
    \tilde{\mu}(t) - \tilde{\mu}(0) \right] \mathrm{d}t,
  \label{eq:Uin2D}
\end{equation}
and the outer expansion, $U^{\text{out}}$, is given by
\begin{multline}
  U^{\text{out}}(y^{\star};\epsilon,\delta) = \frac{1}{2\pi}
  \int_{-\pi}^{-\delta/2} \epsilon^{-1} K(t;\epsilon) \left[
    \tilde{\mu}(t) - \tilde{\mu}(0) \right] \mathrm{d}t + \frac{1}{2\pi} \int_{\delta/2}^{\pi} \epsilon^{-1} K(t;\epsilon)
  \left[ \tilde{\mu}(t) - \tilde{\mu}(0) \right] \mathrm{d}t.
  \label{eq:Uout2D}
\end{multline}
The inner expansion involves integration over a small
  portion of the boundary about $y^{\star}$, whereas the outer
  expansion involves integration over the remaining portion of the
  boundary.

We determine the leading-order asymptotic behaviors of $U^{\text{in}}$
and $U^{\text{out}}$ in the sections below. Then, we combine those
results to obtain the asymptotic approximation for the double-layer
potential in two dimensions, and discuss higher-order asymptotic
approximations. We have developed a \textit{Mathematica} notebook that
contains the presented calculations, available in a GitHub repository
\cite{CKK-2018Codes}.

\subsection{Inner expansion}\label{ssec:inner}

To determine the leading-order asymptotic behavior of $U^{\text{in}}$,
we substitute $t = \epsilon T$ into \eqref{eq:Uin2D}, and obtain
\begin{equation}
  U^{\text{in}}(y^{\star};\epsilon,\delta) = \frac{1}{2\pi}
  \int_{-\delta/2\epsilon}^{\delta/2\epsilon} K(\epsilon T;\epsilon)
  \left[ \tilde{\mu}(\epsilon T) - \tilde{\mu}(0) \right]
  \mathrm{d}T.
  \label{eq:3.6}
\end{equation}
Recognizing that $\tilde{\nu}(\epsilon T) = \nu^{\star} + O(\epsilon)$ and $y_{d}(\epsilon T) = - \epsilon T y'(0) +
  O(\epsilon^{2})$ with $\nu^{\star} \cdot y'(0) = 0$, we find that by expanding
$K(\epsilon T;\epsilon)$ about $\epsilon = 0$ that
\begin{equation}
  K(\epsilon T;\epsilon) = - \frac{\epsilon^{-1} \ell J(0)}{
    T^{2} J^{2}(0) + \ell^{2}} + O(1).
\end{equation}
Using the fact that this leading-order behavior is even in $T$, and
expanding $\tilde{\mu}$ about $\epsilon = 0$, we can substitute into
\eqref{eq:3.6} to get, after expanding about $\delta = 0$,
\begin{align}
  \begin{split}
    U^{\text{in}}&(y^{\star};\epsilon,\delta) = \frac{1}{2\pi}
    \int_{-\delta/2\epsilon}^{\delta/2\epsilon} \left[ -
      \frac{\epsilon^{-1} \ell J(0)}{T^{2} J^{2}(0) + \ell^{2}} + O(1)
    \right] \left[ \tilde{\mu}(\epsilon T) - \tilde{\mu}(0) \right]
    \mathrm{d}T\\
    &= \frac{1}{2\pi} \int_{0}^{\delta/2\epsilon} \left[ -
      \frac{\epsilon^{-1} \ell J(0)}{T^{2} J^{2}(0) + \ell^{2}} + O(1)
    \right] \left[ \tilde{\mu}(\epsilon T) + \tilde{\mu}(-\epsilon T)
      - 2 \tilde{\mu}(0) \right] \mathrm{d}T\\
    &=\frac{1}{2\pi} \int_{0}^{\delta/2\epsilon} \left[ -
      \frac{\epsilon^{-1} \ell J(0)}{T^{2} J^{2}(0) + \ell^{2}} + O(1)
    \right] \left[ \epsilon^{2} T^{2} \tilde{\mu}''(0) +
      O(\epsilon^{4}) \right] \mathrm{d}T\\
    &=\frac{1}{2\pi} \int_{0}^{\delta/2\epsilon} \left[ -
      \frac{\epsilon T^{2} \ell J(0)}{T^{2} J^{2}(0) + \ell^{2}}
      \tilde{\mu}''(0) +
      O(\epsilon^{2}) \right] \mathrm{d}T\\
    &= - \frac{\delta \ell}{4 \pi J(0)} \tilde{\mu}''(0) +
    O(\epsilon).
  \end{split}
  \label{eq:Uin2Dleadingorder}
\end{align}
This result gives the leading-order asymptotic behavior of
$U^{\text{in}}$.

\subsection{Outer expansion} To determine the leading-order asymptotic
behavior of $U^{\text{out}}$, we expand $K(t;\epsilon)$ about
$\epsilon = 0$ and find that
$K(t;\epsilon) = [ \epsilon K_{1}(t) + O(\epsilon^{2}) ] J(t)$, with
\begin{equation}
  K_{1}(t) = \ell \frac{2 ( \tilde{\nu}(t) \cdot y_{d}(t))
    ( \nu^{\star} \cdot y_{d}(t)) - \tilde{\nu}(t) \cdot
    \nu^{\star} |y_{d}(t)|^{2}}{|y_{d}(t)|^{4}}.
    \label{eq:K1-2D}
\end{equation}
Substituting this expansion into \eqref{eq:Uout2D}, we find that
\begin{multline}
  U^{\text{out}}(y^{\star};\epsilon,\delta) = \frac{1}{2\pi}
  \int_{-\pi}^{-\delta/2} K_{1}(t) \left[
    \tilde{\mu}(t) - \tilde{\mu}(0) \right] J(t) \mathrm{d}t  + \frac{1}{2\pi} \int_{\delta/2}^{\pi} K_{1}(t) \left[
    \tilde{\mu}(t) - \tilde{\mu}(0) \right] J(t) \mathrm{d}t +
  O(\epsilon).
  \label{eq:Uout2Dintermediate}
\end{multline}
To eliminate $\delta$ from the integration limits, we rewrite
\eqref{eq:Uout2Dintermediate} as
\begin{equation}
  U^{\text{out}}(y^{\star};\epsilon,\delta) = \frac{1}{2\pi}
  \int_{-\pi}^{\pi} K_{1}(t) \left[ \tilde{\mu}(t) - \tilde{\mu}(0)
  \right] J(t) \mathrm{d}t - V^{\text{out}}(y^{\star};\epsilon,\delta)
  + O(\epsilon)
  \label{eq:Uout2Dintermediate2}
\end{equation}
with
\begin{equation}
  V^{\text{out}}(y^{\star};\epsilon,\delta) = \frac{1}{2\pi}
  \int_{-\delta/2}^{\delta/2} K_{1}(t) \left[ \tilde{\mu}(t) -
    \tilde{\mu}(0) \right] J(t) \mathrm{d}t.
  \label{eq:3.12}
\end{equation}
To determine the leading-order behavior for $V^{\text{out}}$, we proceed exactly as in section \ref{ssec:inner}. We
substitute $t = \epsilon T$ into \eqref{eq:3.12} and obtain
\begin{equation}
  V^{\text{out}}(y^{\star};\epsilon,\delta) = \frac{1}{2\pi}
  \int_{-\delta/2\epsilon}^{\delta/2\epsilon} K_{1}(\epsilon T) \left[
    \tilde{\mu}(\epsilon T) - \tilde{\mu}(0) \right] 
  J(\epsilon T) \epsilon \mathrm{d}T
  \label{eq:Vout2D}
\end{equation}
Again, by recognizing that
$\tilde{\nu}(\epsilon T) = \nu^{\star} + O(\epsilon)$ and $y_{d}(\epsilon T) = - \epsilon T y'(0) +
  O(\epsilon^{2})$ with $\nu^{\star} \cdot y'(0) = 0$, we find that
\begin{equation}
  K_{1}(\epsilon T) = - \frac{\epsilon^{-2} \ell}{T^{2}
    J^{2}(0)} + O(\epsilon^{-1}).
\end{equation}
Using the fact that this leading-order behavior is even in $T$, and
that $J(\epsilon T) = J(0) + O(\epsilon)$, when we substitute it into
\eqref{eq:Vout2D}, we find, after expanding about $\delta = 0$, that
\begin{align}\label{eq:Vout2Dbis}
  \begin{split}
    V^{\text{out}}&(y^{\star};\epsilon,\delta) = \frac{1}{2\pi}
    \int_{-\delta/2\epsilon}^{\delta/2\epsilon} \left[ -
      \frac{\epsilon^{-2} \ell}{T^{2} J^{2}(0)} + O(\epsilon^{-1})
    \right] \left[ \tilde{\mu}(\epsilon T) - \tilde{\mu}(0) \right]
    J(\epsilon T) \epsilon \mathrm{d}T\\
    &= \frac{1}{2\pi} \int_{0}^{\delta/2\epsilon} \left[ -
      \frac{\epsilon^{-2} \ell}{T^{2} J^{2}(0)} + O(\epsilon^{-1})
    \right] \left[ \tilde{\mu}(\epsilon T) + \tilde{\mu}(-\epsilon T)
      - 2
      \tilde{\mu}(0) \right] \left[J(0) + O(\epsilon)\right] \epsilon
    \mathrm{d}T\\
    &= \frac{1}{2\pi} \int_{0}^{\delta/2\epsilon} \left[ -
      \frac{\epsilon^{-1} \ell}{T^{2} J(0)} + O(1) \right] \left[
      \epsilon^{2} T^{2} \tilde{\mu}''(0) + O(\epsilon^{4})
    \right] \mathrm{d}T\\
    &= \frac{1}{2\pi} \int_{0}^{\delta/2\epsilon} \left[
      -\frac{\epsilon \ell}{J(0)} \tilde{\mu}''(0) + O(\epsilon^{2})
    \right]
    \mathrm{d}T\\
    &= - \frac{\delta \ell}{4 \pi J(0)} \tilde{\mu}''(0)+ O(\epsilon).
  \end{split}
\end{align}
Substituting this result into \eqref{eq:Uout2Dintermediate2}, we find
that the leading-order asymptotic behavior for $U^{\text{out}}$ is
given by
\begin{equation}
  U^{\text{out}}(y^{\star};\epsilon,\delta) = \frac{1}{2\pi}
  \int_{-\pi}^{\pi} K_{1}(t) \left[ \tilde{\mu}(t) - \tilde{\mu}(0)
  \right] J(t) \mathrm{d}t + \frac{\delta \ell}{4 \pi J(0)}
  \tilde{\mu}''(0) + O(\epsilon).
  \label{eq:Uout2Dleadingorder}
\end{equation}

\subsection{Two-dimensional asymptotic approximation}

We obtain an asymptotic approximation for $U$ by summing the
leading-order behaviors obtained for $U^{\text{in}}$ and
$U^{\text{out}}$ given in \eqref{eq:Uin2Dleadingorder} and
\eqref{eq:Uout2Dleadingorder}, respectively. Substituting that result
into \eqref{eq:constantapprox}, we obtain the following asymptotic
approximation, 
\begin{equation}
  u(y^{\star} - \epsilon \ell \nu^{\star}) = f(y^{\star}) + \epsilon
  L_{1}[ \mu ] + O(\epsilon^{2}),
  \label{eq:asymptotic2D}
\end{equation}
with
\begin{equation}
  L_{1}[\mu] = \frac{1}{2\pi} \int_{-\pi}^{\pi} K_{1}(t) \left[
    \tilde{\mu}(t) - \tilde{\mu}(0) \right] J(t) \mathrm{d}t,
  \label{eq:L1-2D}
\end{equation}
where $K_{1}$ is given by \eqref{eq:K1-2D}. Naturally, the obtained asymptotic approximation doesn't depend on the arbitrary parameter $\delta$.

Asymptotic approximation \eqref{eq:asymptotic2D} gives an explicit
approximation for the close evaluation of the double-layer potential
in two dimensions. According to the asymptotic analysis, the error of
this approximation is $O(\epsilon^{2})$. It gives the double-layer
potential as the Dirichlet data at the boundary point $y^{\star}$
closest to the evaluation point $x$ plus a nonlocal correction. This
nonlocal correction is consistent with the fact that solutions to
elliptic partial differential equations have a global dependence on
their boundary data. The leading-order asymptotic expansion indicates that the nonlocal correction only comes from the outer expansion, and the inner expansion doesn't contribute to the lower order terms.

\subsection{Higher-order asymptotic approximations}

By continuing on to higher order terms in the expansions for
$U^{\text{in}}$ and $U^{\text{out}}$, we can obtain higher-order
asymptotic approximations.  This process will not be demonstrated
here, because the calculations become unwieldy, details can be found in the developed \textit{Mathematica} notebook~\cite{CKK-2018Codes}. The result from these calculations is
\begin{multline}
  u(y^{\star} - \epsilon \ell \nu^{\star}) = f(y^{\star}) + \epsilon
  L_{1}[ \mu ] + \epsilon^{2} \left[ L_{2}[\mu] - \frac{\ell^{2} y'(0) \cdot
      y''(0)}{4 J^{4}(0)} \tilde{\mu}'(0) + \frac{\ell^{2}}{4
      J^{2}(0)} \tilde{\mu}''(0) \right] + O(\epsilon^{3}),
  \label{eq:asymptotic2D-higherorder}
\end{multline}
with $L_{1}[\mu]$ given in \eqref{eq:L1-2D}, and
\begin{equation}
  L_{2}[\mu] = \frac{1}{2\pi} \int_{-\pi}^{\pi} K_{2}(t) \left[
    \tilde{\mu}(t) - \tilde{\mu}(0) \right] J(t) \mathrm{d}t,
\end{equation}
where
\begin{equation}
  K_{2}(t) = \ell^{2} \frac{(\nu \cdot y_{d}) \left[ 4 ( \nu^{\star}
      \cdot y_{d} )^{2} - |y_{d}|^{2} \right] -2 |y_{d}|^{2} (\nu
    \cdot \nu^{\star})( \nu^{\star} \cdot y_{d})}{|y_{d}|^{6}}.
  \label{eq:K2-2D}
\end{equation}
This asymptotic approximation has an error that is
$O(\epsilon^{3})$. In addition to nonlocal terms, this approximation
includes local contributions made by first and second derivatives of
the density, $\mu$, evaluated at the boundary point $y^{\star}$. The local contributions come from the inner expansion.

\section{Asymptotic analysis in three dimensions}
\label{sec:asymptotics3D}

Suppose $B$ is an analytic, closed, and oriented surface. We introduce
the parameters $s \in [0,\pi]$ and $t \in [-\pi,\pi]$ such that
$y = y(s,t)$ and $y^{\star} = y(0,\cdot)$. In terms of this
parameterization, \eqref{eq:correction} is given by
\begin{equation}
  U(y^{\star};\epsilon) = \frac{1}{4\pi} \int_{-\pi}^{\pi}
  \int_{0}^{\pi} \epsilon^{-1} K(s,t;\epsilon) \left[
    \tilde{\mu}(s,t) - \tilde{\mu}(0,\cdot) \right] \sin(s) \mathrm{d}s
  \mathrm{d}t,
\end{equation}
with $\tilde{\mu}(s,t) = \mu(y(s,t))$, $\tilde{\mu}(0,\cdot) =
\mu(y(0,\cdot))$, and 
\begin{equation}
  K(s,t;\epsilon) = \left[ \frac{\tilde{\nu}(s,t) \cdot ( y_{d}(s,t) -
      \epsilon \ell \nu^{\star})}{| y_{d}(s,t) - \epsilon \ell
      \nu^{\star} |^{3}} - \frac{\tilde{\nu}(s,t) \cdot
      y_{d}(s,t)}{| y_{d}(s,t) |^{3}} \right] J(s,t),
\end{equation} 
with $\tilde{\nu}(s,t) = \nu_y$,
$y_{d}(s,t) = y(0,\cdot) - y(s,t)$, $J(s,t) = | y_{s}(s,t) \times y_{t}(s,t) |/\sin(s)$. Note that $ \nu^\star = \tilde{\nu}(0, \cdot)$.

Just as we have done for the two dimensional problem, we introduce the
small parameter $\delta$ satisfying $0 < \epsilon \ll \delta \ll 1$,
and write
\begin{equation}
  U(y^{\star};\epsilon) = U^{\text{in}}(y^{\star};\epsilon,\delta) +
  U^{\text{out}}(y^{\star};\epsilon,\delta).
  \label{eq:4.3}
\end{equation}
Here, the inner expansion is given by
\begin{equation}
  U^{\text{in}}(y^{\star};\epsilon,\delta) = \frac{1}{4\pi}
  \int_{-\pi}^{\pi} \int_{0}^{\delta} \epsilon^{-1} K(s,t;\epsilon)
  \left[ \tilde{\mu}(s,t) - \tilde{\mu}(0,\cdot) \right]  \sin(s)\mathrm{d}s
  \mathrm{d}t,
  \label{eq:Uin3D}
\end{equation}
and the outer expansion is given by
\begin{equation}
  U^{\text{out}}(y^{\star};\epsilon,\delta) = \frac{1}{4\pi}
  \int_{-\pi}^{\pi} \int_{\delta}^{\pi} \epsilon^{-1}
  K(s,t;\epsilon) \left[ \tilde{\mu}(s,t) - \tilde{\mu}(0,\cdot)
  \right] \sin(s)\mathrm{d}s \mathrm{d}t.
  \label{eq:Uouter3D}
\end{equation}
Again, we determine the leading-order asymptotic behaviors for
$U^{\text{in}}$ and $U^{\text{out}}$ separately. Then, we combine
those results to obtain an asymptotic approximation for the close
evaluation of the double-layer potential in three dimensions and discuss higher-order asymptotic approximations. Some details can be found in the developed \textit{Mathematica} notebook available on GitHub~\cite{CKK-2018Codes}.

\subsection{Inner expansion} \label{ssec:inner3D}

To find the leading-order asymptotic behavior of $U^{\text{in}}$, we
substitute $s = \epsilon S$ into \eqref{eq:Uin3D}, and obtain
\begin{equation}
  U^{\text{in}}(y^{\star};\epsilon,\delta) = \frac{1}{4\pi}
  \int_{-\pi}^{\pi} \int_{0}^{\delta/\epsilon} K(\epsilon
  S,t;\epsilon) \left[ \tilde{\mu}(\epsilon S,t) -
    \tilde{\mu}(0,\cdot) \right] \sin(\epsilon S) \mathrm{d}S
  \mathrm{d}t.
\end{equation}
Recognizing that
$\tilde{\nu}(\epsilon S, t) = \nu^{\star} + O(\epsilon)$, and $y_{d}(\epsilon S, t) = - \epsilon S y_{s}(0,\cdot) +
  O(\epsilon^{2})$ with the vector $y_{s}(0,\cdot)$ lying on the plane tangent to $B$ at
$y^{\star}$, we find by expanding $K(\epsilon S, t ; \epsilon)$ about
$\epsilon = 0$ that
\begin{equation}
  K(\epsilon S, t ; \epsilon) = - \frac{\epsilon^{-2} \ell J(0,\cdot)}{(S^{2} |
    y_{s}(0,\cdot) |^{2} +  \ell^{2})^{3/2}} + O(\epsilon^{-1}).
  \end{equation}
Since this leading-order asymptotic behavior for $K(\epsilon S, t; \epsilon)$ is
independent of $t$, we write
\begin{multline}
  U^{\text{in}}(y^{\star};\epsilon,\delta) = \frac{1}{4\pi}
  \int_{0}^{\pi} \int_{0}^{\delta/\epsilon} \left[ -
    \frac{\epsilon^{-2} \ell J(0,\cdot)}{(S^{2} | y_{s}(0,\cdot) |^{2}
      + \ell^{2})^{3/2}} + O(\epsilon^{-1}) \right]\\
  \left[ \tilde{\mu}(\epsilon S,t) + \tilde{\mu}(\epsilon S, t + \pi)
    - 2 \tilde{\mu}(0,\cdot) \right] \sin(\epsilon S) \mathrm{d}S
  \mathrm{d}t.
  \label{eq:4.8}
\end{multline}
Next, we use the regularity of $\tilde{\mu}$ over the north pole to
substitute
$\tilde{\mu}(\epsilon S, t + \pi) = \tilde{\mu}(-\epsilon S, t)$, so
that
\begin{align}
  \begin{split}
    \tilde{\mu}(\epsilon S,t) + \tilde{\mu}(\epsilon S, t + \pi) - 2
    \tilde{\mu}(0,\cdot) &= \tilde{\mu}(\epsilon S,t) +
    \tilde{\mu}(-\epsilon S,t) - 2 \tilde{\mu}(0,\cdot)\\
    &= \epsilon^{2} S^{2} \tilde{\mu}_{ss}(0,\cdot) + O(\epsilon^{4}).
  \end{split}
      \label{eq:mu_ss}
\end{align}
Thus, we find after substituting \eqref{eq:mu_ss} and
$\sin(\epsilon S) = \epsilon S + O(\epsilon^{3})$ into \eqref{eq:4.8}
that
\begin{align}
  \begin{split}
    U^{\text{in}}(y^{\star};\epsilon,\delta) &= \frac{1}{4\pi}
    \int_{0}^{\pi} \int_{0}^{\delta/\epsilon} \left[ - \frac{\epsilon
        S^{3} \ell J(0,\cdot)}{(S^{2} | y_{s}(0,\cdot) |^{2} +
        \ell^{2})^{3/2}} \tilde{\mu}_{ss}(0,\cdot) + O(\epsilon^{2})
    \right] \mathrm{d}S \mathrm{d}t\\
    &= - \frac{\ell J(0,\cdot)}{8} \varDelta_{S^{2}} \mu(y^{\star})
    \int_{0}^{\delta/\epsilon} \left[ \frac{\epsilon S^{3}}{(S^{2} |
        y_{s}(0,\cdot) |^{2} + \ell^{2})^{3/2}} + O(\epsilon^{2})
    \right] \mathrm{d}S,
    \end{split}
\end{align}
where we have used the fact that
\begin{equation}
  \frac{1}{\pi} \int_{0}^{\pi} \tilde{\mu}_{ss}(0,\cdot)\mathrm{d}t =
  \frac{1}{2} \varDelta_{S^{2}} \mu(y^{\star}),
  \label{eq:spherical-laplacian}
\end{equation}
with $\varDelta_{S^{2}} \mu(y^{\star})$ denoting the spherical
Laplacian of $\mu$ evaluated at $y^{\star}$ (see Appendix
\ref{sec:spherical-laplacian}). Furthermore, when expanding about $\delta = 0$ we have
\begin{equation}
  \int_{0}^{\delta/\epsilon} \left[ \frac{\epsilon S^{3}}{(S^{2}
    |y_{s}(0,\cdot)|^{2} + \ell^{2})^{3/2}} + O(\epsilon^{2}) \right]
  \mathrm{d}S  = \frac{\delta}{|y_{s}(0,\cdot)|^{3}} + O(\epsilon),
\end{equation}
we determine that
\begin{equation}
  U^{\text{in}}(y^{\star};\epsilon,\delta) = - \frac{\delta \ell
    J(0,\cdot)}{8 |y_{s}(0,\cdot)|^{3}} \varDelta_{S^{2}}
  \mu(y^{\star}) + O(\epsilon).
  \label{eq:Uinner3D}
\end{equation}
This result gives the leading-order asymptotic behavior of
$U^{\text{in}}$.

\subsection{Outer expansion}

To determine the leading-order asymptotic behavior of
$U^{\text{out}}$, we expand $K(s,t;\epsilon)$ about $\epsilon = 0$ and
find $K(s,t;\epsilon) = $ $\left[ \epsilon K_{1}(s,t) + O(\epsilon^{2}) \right] J(s,t)$,
with
\begin{equation}
  K_{1}(s,t) = \ell \frac{3 ( \tilde{\nu}(s,t) \cdot y_{d}(s,t) ) ( 
    \nu^{\star} \cdot y_{d}(s,t) ) - |y_{d}(s,t)|^{2} \tilde{\nu}(s,t) \cdot
    \nu^{\star}}{| y_{d}(s,t) |^{5}}.
  \label{eq:K1-3D}
\end{equation}
Substituting this expansion into \eqref{eq:Uouter3D}, we obtain
\begin{equation}
  U^{\text{out}}(y^{\star};\epsilon,\delta) = \frac{1}{4\pi}
  \int_{-\pi}^{\pi} \int_{\delta}^{\pi} K_{1}(s,t) \left[
    \tilde{\mu}(s,t) - \tilde{\mu}(0,\cdot)
  \right] J(s,t)\sin(s) \mathrm{d}s \mathrm{d}t + O(\epsilon).
  \label{eq:4.16}
\end{equation}
To eliminate $\delta$ as a limit of integration in \eqref{eq:4.16}, we
write
\begin{multline}
  U^{\text{out}}(y^{\star};\epsilon,\delta) = \frac{1}{4\pi}
  \int_{-\pi}^{\pi} \int_{0}^{\pi} K_{1}(s,t) \left[ \tilde{\mu}(s,t)
    - \tilde{\mu}(0,\cdot)
  \right] J(s,t)\sin(s) \mathrm{d}s \mathrm{d}t  - V^{\text{out}}(y^{\star};\epsilon,\delta) + O(\epsilon),
  \label{eq:Uouter3D-expand}
\end{multline}
with
\begin{equation}
  V^{\text{out}}(y^{\star};\epsilon,\delta) = \frac{1}{4\pi}
  \int_{-\pi}^{\pi} \int_{0}^{\delta} K_{1}(s,t) \left[
    \tilde{\mu}(s,t) - \tilde{\mu}(0,\cdot)
  \right] J(s,t)\sin(s) \mathrm{d}s \mathrm{d}t.
  \label{eq:4.17}
\end{equation}

To determine the leading-order asymptotic behavior of
$V^{\text{out}}(y^{\star};\epsilon,\delta)$, we proceed as in section \ref{ssec:inner3D}. We substitute
$s = \epsilon S$ into \eqref{eq:4.17}, and obtain
\begin{equation}
  V^{\text{out}}(y^{\star};\epsilon,\delta) = \frac{1}{4\pi}
  \int_{-\pi}^{\pi} \int_{0}^{\delta/\epsilon} K_{1}(\epsilon S,t) \left[
    \tilde{\mu}(\epsilon S,t) - \tilde{\mu}(0,\cdot)
  \right]J(\epsilon S,t) \sin(\epsilon S) \epsilon \mathrm{d}S \mathrm{d}t.
\end{equation}
Recognizing that
$\tilde{\nu}(\epsilon S, t) = \nu^{\star} + O(\epsilon)$, and $y_{d}(\epsilon S, t) = - \epsilon S y_{s}(0,\cdot) +
  O(\epsilon^{2})$ with the vector $y_{s}(0,\cdot)$ lying on the plane tangent to $B$ at
$y^{\star}$, we find by expanding $K_{1}(\epsilon S, t)$ about
$\epsilon = 0$ that
\begin{equation}
  K_{1}(\epsilon S, t) = -\frac{\epsilon^{-3} \ell }{S^{3} |
    y_{s}(0,\cdot)|^{3}} + O(\epsilon^{-2}).
\end{equation}
Since the leading-order behavior of $K_{1}$ is independent of $t$, we
use \eqref{eq:mu_ss}, plus knowing that $J(\epsilon S,t) = J(0,\cdot) + O(\epsilon)$ and
$\sin(\epsilon S) = \epsilon S + O(\epsilon^{3})$ to obtain, after expanding about $\delta = 0$
\begin{align}
  \begin{split}
    V^{\text{out}}(y^{\star};\epsilon,\delta) &= \frac{1}{4\pi}
    \int_{0}^{\pi} \int_{0}^{\delta/\epsilon} \left[ - \frac{\epsilon
        \ell J(0,\cdot)}{|y_{s}(0,\cdot)|^{3}}
      \tilde{\mu}_{ss}(0,\cdot) + O(\epsilon^{2}) \right] \mathrm{d}S
    \mathrm{d}t\\
    &= - \frac{\delta \ell J(0,\cdot)}{8 |y_{s}(0,\cdot)|^{3}}
    \varDelta_{S^{2}} \mu(y^{\star}) + O(\epsilon).
  \end{split}
    \label{eq:Vouter-3D}
\end{align}
Note that we have used \eqref{eq:spherical-laplacian} in the last
step. Substituting this result into \eqref{eq:Uouter3D-expand}, we
find that
\begin{multline}
  U^{\text{out}}(y^{\star};\epsilon,\delta) = \frac{1}{4\pi}
  \int_{-\pi}^{\pi} \int_{0}^{\pi} K_{1}(s,t) \left[ \tilde{\mu}(s,t)
    - \tilde{\mu}(0,\cdot)
  \right] J(s,t)\sin(s) \mathrm{d}s \mathrm{d}t  + \delta \frac{\ell J(0,\cdot)}{8 |y_{s}(0,\cdot)|^{3}}
  \varDelta_{S^{2}} \mu(y^{\star}) + O(\epsilon).
  \label{eq:Uout3D-leadingorder}
\end{multline}
This result gives the leading-order asymptotic behavior of
$U^{\text{out}}$.

\subsection{Three-dimensional asymptotic approximation}

We obtain an asymptotic approximation for $U$ by summing the
leading-order behaviors obtained for $U^{\text{in}}$ and
$U^{\text{out}}$ given in \eqref{eq:Uinner3D} and
\eqref{eq:Uout3D-leadingorder}, respectively. Substituting that result
into \eqref{eq:constantapprox}, we obtain the following asymptotic
approximation
\begin{equation}
  u(y^{\star} - \epsilon \ell \nu^{\star}) = f(y^{\star}) + \epsilon
  L_{1}[\mu] + O(\epsilon^{2}),
  \label{eq:asymptotic3D}
\end{equation}
with
\begin{equation}
  L_{1}[\mu] = \frac{1}{4\pi} \int_{-\pi}^{\pi} \int_{0}^{\pi}
  K_{1}(s,t) \left[ \tilde{\mu}(s,t) - \tilde{\mu}(0,\cdot)
  \right] J(s,t)\sin(s) \mathrm{d}s \mathrm{d}t.
\end{equation}
The structure of this asymptotic approximation for the close
evaluation of the double-layer potential in three dimensions is
exactly the same as what we found for the two-dimensional case: the leading-order asymptotic approximation is composed of the Dirichlet data and a non-local term coming from the outer expansion. Similarly, high-order asymptotic approximations could be obtained by continuing on to higher order terms in the expansions $U^{\text{in}}$ and $U^{\text{out}}$ after cumbersome calculations. 


\section{Numerical methods}
\label{sec:numerics}

Numerical methods to compute the asymptotic approximations for the
close evaluation of the double-layer potential must be sufficiently
accurate in comparison to  $O(\epsilon)$. Otherwise, the error made by the
numerical method will dominate over the error of the asymptotic
approximation. On the other hand, if the numerical method requires
restrictively high resolution to compute the asymptotic approximation
to sufficient accuracy, the numerical method suffers from the very
issue of the close evaluation problem. In what follows, we
describe numerical methods to compute the asymptotic approximations
derived above at high accuracy with modest resolution requirements.

\subsection{Two dimensions}

Suppose we have parameterized $B$ by $y = y(\varphi)$ with
$-\pi \le \varphi \le \pi$ with $y^{\star} = y(\varphi^{\star})$. For
that case, we need to compute
\begin{equation}
  \mathscr{U}_{1}(y^{\star}) = \frac{1}{2\pi} \int_{-\pi}^{\pi}
  F_{1}(\varphi;\varphi^{\star}) \mathrm{d}\varphi,
  \label{eq:U1-2D}
\end{equation}
with
\begin{equation}
  F_{1}(\varphi;\varphi^{\star}) = K_{1}(\varphi)
  J(\varphi) \left[ \tilde{\mu}(\varphi) -
    \tilde{\mu}(\varphi^{\star}) \right],
  \label{eq:F-2D}
\end{equation}
where $K_{1}$ is given in \eqref{eq:K1-2D}.  The function, $K_{1}$, is
singular at $\varphi = \varphi^{\star}$. Consequently, applying a high
order accurate numerical quadrature rule to compute $\mathscr{U}_{1}$
will be limited in its accuracy even though $F_{1}$ vanishes
identically at $\varphi = \varphi^{\star}$ due to the factor of
$\tilde{\mu}(\varphi) - \tilde{\mu}(\varphi^{\star})$.  To improve the
accuracy of a numerical evaluation of \eqref{eq:U1-2D}, we revisit the
asymptotic expansion obtained for
$V^{\text{out}}(y^{\star};\epsilon,\delta)$ in \eqref{eq:Vout2D}-\eqref{eq:Vout2Dbis}. By
rewriting that result for the present context, we find
\begin{equation}
  \frac{1}{2\pi} \int_{-\delta/2}^{\delta/2}
  F_{1}(\varphi;\varphi^{\star}) \mathrm{d}\varphi = -\frac{\delta}{4
    \pi} \frac{\ell \tilde{\mu}''(\varphi^{\star}) }{J(\varphi^{\star})}
  + O(\epsilon).
  \label{eq:5.3}
\end{equation}

This result suggests the following method to compute
$\mathscr{U}_{1}(y^{\star})$ numerically using the $N$-point periodic
trapezoid rule (PTR). Suppose we are given the grid function,
$\tilde{\mu}(\varphi_{j})$ for $j = 1, \cdots, N$ with
$\varphi_{j} = -\pi + 2 \pi (j - 1) /N$, and suppose
$\varphi^{\star} = \varphi_{k}$ is one of the quadrature points. We
introduce the numerical approximation
\begin{equation}
  \mathscr{U}_{1}(y^{\star}) \approx U_{1}^{N}(y^{\star}) =
  \frac{1}{N} \sum_{j \neq k} F_{1}(\varphi_{j};\varphi_{k}) -
  \frac{\ell \tilde{\mu}''(\varphi_{k})}{2 N J(\varphi_{k})}.
\end{equation}
where we have replaced the quadrature around $\varphi^{\star}$ with 
\eqref{eq:5.3} where $\delta = 2 \pi/N$. 
We compute $\tilde{\mu}''(\varphi_{k})$ with spectral accuracy using
Fast Fourier transform methods.  Using this numerical approximation,
we compute the $O(\epsilon^{2})$ asymptotic approximation for the
close evaluation of the double-layer potential in two dimensions
through evaluation of
\begin{equation}
  u(y^{\star} - \epsilon \ell \nu^{\star}) \approx f(y^{\star})
  + \epsilon U_{1}^{N}(y^{\star}) + O(\epsilon^{2}).
  \label{eq:numerical-method2D}
\end{equation}

To compute the $O(\epsilon^{3})$ asymptotic approximation, in addition
to $\mathscr{U}_{1}$, we need to compute
\begin{equation}
  \mathscr{U}_{2}(y^{\star}) = \frac{1}{2\pi} \int_{-\pi}^{\pi}
  F_{2}(\varphi;\varphi^{\star}) \mathrm{d}\varphi
\end{equation}
where
\begin{equation}
  F_{2}(\varphi;\varphi^{\star}) = K_{2}(\varphi) J(\varphi) \left[
    \tilde{\mu}(\varphi) - \tilde{\mu}(\varphi^{\star}) \right],
\end{equation}
and $K_{2}$ given in \eqref{eq:K2-2D}. By using the higher-order
asymptotic expansion for $V^{\text{out}}$ (computed on the {\it
  Mathematica} notebook available on the GitHub
repository~\cite{CKK-2018Codes}), we apply the same method used for
$\mathscr{U}_{1}(y^{\star})$ and arrive at
\begin{equation}
  \mathscr{U}_{2}(y^{\star}) \approx U_{2}^{N} = \frac{1}{N} \sum_{j
    \neq k} F_{2}(\varphi_{j};\varphi_{k}) - \frac{\ell^{2}
    \kappa^{\star} \tilde{\mu}''(\varphi_k)}{4 N J(\varphi_k)},
\end{equation}
with $\kappa^{\star}$ denoting the signed curvature at
$y^{\star}$. Using this numerical approximation, we compute the
$O(\epsilon^{3})$ asymptotic approximation for the close evaluation of
the double-layer potential in two dimensions through evaluation of
\begin{multline}
  u(y^{\star} - \epsilon \ell \nu^{\star}) \approx f(y^{\star})
  + \epsilon U_{1}^{N}(y^{\star}) + \epsilon^{2} \left[ U_{2}^{N}(y^{\star}) - \frac{\ell^{2}
      y'(\varphi^{\star}) \cdot y''(\varphi^{\star})}{4
      J^{4}(\varphi^{\star})} \tilde{\mu}'(\varphi^{\star}) +
    \frac{\ell^{2} \tilde{\mu}''(\varphi^{\star})}{4
      J^{2}(\varphi^{\star})} \right] + O(\epsilon^{3}).
  \label{eq:numerical-method2D-higher-order}
\end{multline}
Since the boundary is given, we are able to compute $y'(\varphi)$ and
$y''(\varphi)$, explicitly. We use Fast Fourier transform methods to
compute $\tilde{\mu}'(\varphi)$ and $\tilde{\mu}''(\varphi)$ with
spectral accuracy.

\subsection{Three dimensions}

Suppose we have parameterized $B$ by $y = y(\theta,\varphi)$ with
$\theta \in [0,\pi]$ and $\varphi \in [-\pi,\pi]$ with
$y^{\star} = y(\theta^{\star},\varphi^{\star})$. For that case, we
seek to compute
\begin{equation}
  \mathscr{U}_{1}(y^{\star}) = \frac{1}{4\pi} \int_{-\pi}^{\pi}
  \int_{0}^{\pi} F_{1}(\theta,\varphi;\theta^{\star},\varphi^{\star})
  \sin(\theta) \mathrm{d}\theta \mathrm{d}\varphi,
  \label{eq:U1-3D}
\end{equation}
with
\begin{equation}
  F_{1}(\theta,\varphi;\theta^{\star},\varphi^{\star}) =
  K_{1}(\theta,\varphi) J(\theta,\varphi) \left[
    \tilde{\mu}(\theta,\varphi) -
    \tilde{\mu}(\theta^{\star},\varphi^{\star}) \right],
  \label{eq:F-3D}
\end{equation}
where $K_{1}$ is given in \eqref{eq:K1-3D}.  Just as with the
two-dimensional case, the function $K_{1}$ is singular at
$(\theta,\varphi) = (\theta^{\star},\varphi^{\star})$, so any attempt
to apply a quadrature rule to compute $\mathscr{U}_{1}$ will be
limited in its accuracy even though $F_{1}$ vanishes identically at
$(\theta^{\star}, \varphi^{\star})$ due to the factor of
$\tilde{\mu}(\theta,\varphi) -
\tilde{\mu}(\theta^{\star},\varphi^{\star})$.

To numerically evaluate \eqref{eq:U1-3D}, we apply a three-step method
developed by the authors~\cite{carvalho2018_3D}. This method has been
shown to be effective for computing the modified double-layer
potential in three dimensions resulting from the subtraction
method. We first rotate this integral to another spherical coordinate
system in which $y^{\star}$ is aligned with the north pole. The
details of this rotation are given in Appendix \ref{sec:rotation} and
lead to $\theta = \theta(s,t)$ and $\varphi = \varphi(s,t)$ with
$s \in [0,\pi]$ and $t \in [-\pi,\pi]$ where
$\theta^{\star} = \theta(0,\cdot)$ and
$\varphi^{\star} = \varphi(0,\cdot)$. We apply this rotation and find
that
\begin{equation}
  \mathscr{U}_{1}(y^{\star}) = \frac{1}{4\pi} \int_{-\pi}^{\pi}
  \int_{0}^{\pi} \tilde{F}(s,t) \sin(s) \mathrm{d}s \mathrm{d}t,
  \label{eq:5.10}
\end{equation}
with
$\tilde{F}(s,t) = F_{1}(\theta(s,t), \varphi(s,t); \theta^{\star},
\varphi^{\star})$. Now $\tilde{K}_{1}(\theta(s,t),\varphi(s,t))$ is
singular at the north pole of this rotated coordinate system
corresponding to $s = 0$.  To improve the accuracy of a numerical
evaluation of \eqref{eq:5.10}, we revisit the asymptotic expansion
obtained for $V^{\text{out}}(y^{\star};\epsilon,\delta)$ in
\eqref{eq:Vouter-3D}. By rewriting that result for the present
context, we find
\begin{equation}
  \frac{1}{2} \int_{0}^{\delta}  \left[ \frac{1}{2\pi}
    \int_{-\pi}^{\pi} \tilde{F}(s,t) \mathrm{d}t \right] \sin(s)
  \mathrm{d}s = \int_{0}^{\delta} \left[ - \frac{\ell J(0,\cdot)}{8 |
      y_{s}(0,\cdot) |^{3}} \varDelta_{S^{2}} \mu(y^{\star})  \right]
  \mathrm{d}s + O(\epsilon).
    \label{eq:L1-limit}
\end{equation}
Suppose we compute
\begin{equation}
  \bar{F}(s) = \frac{1}{2\pi} \int_{-\pi}^{\pi} \tilde{F}(s,t)
  \mathrm{d}t.
\end{equation}
The result in \eqref{eq:L1-limit} suggests that $\bar{F}(s)$ smoothly
limits to a finite value as $s \to 0$.  Although we could use this
result to evaluate $\bar{F}(s)$ in a numerical quadrature scheme, it
will suffice to consider an open quadrature rule for $s$ that does not
include the point $s = 0$ such as the Gauss-Legendre quadrature.  This
result suggests the following three-step method to compute
$\mathscr{U}_{1}(y^{\star})$ numerically.

Let $t_{k} = -\pi + \pi (k - 1)/N$ for
$k = 1, \cdots, 2N$, and let $z_{j}$ and $w_{j}$ for
$j = 1, \cdots, N$ denote the $N$-point Gauss-Legendre quadrature
abscissas and weights such that
\begin{equation}
  \int_{-1}^{1} f(x) \mathrm{d}x \approx \sum_{j = 1}^{N} f(z_{j})
  w_{j}.
\end{equation}
We perform the mapping: $s_{j} = \pi ( z_{j} + 1 )/2$ for
$j = 1, \cdots, N$, and make appropriate adjustments to the weights as
will be shown below. For the first step, we rotate the spherical
coordinate system so that $y^{\star}$ is aligned with its north pole
as described in Appendix \ref{sec:rotation}.  For the second step, we
compute
\begin{equation}
  \bar{F}(s_{j}) \approx \bar{F}^{N}_{j} = \frac{1}{2N} \sum_{k =
    1}^{2N} \tilde{F}(s_{j},t_{k}), \quad j = 1, \cdots, N.
\end{equation}
For the third step, we compute the numerical approximation
\begin{equation}
  \mathscr{U}_{1}(y^{\star}) \approx U_{1}^{N}(y^{\star}) =
  \frac{\pi}{4} \sum_{j = 1}^{N} \bar{F}_{j}^{N} w_{j}.
  \label{eq:L1-3D-numerics}
\end{equation}
In \eqref{eq:L1-3D-numerics}, a factor of $\pi/2$ is introduced to
scale the quadrature weights due to the mapping from $z_{j}$ to
$s_{j}$, and a factor of $1/2$ remains from the factor of $1/4\pi$ in
\eqref{eq:U1-3D}. 

Using the numerical approximation $U_{1}^{N}$, we compute the
$O(\epsilon^{2})$ asymptotic approximation for the close evaluation of
the double-layer potential in three dimensions through evaluation of
\begin{equation}
  u(y^{\star} - \epsilon \ell \nu^{\star}) \approx f(y^{\star})
  + \epsilon U_{1}^{N}(y^{\star}) + O(\epsilon^{2}).
  \label{eq:numerical-method3D}
\end{equation}

\section{Numerical results}
\label{sec:results}

We present results that show the accuracy and efficiency of
the asymptotic approximation and the corresponding numerical method
 for the close evaluation of the
double-layer potential. For all of the examples shown, we
prescribe Dirichlet data corresponding to a particular harmonic
function. With that Dirichlet data, we solve the boundary integral
equation \eqref{eq:bie} numerically to obtain the density, $\mu$. We
use that density to compute the double-layer
potential using different methods, for comparison. The results below
show the error made in computing the harmonic function at close
evaluation points. The Matlab codes used to compute all of the
following examples are available in a GitHub
repository~\cite{CKK-2018Codes}.

\subsection{Two dimensions}

For the two-dimensional examples, we use the harmonic
function,
\begin{equation}
  u(x) = -\frac{1}{2\pi} \log | x - x_{0} |,
\end{equation}
with $x_{0} \in \mathbb{R}^2 \setminus \bar{D}$
and prescribe Dirichlet data by evaluating
this function on the boundary. We solve the boundary integral equation
\eqref{eq:bie} using the Nystr\"om method with the $N$-point Periodic
Trapezoid Rule (PTR) resulting in the numerical approximation for the density,
$\tilde{\mu}_{j} \approx \tilde{\mu}(\varphi_{j})$ with
$\varphi_{j} = -\pi + 2 (j-1) \pi/N$ for $j = 1, \cdots, N$.

We compute the close evaluation of the double-layer potential at
points, $x = y^{\star} - \epsilon \nu^{\star}$, using the following
four methods:
\begin{enumerate}

\item {\bf PTR method} -- Compute the double-layer potential,
  \begin{equation*}
    u(y^{\star} - \epsilon \nu^{\star}) = \frac{1}{2\pi}
    \int_{-\pi}^{\pi} \frac{\nu_{y} \cdot ( y^{\star} - \epsilon
      \nu^{\star} - y)}{| y^{\star} - \epsilon \nu^{\star} - y|^{2}}
    \mu(y) \mathrm{d}\sigma_{y},
  \end{equation*}
  using the same $N$-point PTR used to solve \eqref{eq:bie} .

\item {\bf Subtraction method} -- Compute the modified
  double-layer potential,
  \begin{equation*}
    u(y^{\star} - \epsilon \nu^{\star}) = -\mu(y^{\star}) + \frac{1}{2\pi}
    \int_{-\pi}^{\pi} \frac{\nu_{y} \cdot ( y^{\star} - \epsilon
      \nu^{\star} - y)}{| y^{\star} - \epsilon \nu^{\star} - y|^{2}}
    \left[ \mu(y) - \mu(y^{\star}) \right] \mathrm{d}\sigma_{y},
  \end{equation*}
  using the same $N$-point PTR used to solve \eqref{eq:bie} and as in the first method.

\item {\bf $O(\epsilon^{2})$ asymptotic approximation} -- Compute the
  $O(\epsilon^{2})$ asymptotic approximation given by
  \eqref{eq:asymptotic2D} using the new numerical method given in
  \eqref{eq:numerical-method2D} using the same $N$-point PTR used to
  solve \eqref{eq:bie}.
  
\item {\bf $O(\epsilon^{3})$ asymptotic approximation} -- Compute the
  $O(\epsilon^{3})$ asymptotic approximation given by
  \eqref{eq:asymptotic2D-higherorder} using the new numerical method given
  in \eqref{eq:numerical-method2D-higher-order} using the same
  $N$-point PTR used to solve \eqref{eq:bie}.
  
\end{enumerate}

We consider two different domains: 
 \begin{itemize}
 \item {\bf A kite domain} whose boundary is
given by
\begin{equation}\label{eq:kite}
  y(t) = ( \cos t + 0.65 \cos 2t - 0.65, 1.5 \sin t), \quad -\pi \le t
  \le \pi.
\end{equation}
 \item {\bf A star domain} whose boundary is
given by
\begin{equation}\label{eq:star}
  y(t) = r(t) ( \cos t, \sin t ), \quad r(t) = 1 + 0.3 \cos 5 t, \quad
  -\pi \le t \le \pi.
\end{equation}
 \end{itemize} 
For both examples we pick $x_0 = (1.85, 1.65)$ which lies outside the domains. 
We consider $N$ fixed, here $N = 128$, and study the dependence of the error on
$\epsilon$ as $\epsilon \to 0$.

\begin{figure}[!h]
  \centering
  \includegraphics[width=0.44\linewidth]{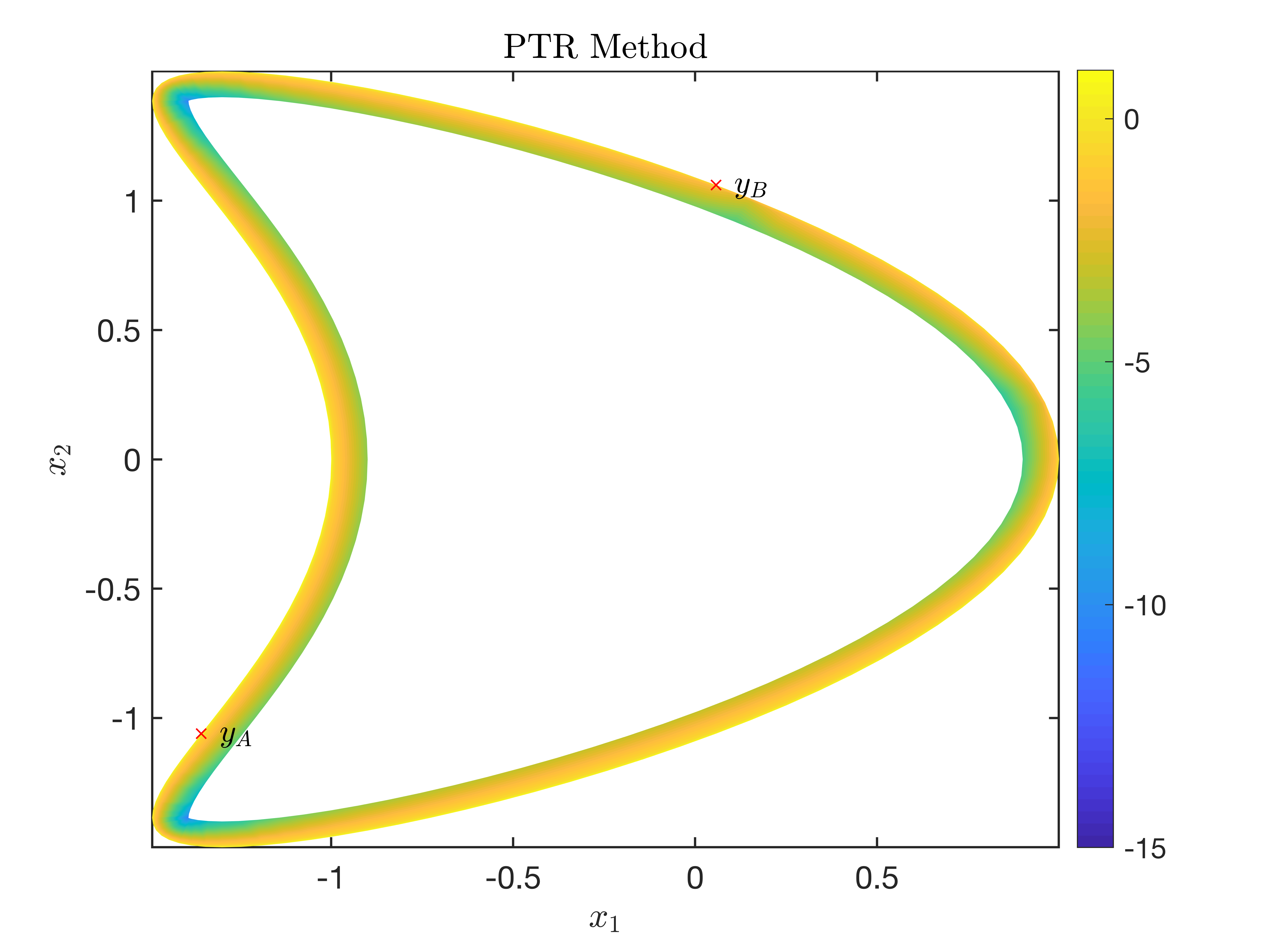}
  \includegraphics[width=0.44\linewidth]{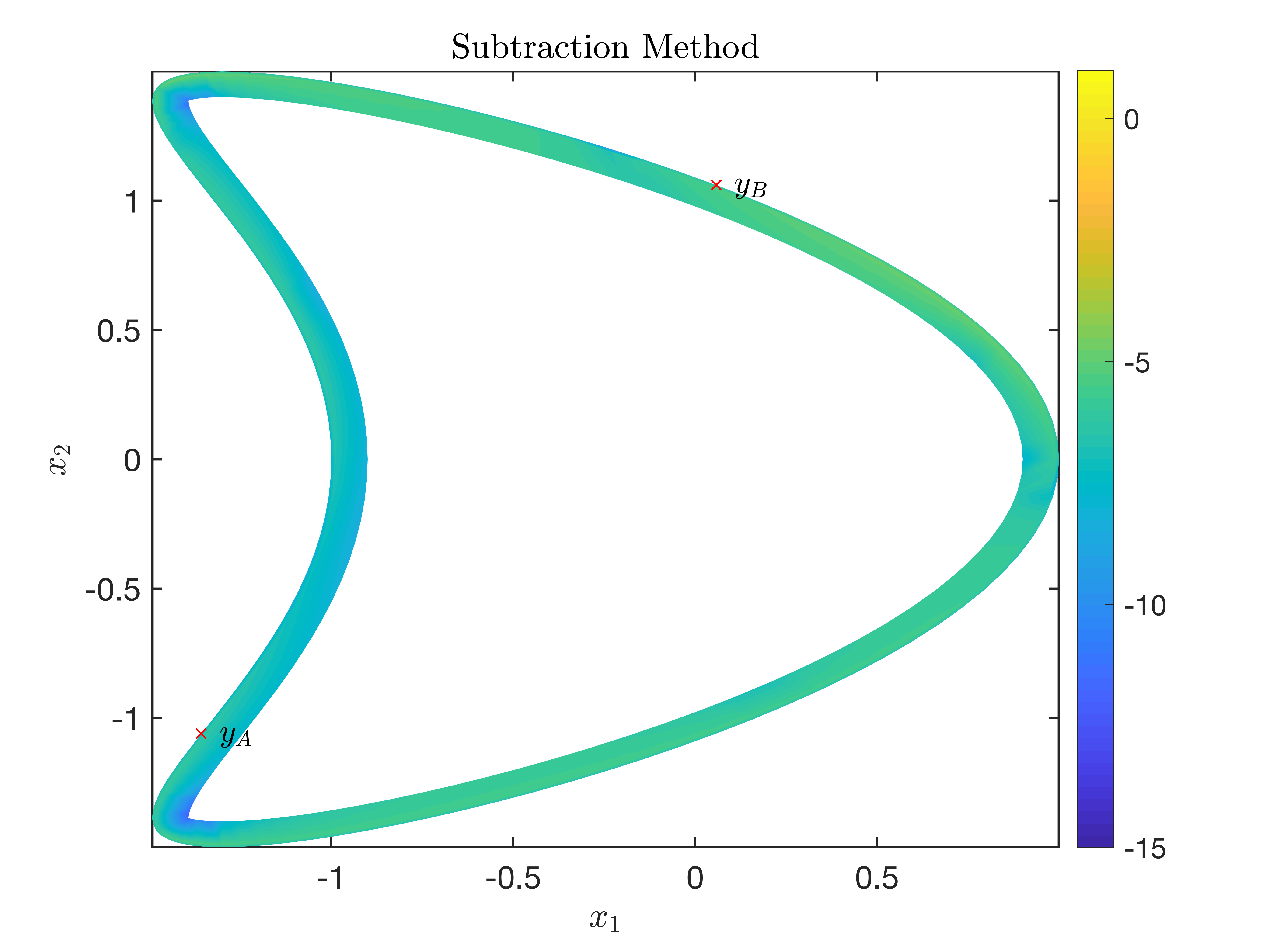}\\
  \includegraphics[width=0.44\linewidth]{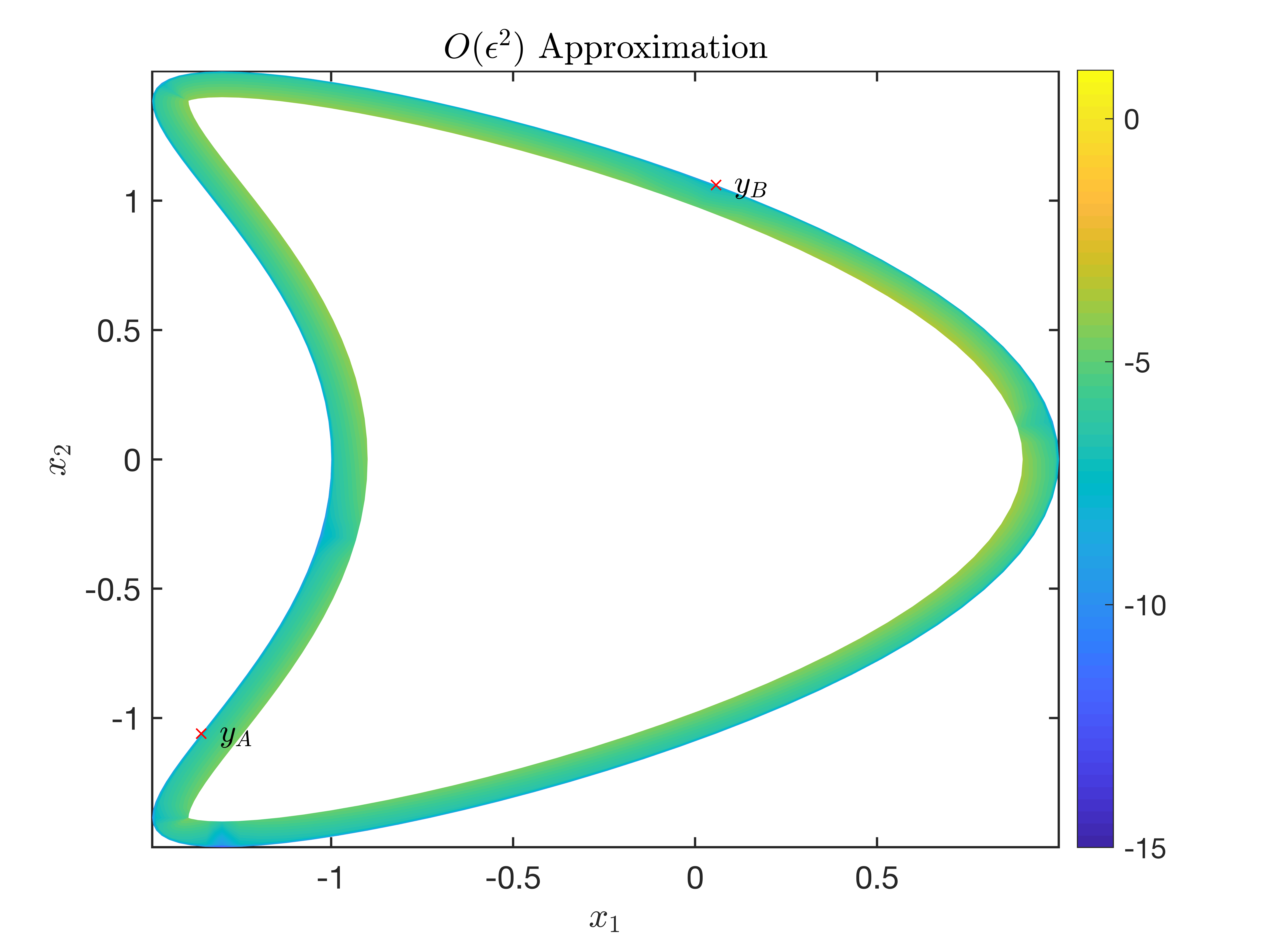}
  \includegraphics[width=0.44\linewidth]{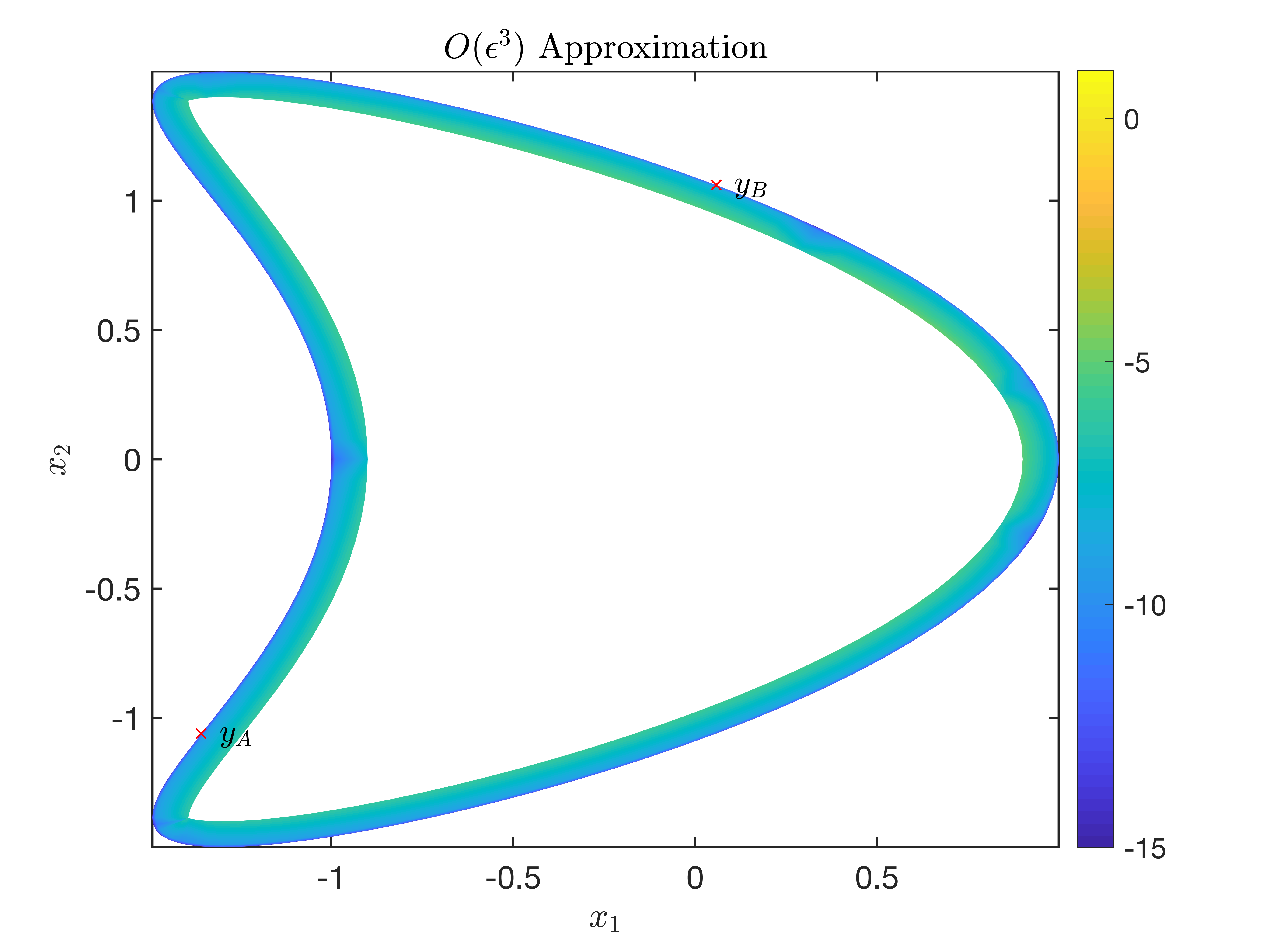}
  \caption{Plots of $\log_{10}$ of the error for the evaluation of the double-layer potential in the kite domain defined by \eqref{eq:kite} using four methods:
    the PTR method (upper left), the subtraction method (upper right), the $O(\epsilon^{2})$ asymptotic
    approximation method (lower left), and 
    the $O(\epsilon^{3})$ asymptotic approximation method (lower right). In each of these
    plots, boundary points $y_{A} = (-1.3571, -1.0607)$ and
    $y_{B} = (0.0571, 1.0607)$ are plotted as red $\times$'s.}
  \label{kite-results}

  \bigskip

  \centering
  \includegraphics[width=0.48\linewidth]{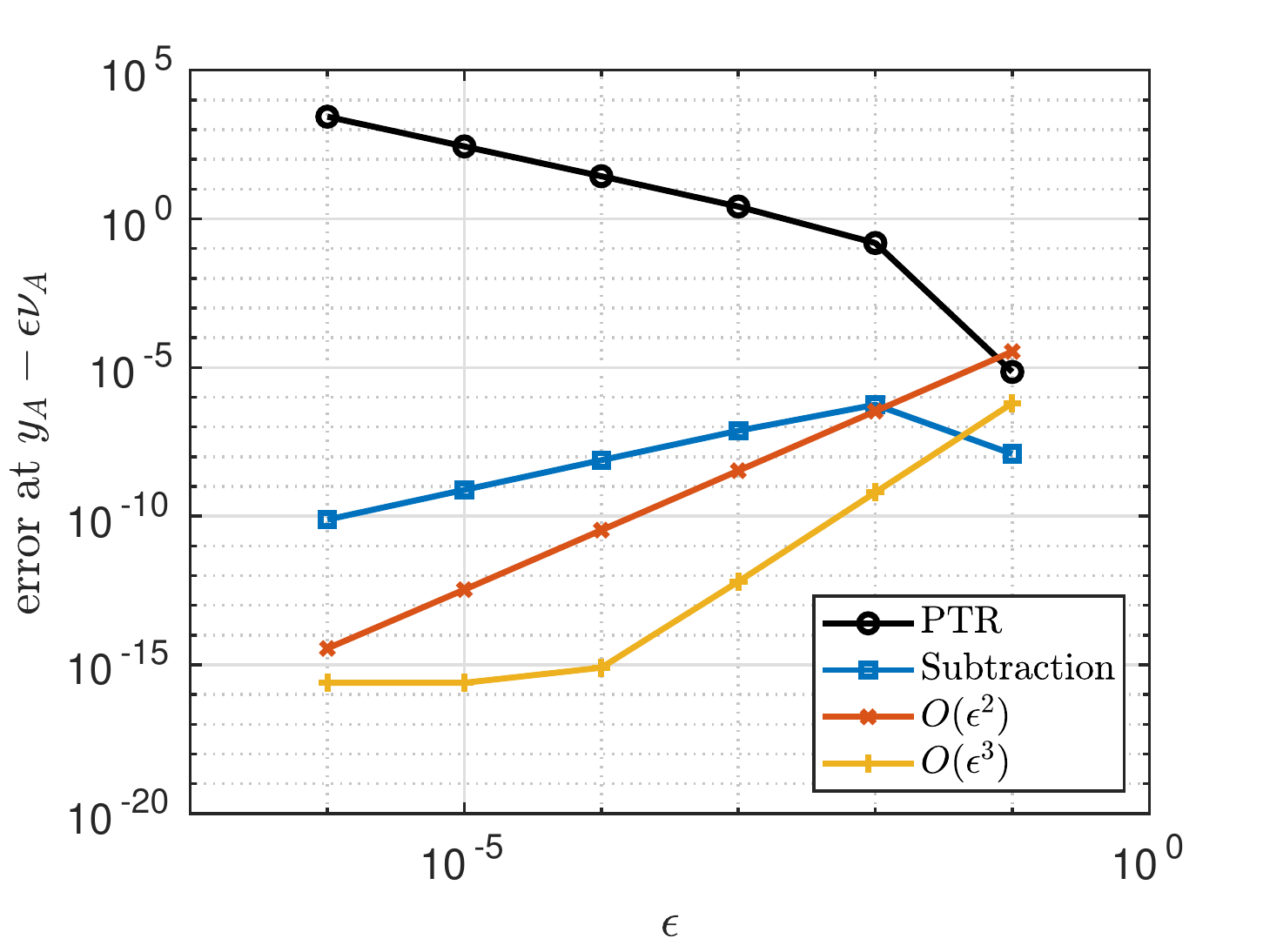}
  \includegraphics[width=0.48\linewidth]{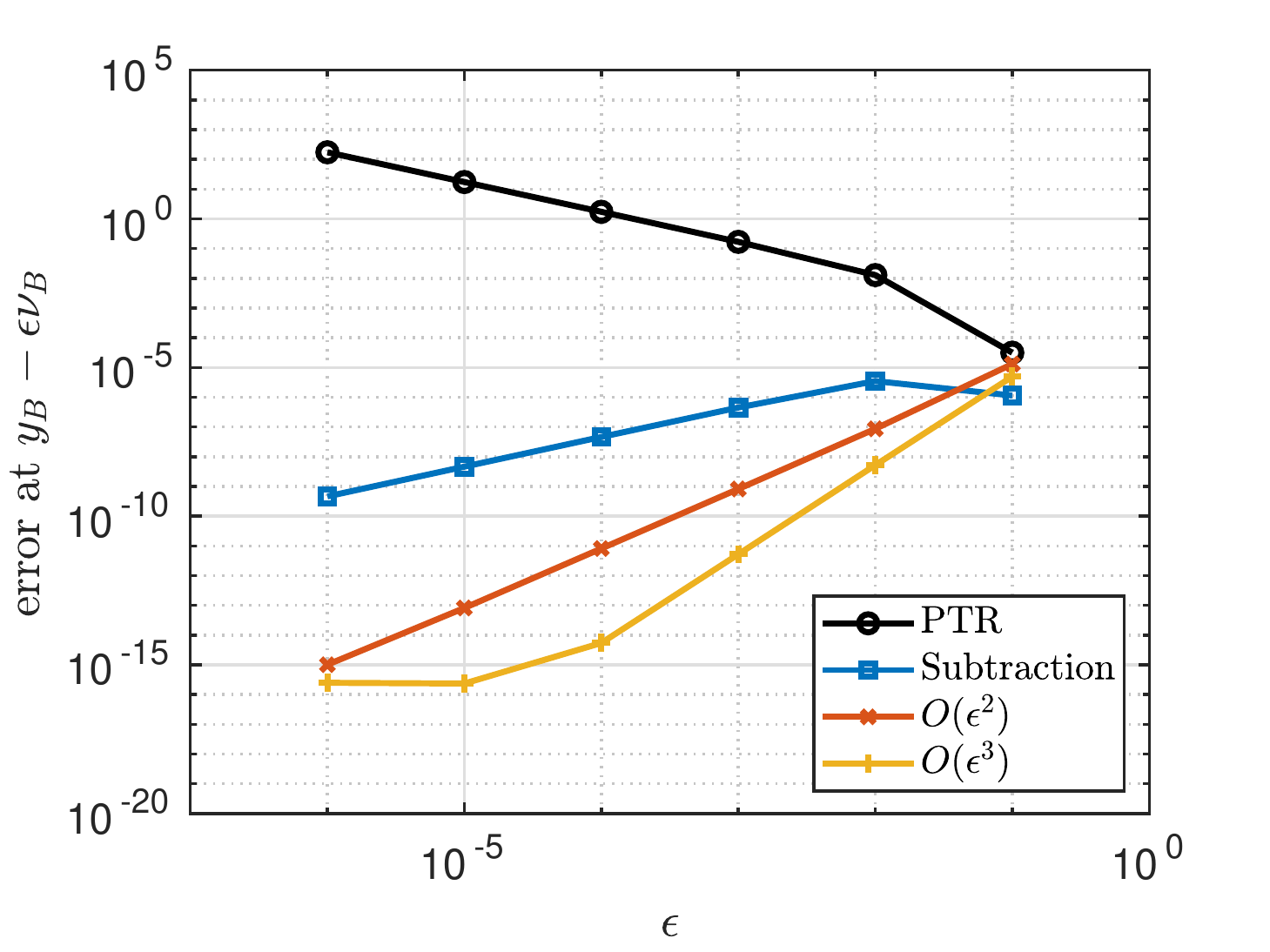}
  \caption{Log-log plots of the errors made in computing the double-layer potential by the four methods shown
    in Fig.~\ref{kite-results} at $y_{A} - \epsilon \nu_{A}$ (left)
    and at $y_{B} - \epsilon \nu_{B}$ (right) for
    $10^{-6} \le \epsilon \le 10^{-1}$. }    
  \label{kite-errors}
\end{figure}

In Fig.~\ref{kite-results} we show results for the kite domain. The error, using a $\log$ scale, is presented for 
each of the four methods described above. The results show that the PTR method exhibits an $O(1)$ error
as $\epsilon \to 0$. The subtraction method and the asymptotic approximations all show
substantially smaller errors. 


\begin{figure}[htb]
  \centering
  \includegraphics[width=0.48\linewidth]{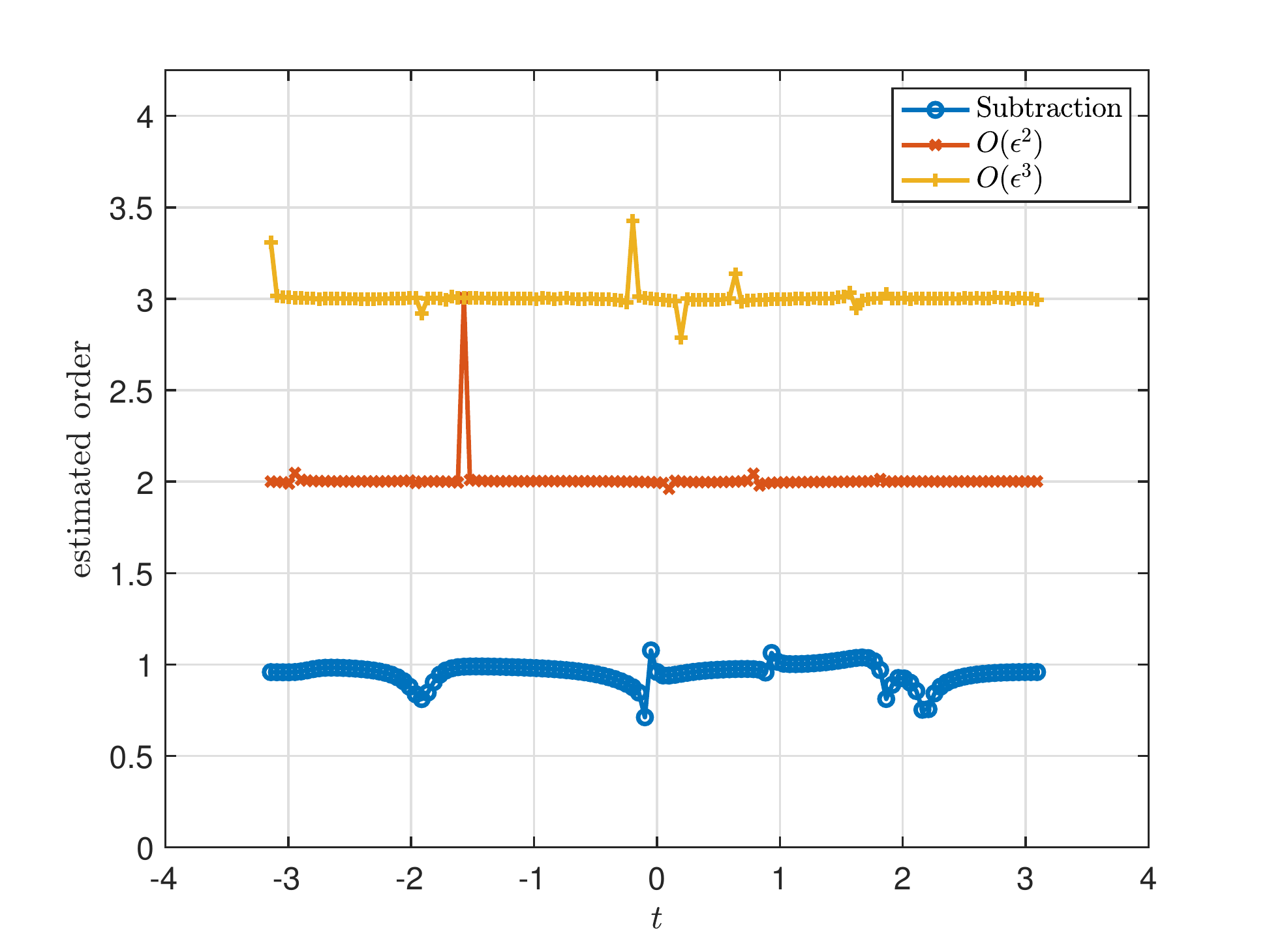}
  \caption{Estimated order of accuracy in computing the double-layer
    potential in the kite domain when using the subtraction method
    (blue $\circ$), the $O(\epsilon^{2})$ asymptotic approximation method
    (red $\times$), and the $O(\epsilon^{3})$ asymptotic approximation method
    (yellow $+$) for varying values of $t$.}
  \label{kite-order}
\end{figure}

To compare the four methods more quantitatively, in
Fig.~\ref{kite-errors} we plot the errors made by the four methods at
$y_{A} - \epsilon \nu_{A}$ (left) and $y_{B} - \epsilon \nu_{B}$
(right) with $10^{-6} \le \epsilon \le 10^{-1}$ where $\nu_{A}$ and
$\nu_{B}$ are the unit outward normals at $y_{A}=(-1.3571, -1.0607)$,
and $y_{B}=(0.0571, 1.0607)$, respectively.  The points $y_A$ and
$y_B$ are shown in each plot of Fig. \ref{kite-results}.  From the
results in Fig.~\ref{kite-errors} we observe that, the error when
using the PTR method increases as $\epsilon \to 0^{+}$, while the
error in the other three methods decreases.  The errors made by the
asymptotic approximations are monotonically decreasing as
$\epsilon \to 0^{+}$. However, the error made by the subtraction
method presents a different behavior: it reaches a maximum at
$\epsilon \approx 10^{-2}$ after which it decreases as $\epsilon$
increases. For larger values of $\epsilon$, the double-layer potential
is no longer nearly singular, so the $N$-point PTR (and therefore
methods 1 and 2) become more accurate. The error is at a maximum for
the subtraction method when $\epsilon = O(1/N)$, which is why we
observe the maximum error occuring at $\epsilon \approx 10^{-2}$.  The
results in Fig.~\ref{kite-errors} show a clear difference in the rate
at which the errors vanish as $\epsilon \to 0^{+}$ between the
subtraction method and the asymptotic approximation methods. The
$O(\epsilon^{3})$ asymptotic approximation decays the fastest,
followed by the $O(\epsilon^{2})$ asymptotic approximation, and then
the subtraction method.  For $\epsilon < 10^{-4}$, the error incurred
by the $O(\epsilon^{3})$ asymptotic approximation levels out at
machine precision.  We estimate the rate at which the subtraction
method and the asymptotic approximation methods decay with respect to
$\epsilon$ from the slope of the best fit line through the $\log-\log$
plot of the error versus $\epsilon$ in Fig.~\ref{kite-order}.  We
compute the slope for each evaluation point $y(t_{j})- \epsilon \tilde{\nu}(t_{j})$ where
$t_{j} = -\pi + 2 (j-1) \pi/N$ for $j = 1, \cdots, N$, and we vary
$\epsilon$.  For the
subtraction method, we consider $\epsilon$ values such that
$10^{-6} \le \epsilon \le 10^{-2}$, and for the $O(\epsilon^{2})$ and
$O(\epsilon^{3})$ asymptotic approximation methods, we consider the
same range but only include values where the error is greater than
$10^{-15}$. The results shown in Fig.~\ref{kite-order} indicate that
the subtraction method decays linearly with $\epsilon$, and the rates
of the asymptotic approximations are consistent with the theory
presented in Section \ref{sec:asymptotics2D}.

\begin{figure}[h!]
  \centering
  \includegraphics[width=0.44\linewidth]{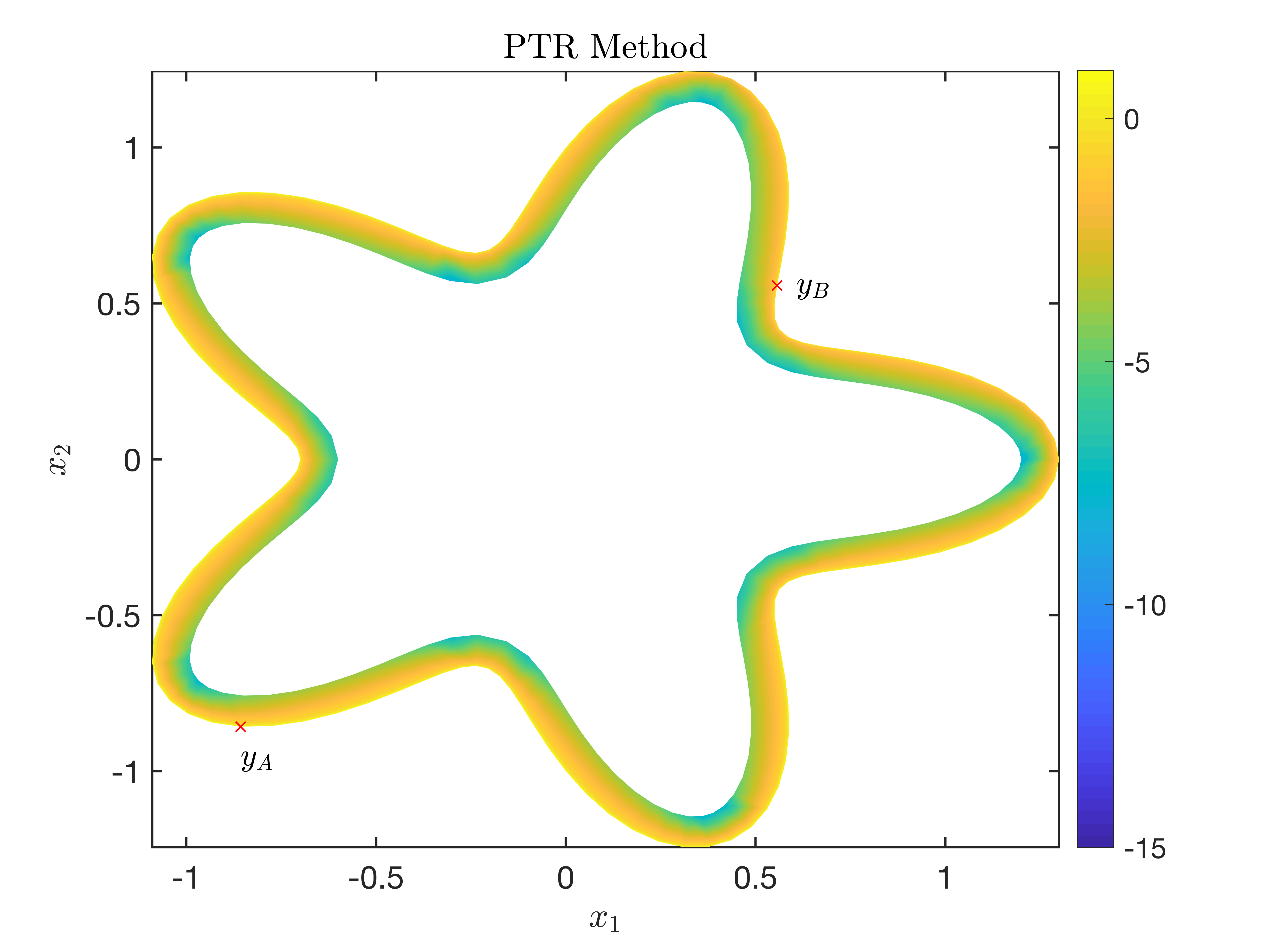}
  \includegraphics[width=0.44\linewidth]{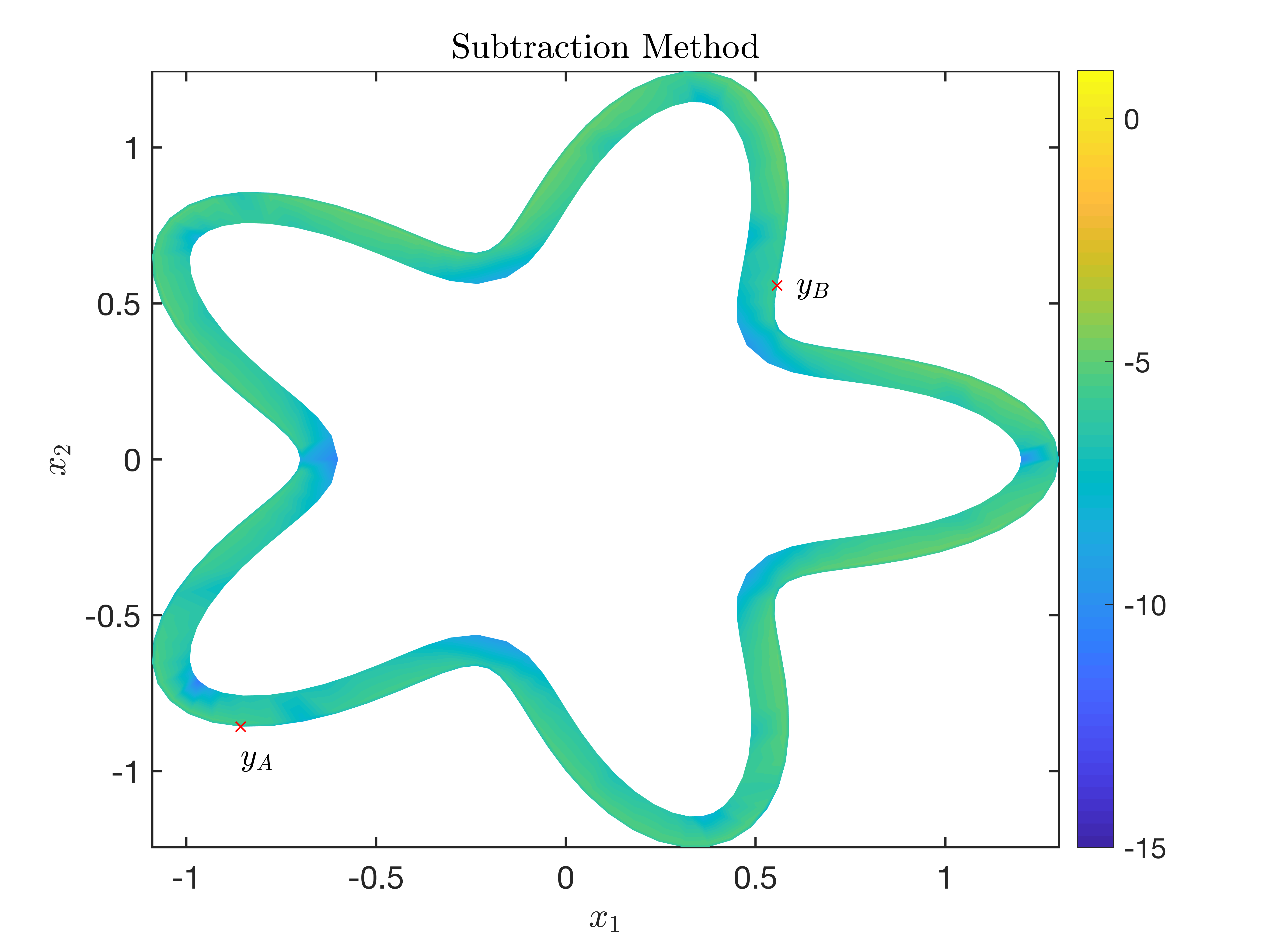}\\
  \includegraphics[width=0.44\linewidth]{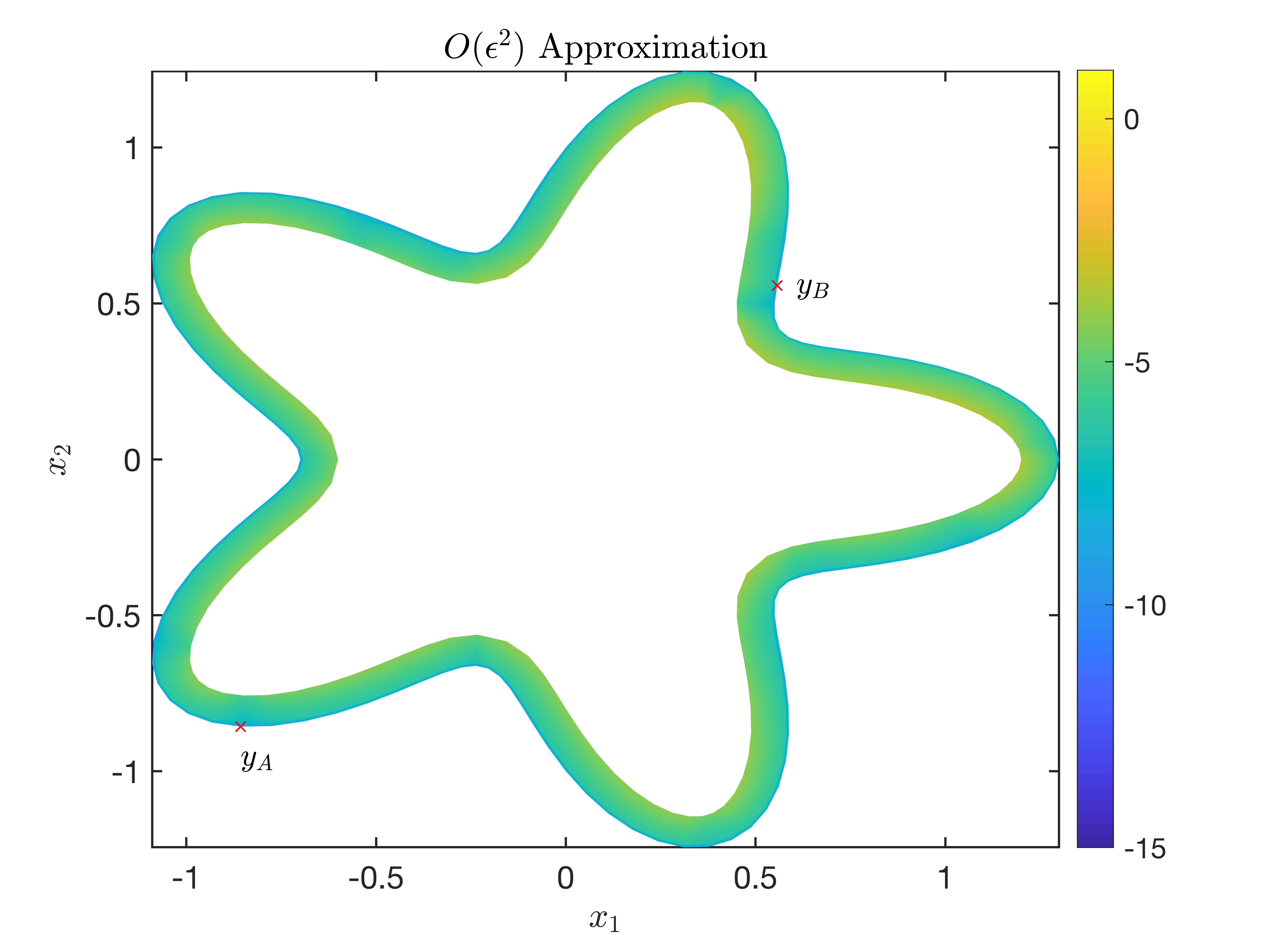}
  \includegraphics[width=0.44\linewidth]{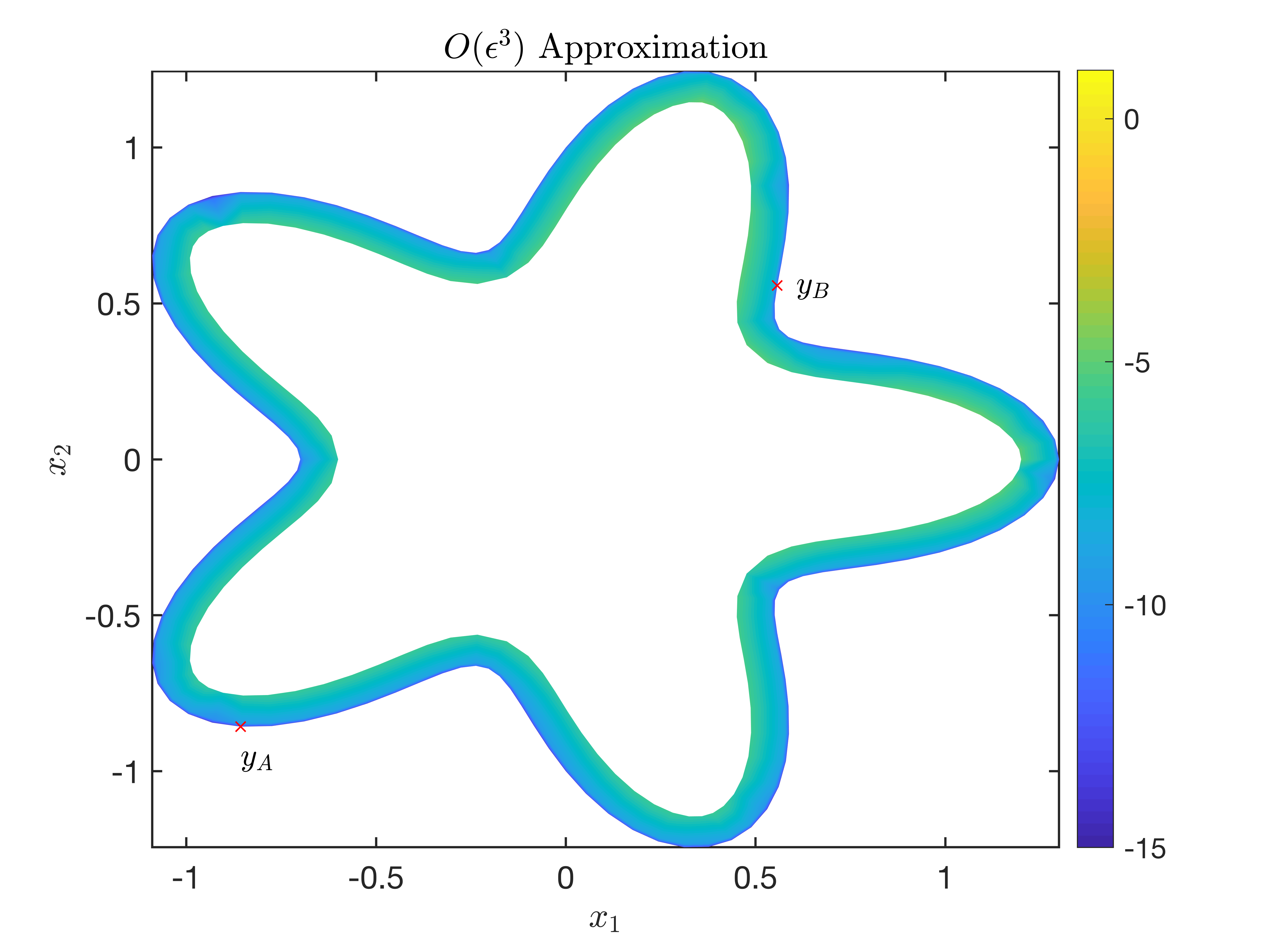}
  \caption{Plots of $\log_{10}$ of the error for the evaluation of the
    double-layer potential in the star domain defined by
    \eqref{eq:star} using four methods: the PTR method (upper left),
    the subtraction method (upper right), the $O(\epsilon^{2})$
    asymptotic approximation method (lower left), and the
    $O(\epsilon^{3})$ asymptotic approximation method (lower
    right). In each of these plots, boundary points
    $y_{A} = (-1.3571, -1.0607)$ and $y_{B} = (0.0571, 1.0607)$ are
    plotted as red $\times$'s.}
  \label{star-results}

  \bigskip

  \centering
  \includegraphics[width=0.48\linewidth]{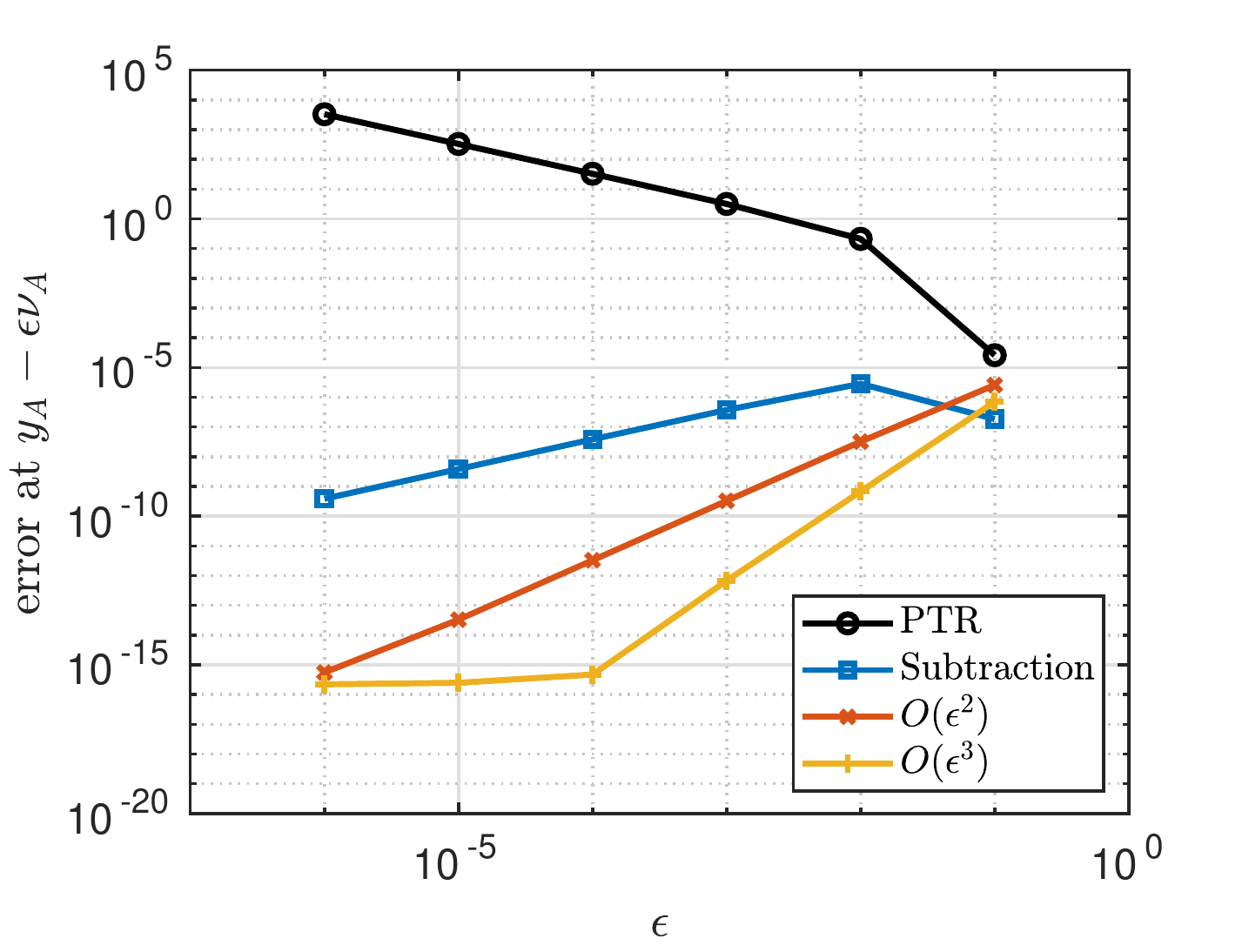}
  \includegraphics[width=0.48\linewidth]{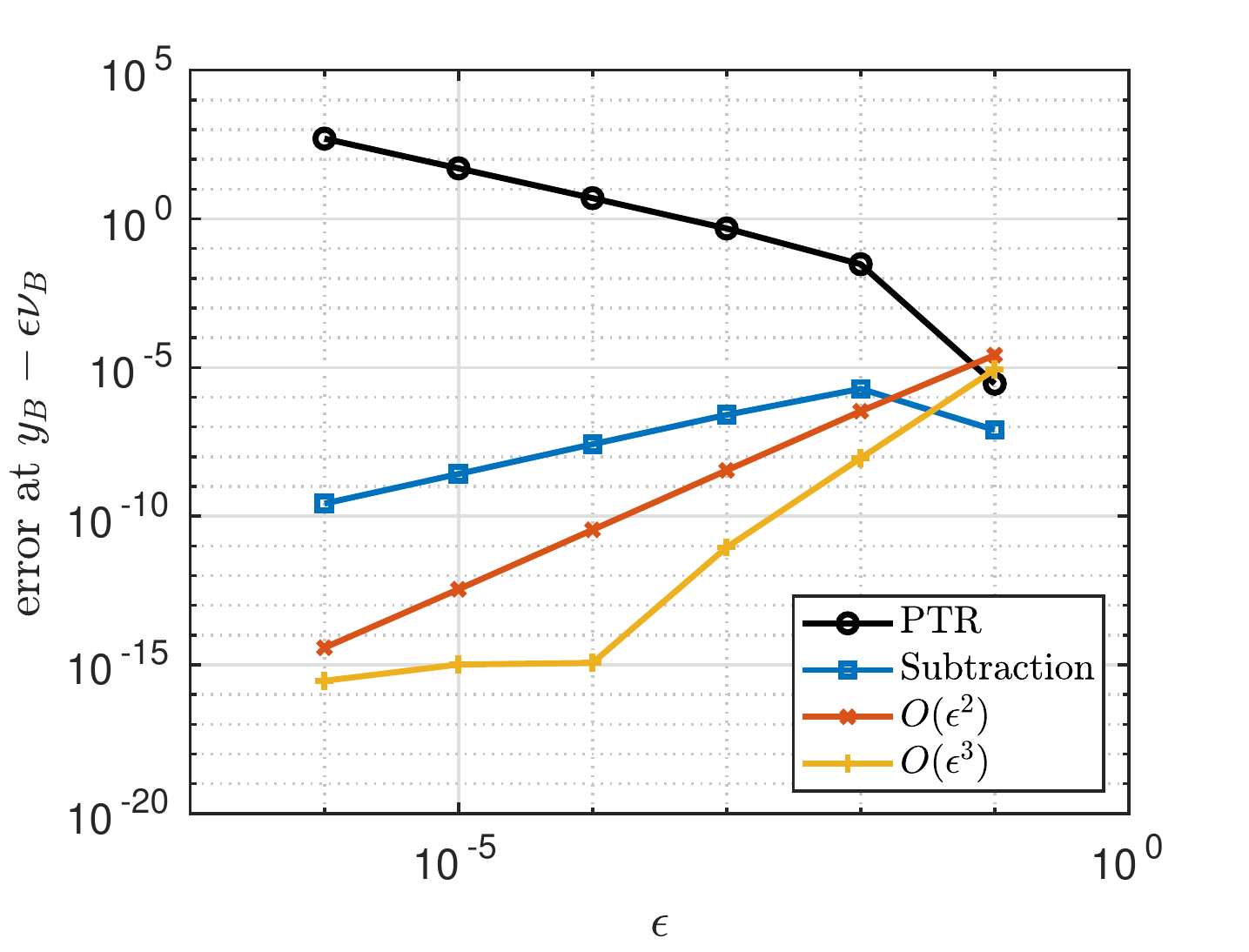}
   \caption{Log-log plots of the errors made in computing the double-layer potential by the four methods shown
    in Fig.~\ref{star-results} at $y_{A} - \epsilon \nu_{A}$ (left)
    and at $y_{B} - \epsilon \nu_{B}$ (right) for
    $10^{-6} \le \epsilon \le 10^{-1}$. }    
  \label{star-errors}

\end{figure}

\begin{figure}[h!]
  \centering
  \includegraphics[width=0.48\linewidth]{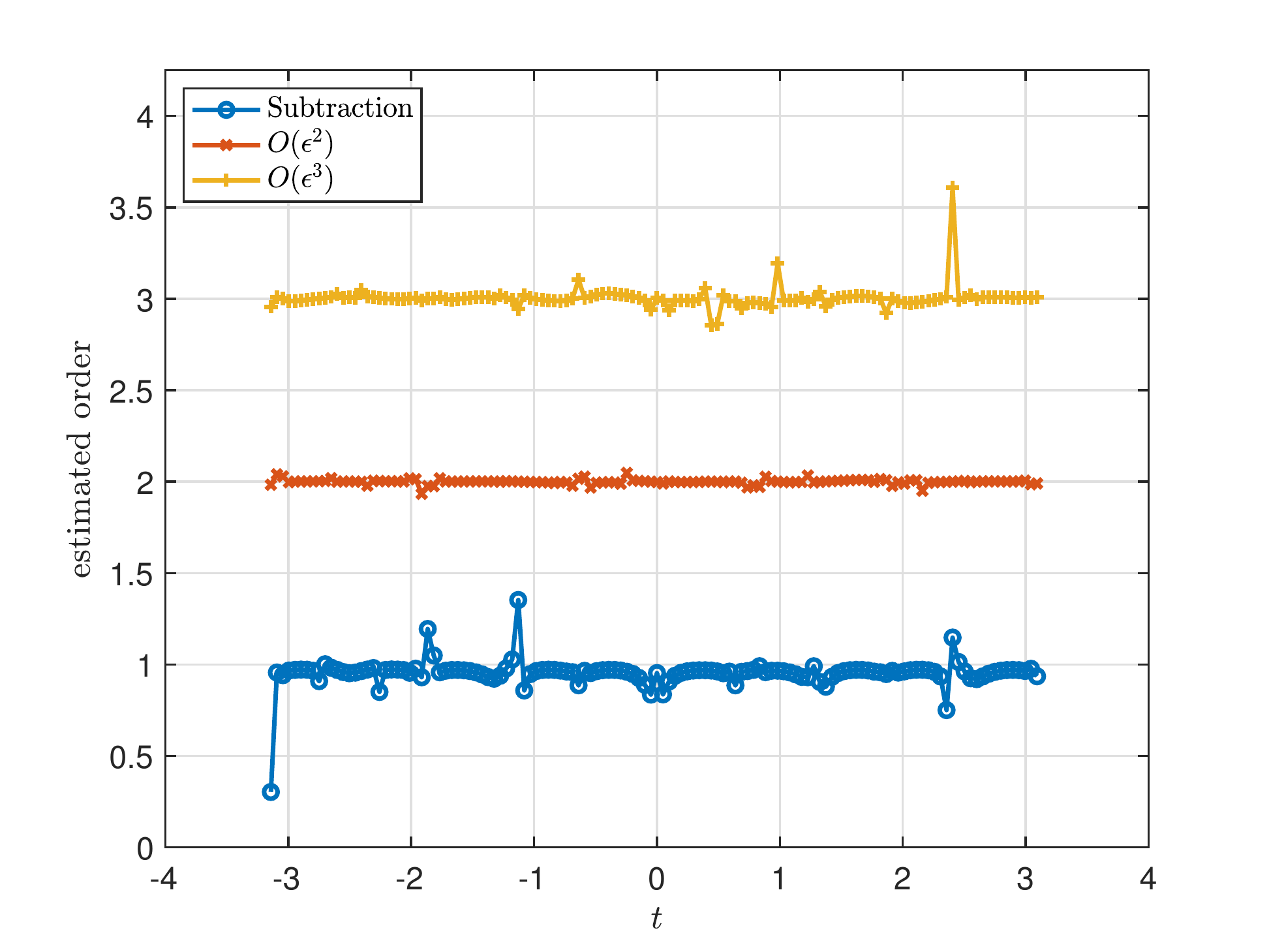}
    \caption{Estimated order of accuracy in computing the double-layer potential in the star domain when using the
    subtraction method (blue $\circ$), the $O(\epsilon^{2})$ asymptotic
    approximation method (red $\times$), and the $O(\epsilon^{3})$ asymptotic
    approximation method (yellow $*$) for varying values of $t$.}
  \label{star-order}
\end{figure}

For the second example of computing the double-layer potential in the star domain, Figures~\ref{star-results}, \ref{star-errors}, and \ref{star-order}
 are analogous to Figures~\ref{kite-results},
\ref{kite-errors}, and \ref{kite-order} for the kite domain. The characteristics of the errors for this second domain
are exactly the same as described for the kite domain. 

\subsubsection*{Summary of the results}

In the case of two dimensional-problems, the subtraction method yields
a method whose error decays linearly with the distance away from the
boundary. The $O(\epsilon^{2})$ and $O(\epsilon^{3})$ asymptotic
approximations methods are much more accurate for close evaluation
points. Moreover, only relatively modest resolution is required for
these asymptotic approximations to be effective. However, the error
for these asymptotic approximations are monotonically increasing with
the distance to the boundary, so they are not accurate for points
further away from the boundary.  The error estimates provided by the
asymptotic theory provide guidance on where to apply these asymptotic
approximations effectively.  For all of these reasons, we find that
the asymptotic approximations and corresponding numerical methods are
quite useful for two-dimensional problems.

\subsection{Three dimensions}

\begin{figure}[h!]
  \begin{minipage}{0.99\textwidth}
    \centering
    \raisebox{-0.5\height}{\includegraphics[height=1.8in]{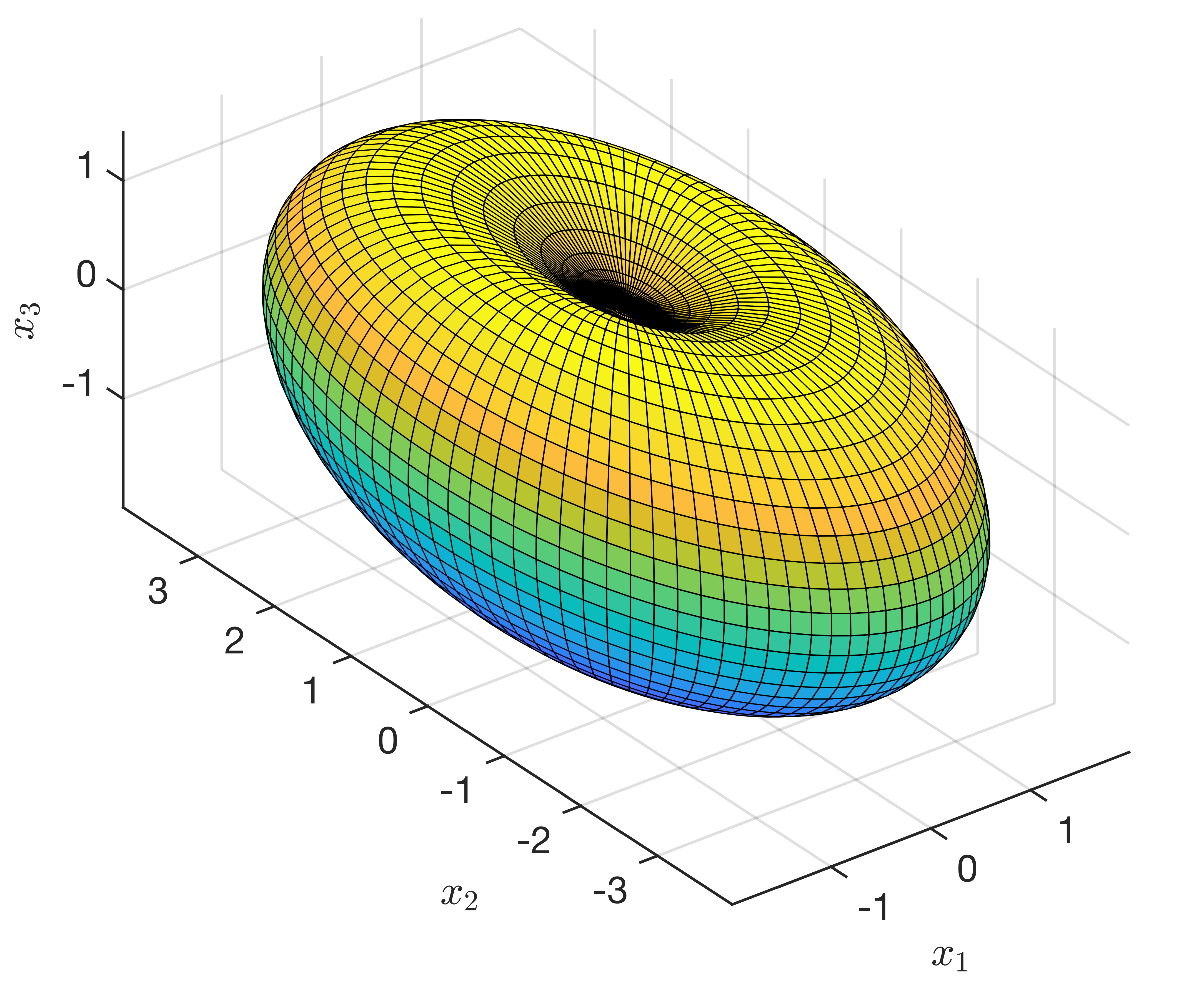}}
    \raisebox{-0.5\height}{\includegraphics[height=1.3in]{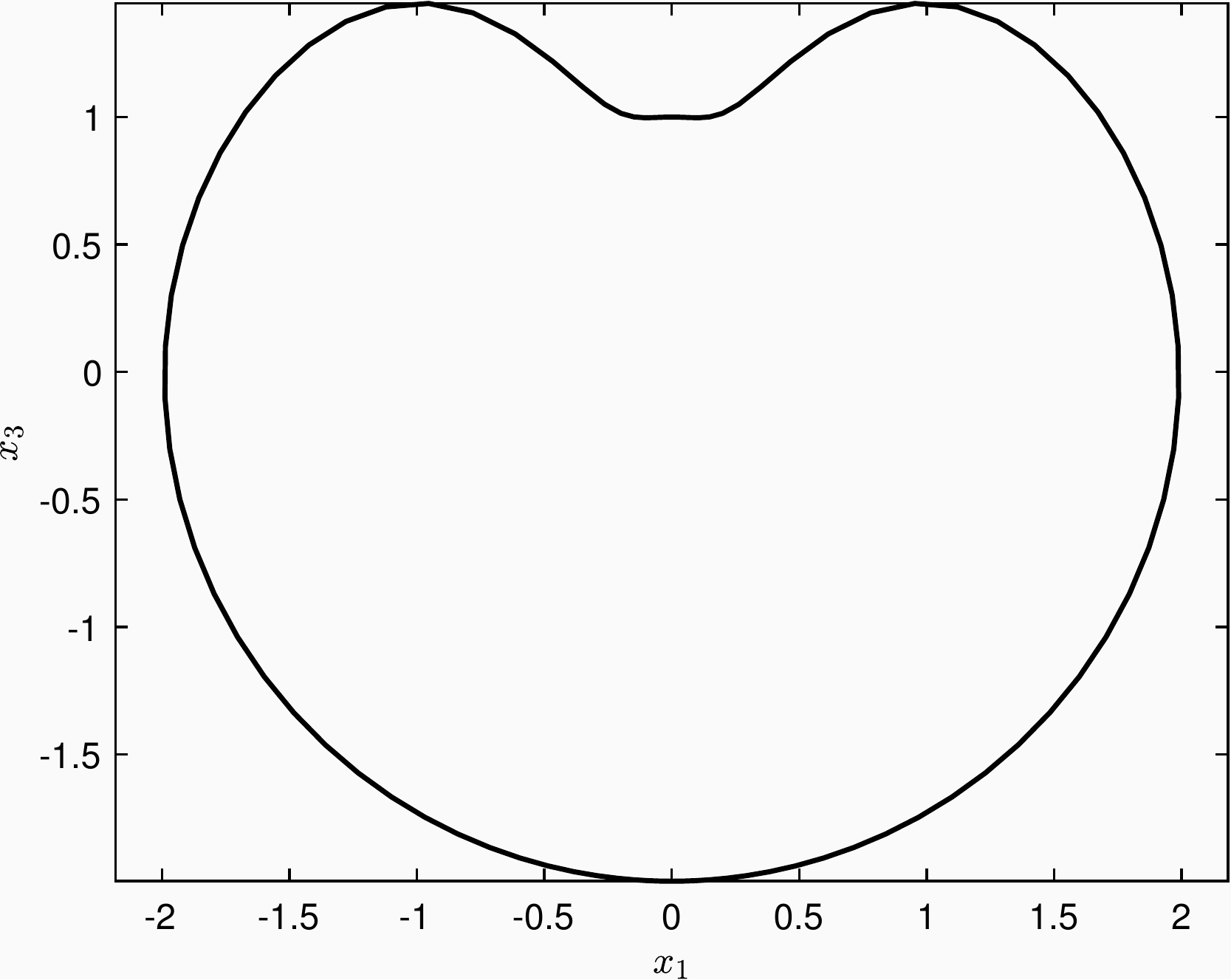}}
    \raisebox{-0.5\height}{\includegraphics[height=1.3in]{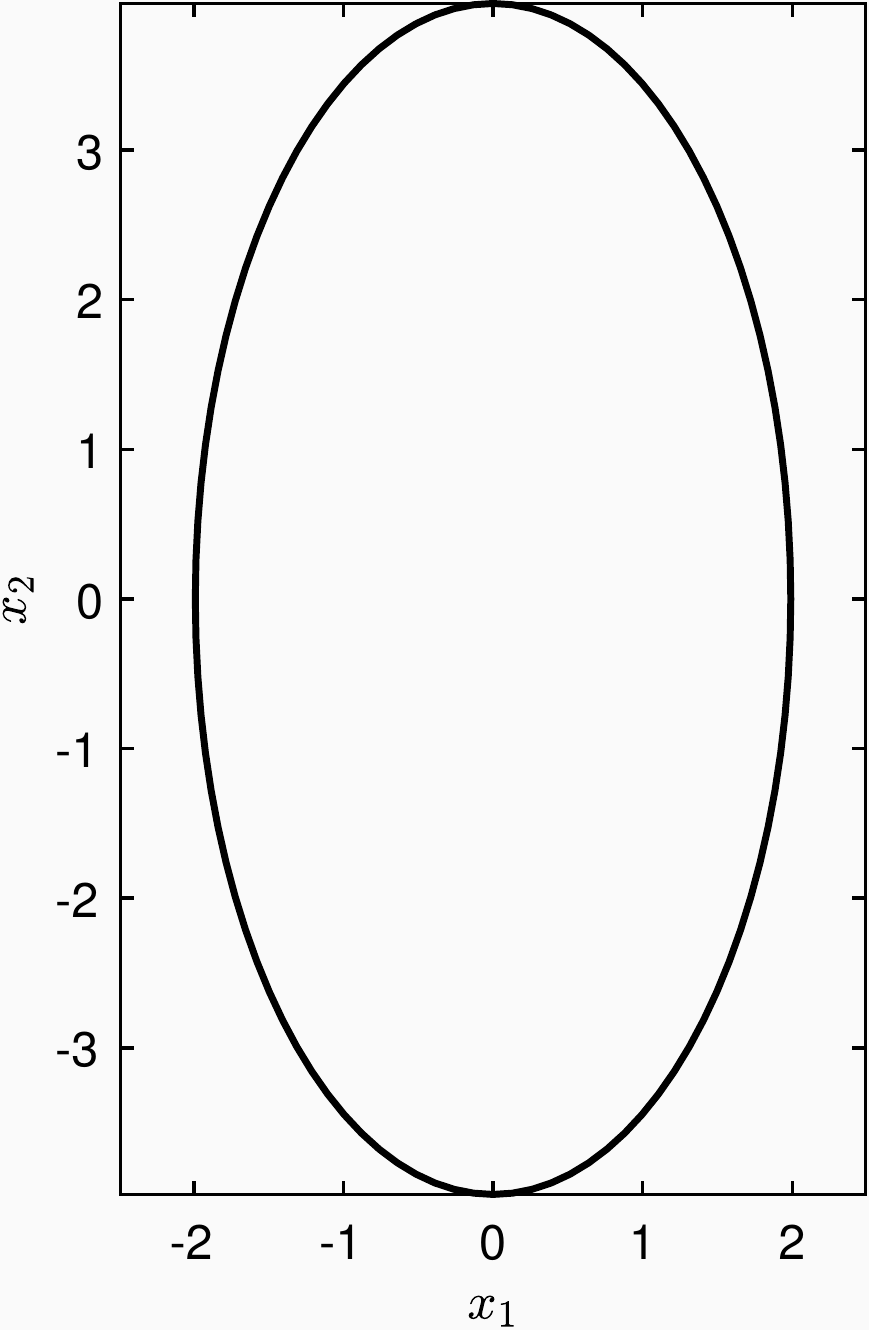}}
  \end{minipage}
  \caption{The boundary surface defined by
    \eqref{eq:mushroom}-\eqref{eq:mushroom2} that is used to exemplify
    the evaluation of the double layer potential in three dimensions
    (left),
    and the intersections of this boundary with the $x_{1}x_{3}$-plane
    (center), and $x_{1}x_{2}$-plane (right).}
  \label{mushroom}
\end{figure}

Let $(x_{1},x_{2},x_{3})$ denote an ordered triple in a Cartesian
coordinate system. To study the computation of the double-layer
potential in three dimensions, we consider the harmonic function,
\begin{equation}
  u(x_{1},x_{2},x_{3}) = \frac{1}{\sqrt{ (x_{1} - 5)^{2} + (x_{2} - 4 )^{2} +
    (x_{3} - 3)^{2}}}
\end{equation}
in the domain whose boundary is given by
\begin{equation}\label{eq:mushroom}
  y(\theta,\varphi) = R(\theta) ( \sin\theta \cos\varphi, 2 \sin\theta
  \sin\varphi, \cos\theta), \quad 0 \le \theta \le \pi, \quad -\pi \le
  \varphi \le \pi,
\end{equation}
with
\begin{equation}\label{eq:mushroom2}
  R(\theta) = 2 - \frac{1}{1 + 100 (1 - \cos\theta)^{2}}.
\end{equation}
This boundary surface is shown in Fig.~\ref{mushroom} (left) along
with its intersection with the vertical $x_{1}x_{3}$-plane (center)
and the horizontal $x_{1}x_{2}$-plane (right).

We solve boundary integral equation \eqref{eq:bie} using the Galerkin
method~\cite{atkinson1982laplace, atkinson1985algorithm,
  atkinson1990survey, atkinson1997numerical}. The Galerkin method
approximates the density according to
\begin{equation}
  \tilde{\mu}(\theta,\varphi) \approx \tilde{\mu}^{N}(\theta,\varphi)
  = \sum_{n = 0}^{N-1} \sum_{m = -n}^{n} \hat{\mu}_{nm}
  Y_{nm}(\theta,\varphi),
  \label{eq:mu-galerkin}
\end{equation}
with $\{ Y_{nm} \}$ denoting the orthonormal set of spherical
harmonics. For these results, we have set $N = 48$. 
We have computed the close evaluation of the double-layer potential at
points $x = y^{\star} - \epsilon \nu^{\star}$, using the following two
different methods for comparison.
\begin{enumerate}

\item {\bf Numerical approximation} -- Compute the modified
  double-layer potential,
  \begin{equation*}
    u(y^{\star} - \epsilon \nu^{\star}) = - \mu(y^{\star}) +
    \frac{1}{4\pi} \int_{B} \frac{\nu_{y} \cdot
      ( y^{\star} - \epsilon \nu^{\star} - y )}{| y^{\star} - \epsilon
      \nu^{\star} - y |^{3}} \left[ \mu(y) - \mu(y^{\star}) \right]
    \mathrm{d}\sigma_{y},
  \end{equation*}
  using the three-step numerical method given by Carvalho {\it et
    al.}~\cite{carvalho2018_3D}. For the first step, the modified
  double-layer potential is written as a spherical integral that has
  been rotated so that $y^{\star}$ is aligned with the north pole. For
  the second step, the $2N$-point PTR is used to compute the integral
  in the azimuthal angle. For the third step, the $N$-point
  Gauss-Legendre quadrature rule mapped to $[0,\pi]$ is used to
  compute the integral in the polar angle.

\item {\bf $O(\epsilon^{2})$ Asymptotic approximation} -- Compute the
  asymptotic approximation given by \eqref{eq:asymptotic3D} using the
  method given in \eqref{eq:numerical-method3D} (with $N = 48$).

\end{enumerate}

The error of the numerical method has been shown to decay
quadratically with $\epsilon$ when
$\epsilon \ll 1/N$~\cite{carvalho2018_3D}. This quadratic error decay
occurs because, in the rotated coordinate system, the azimuthal
integration acts as an averaging operation yielding a smooth function
of the polar angle that is computed to high order using Gaussian
quadrature. However, this asymptotic error estimate is valid only when
the numerical approximation of the density is sufficiently
resolved. If $N$ in \eqref{eq:mu-galerkin} is not sufficiently large
that $|\hat{\mu}_{nm}|$ for $n > N$ is negligibly small, then the
truncation error associated with \eqref{eq:mu-galerkin} may interrupt
this quadratic error decay. For the domain here, with $N = 48$, we
find that the estimated truncation error for \eqref{eq:mu-galerkin} is
approximately $10^{-8}$. While this error is relatively small, it is
not small enough to observe the error's quadratic rate of decay. We
would have to consider a much larger value of $N$ to observe that
decay rate. However, computing the numerical solution of boundary
integral equation \eqref{eq:bie} with $N > 48$ becomes restrictively
large. Hence, we evaluate below what the subtraction method and the
$O(\epsilon^{2})$ asymptotic approximation do in this limited
resolution situation.

\begin{figure}[!htb]
  \centering
  \includegraphics[width=0.44\linewidth]{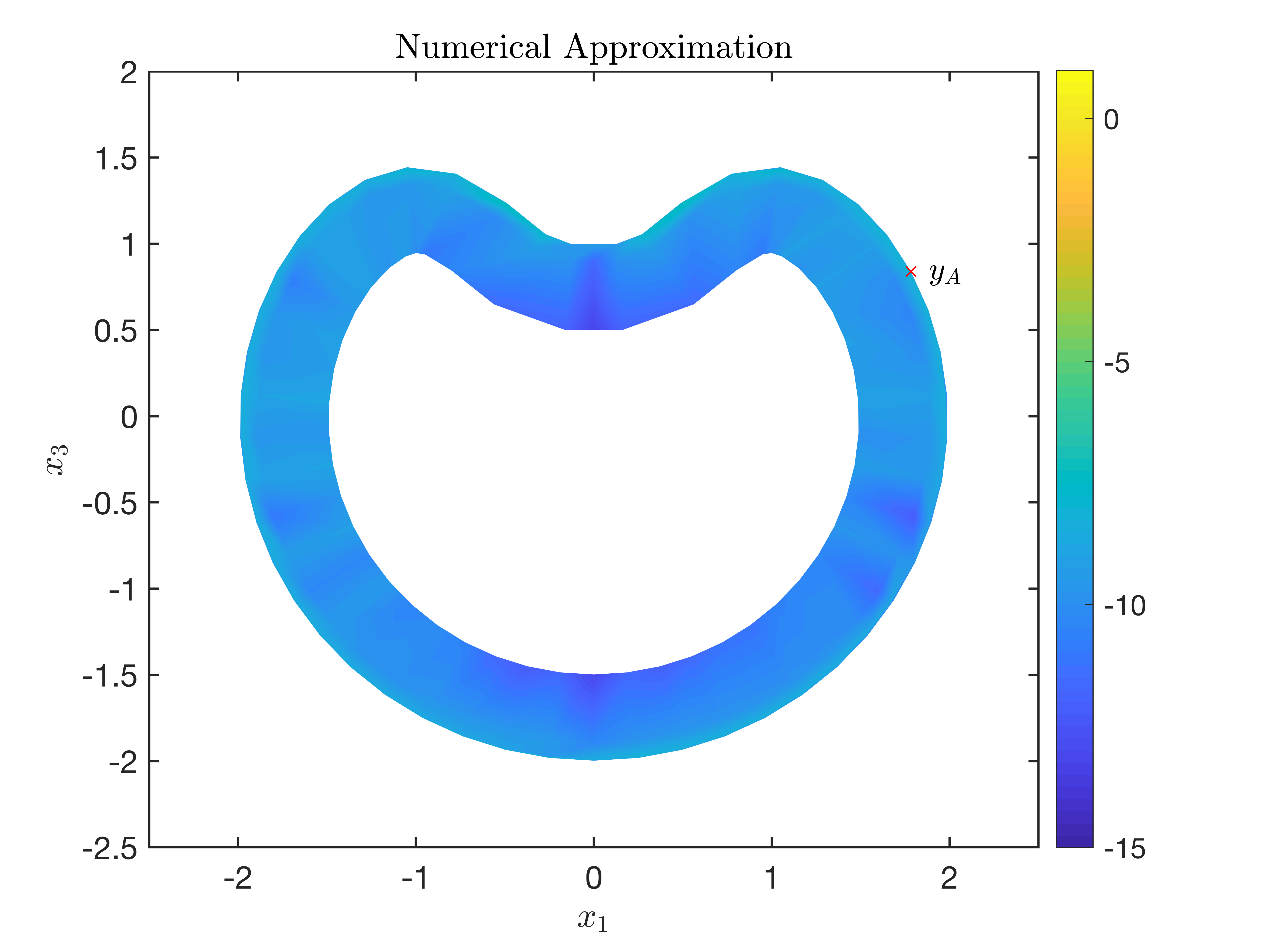}
  \includegraphics[width=0.44\linewidth]{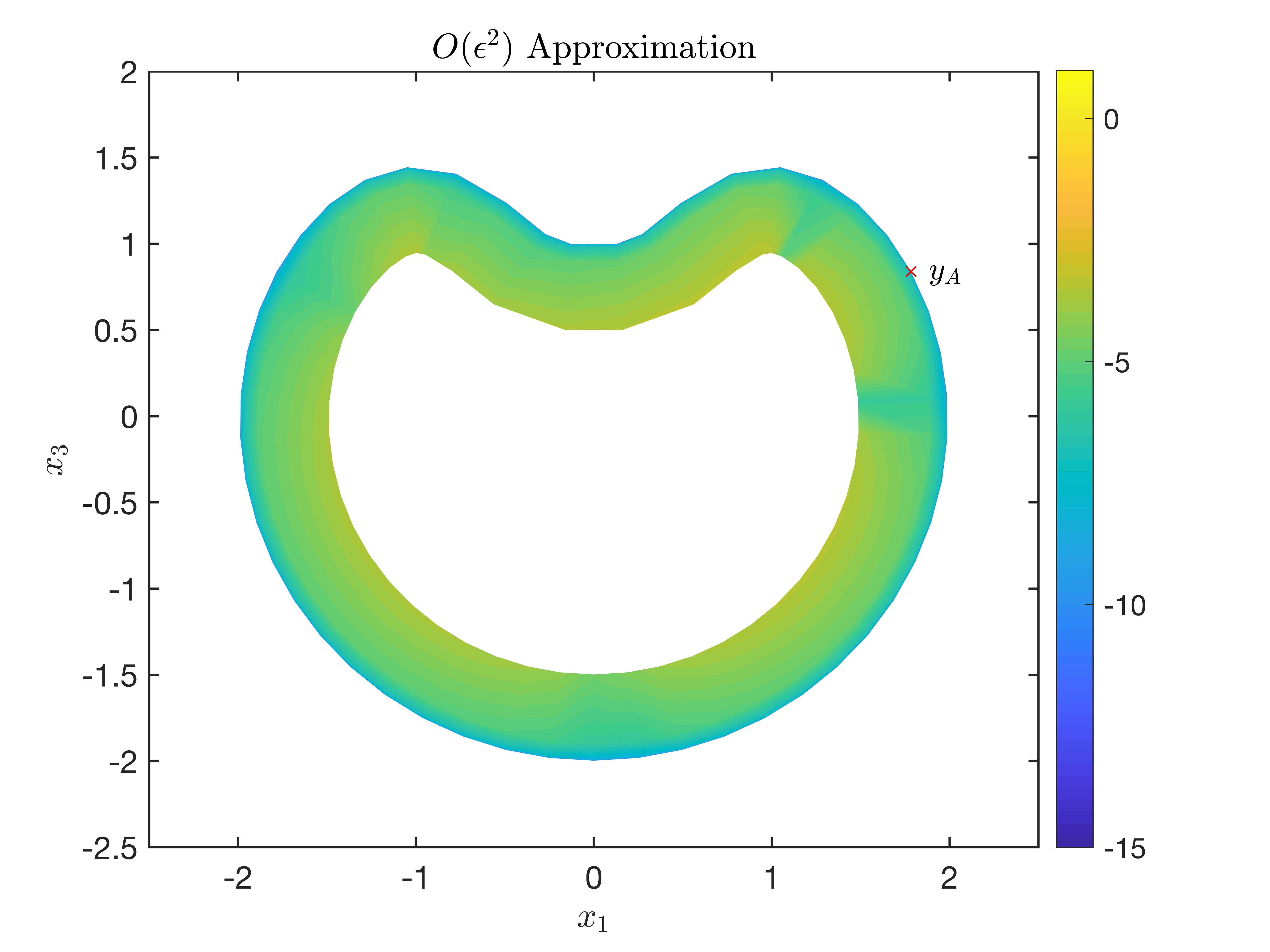}\\
  \includegraphics[width=0.44\linewidth]{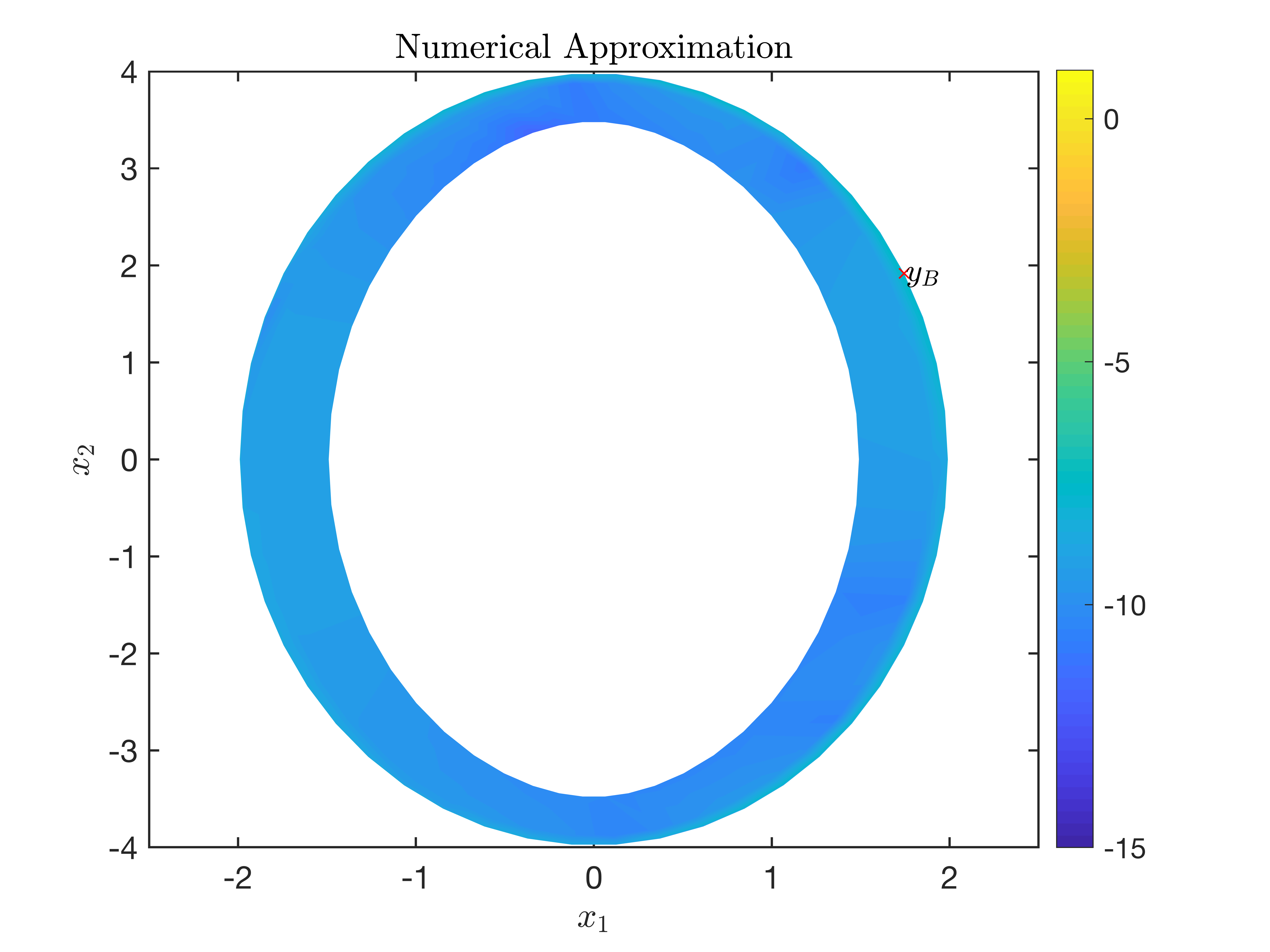}
  \includegraphics[width=0.44\linewidth]{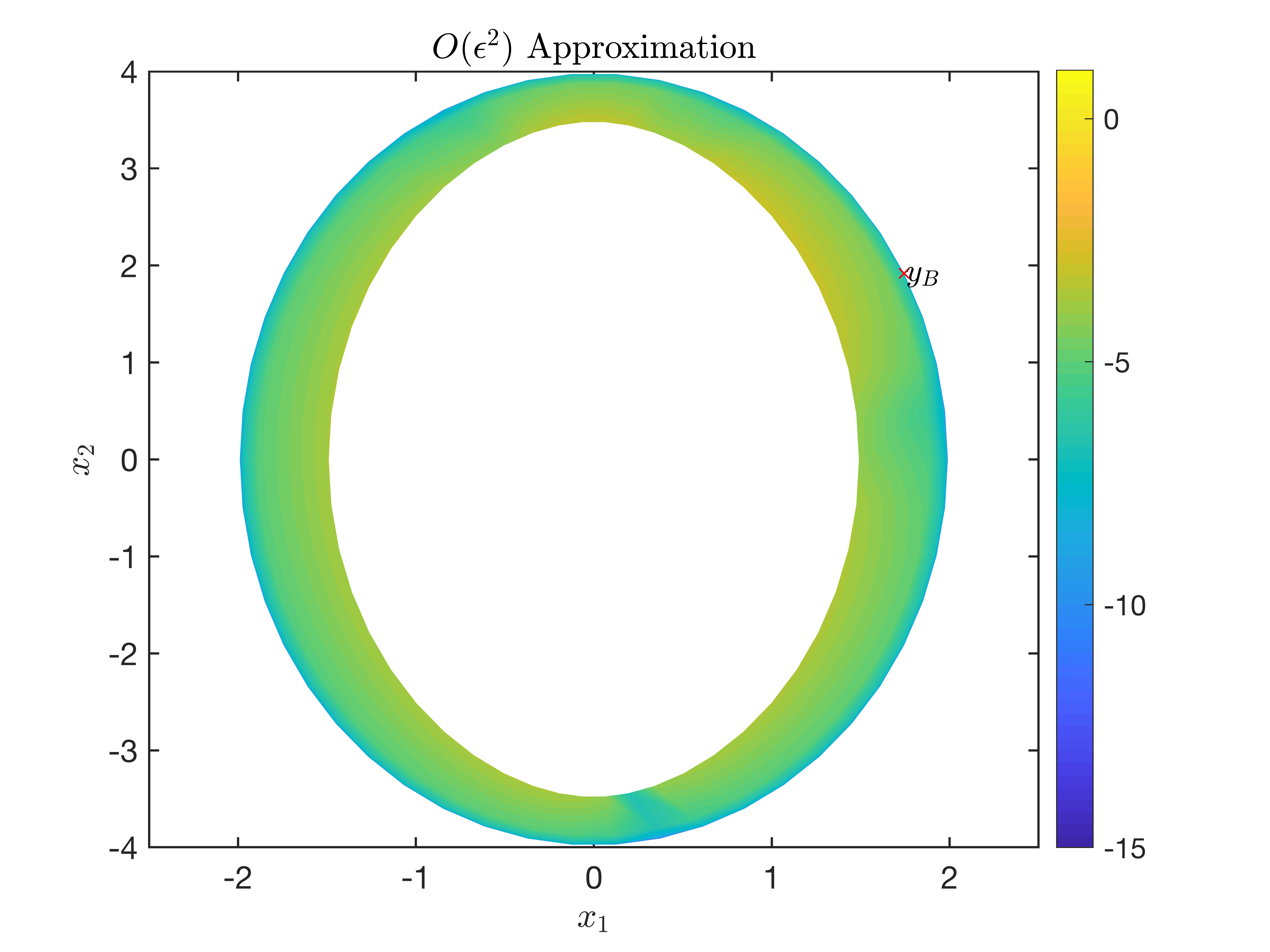}
  \caption{Plots of $\log_{10}$ of the error made in computing the
    double-layer potential by the numerical approximation (left
    column) and the $O(\epsilon^{2})$ asymptotic approximation (right
    column) in the domain whose boundary is shown in
    Fig.~\ref{mushroom}. The top row of plots show the error on the
    $x_{1}x_{3}$-plane, and the bottom row of plots show the error on
    the $x_{1}x_{2}$-plane. In the top row of plots, the point
    $y_{A} = ( 1.7830, 0, 0.8390 )$ is plotted as a red $\times$, and
    in the bottom row of plots, the point
    $y_{B} = ( 1.7439, 1.19175, 0 )$ is plotted as a red $\times$.}
  \label{mushroom-results}

  \bigskip

  \centering
  \includegraphics[width=0.48\linewidth]{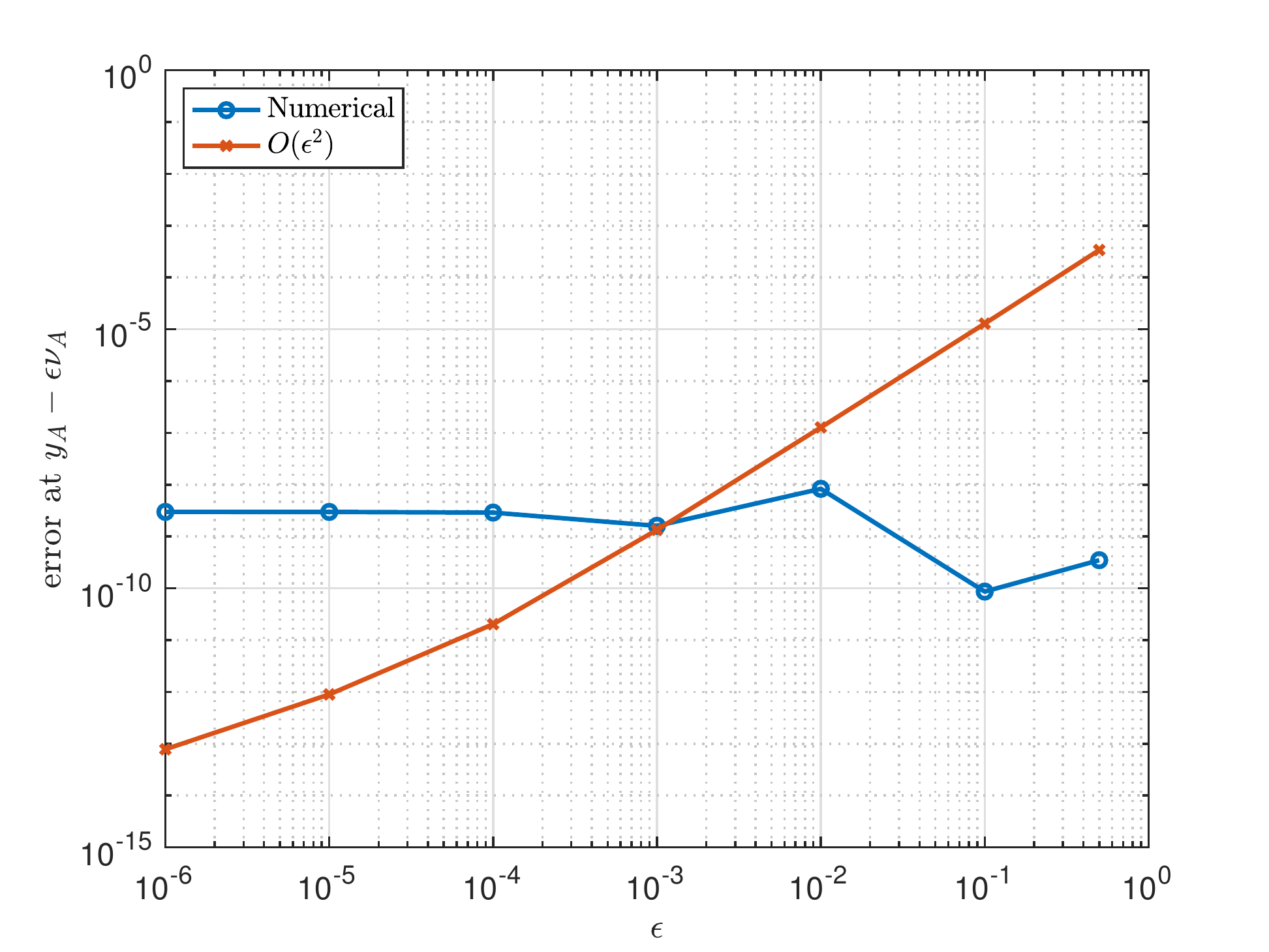}
  \includegraphics[width=0.48\linewidth]{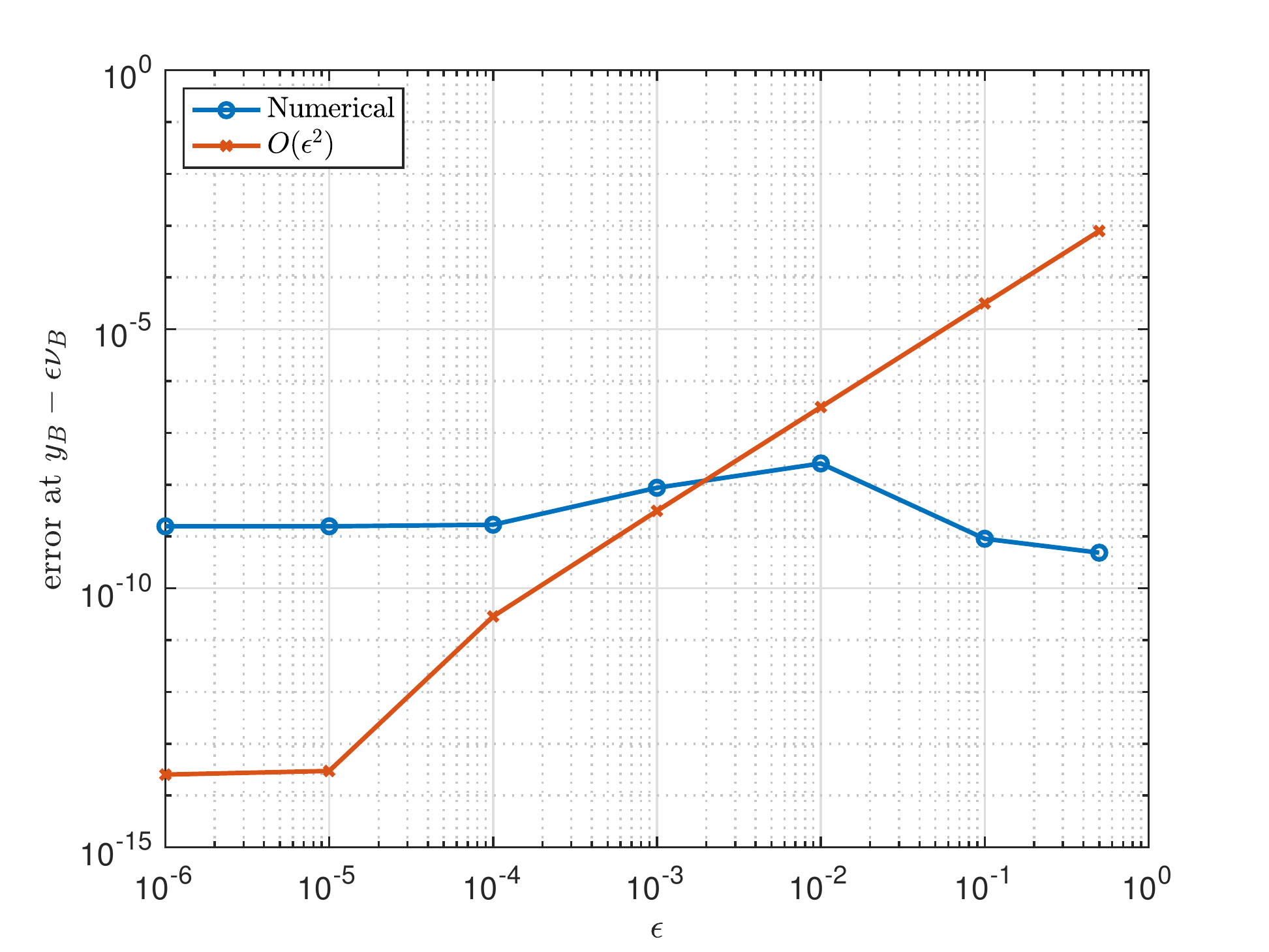}
  \caption{Log-log plots of the errors made by the two methods shown
    in Fig.~\ref{mushroom-results} at $y_{A} - \epsilon \nu_{A}$
    (left) and at $y_{B} - \epsilon \nu_{B}$ (right) for
    $10^{-6} \le \epsilon \le 0.5$.} 
  \label{mushroom-errors}
\end{figure}

Error results for the computation of the double-layer potential in this domain for each of the two methods described
above appear in Fig.~\ref{mushroom-results}. The top
row shows the error on the slice of the domain through the
vertical $x_{1}x_{3}$-plane for the numerical method (left) and
the $O(\epsilon^{2})$ asymptotic approximation method (right).The point $y_{A} = ( 1.7830, 0, 0.8390 )$ is
plotted as a red $\times$ symbol in both plots. The bottom row shows the
errors of the same methods (left for the numerical method, right for the $O(\epsilon^2)$ asymptotic approximation) on the slice of the domain through the horizontal
$x_{1}x_{2}$-plane. The point
$y_{B} = ( 1.7439, 1.19175, 0 )$ is plotted as a red $\times$ symbol in both plots.

In Fig.~\ref{mushroom-errors}, we show the errors computed at
$y_{A} - \epsilon \nu_{A}$ (left) and $y_{B} - \epsilon \nu_{B}$
(right) for varying $\epsilon$, where $\nu_{A}$ and $\nu_{B}$ are the
unit outward normals at $y_{A}$, and $y_{B}$, respectively. In
contrast to the two-dimensional results, we find that the error for
the numerical method is approximately $10^{-8}$ for all values of
$\epsilon$. This error is due to the truncation error made by the
Galerkin method. Because the truncation error dominates at this
resolution, we are not able to see its quadratic decay as
$\epsilon \to 0^{+}$. If a higher resolution computation was used to
solve the boundary integral equation, the error of the numerical
method would exhibit a similar behavior to that made by the
subtraction method for the two-dimensional examples. In particular,
the error would have a maximum at $\epsilon = O(1/N)$ about which the
error decays. 
We observe that
the $O(\epsilon^{2})$ asymptotic method decays monotonically with
$\epsilon$ even when $N = 48$.

\begin{figure}[htb]
  \centering
  \includegraphics[width=0.48\linewidth]{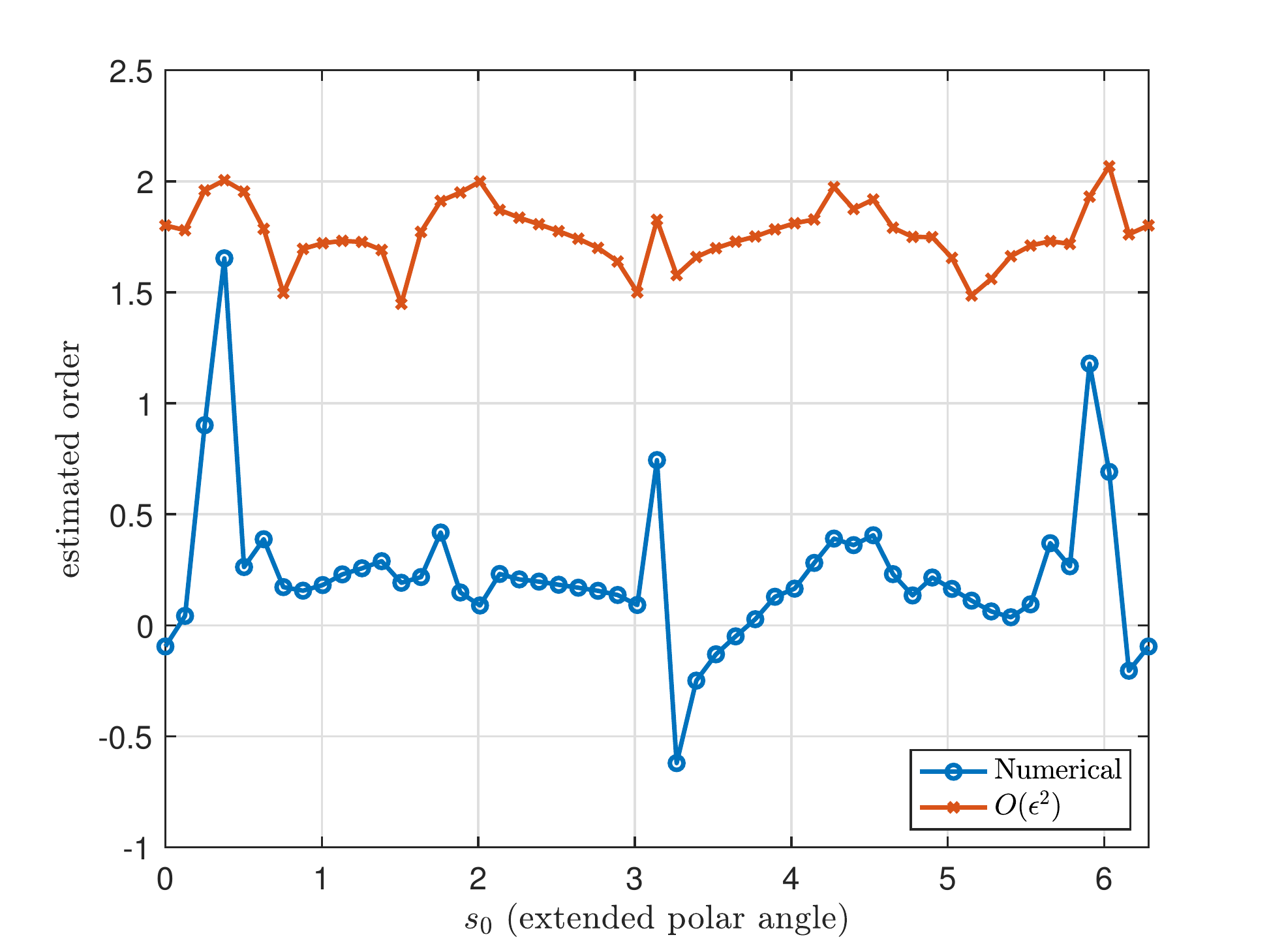}
  \includegraphics[width=0.48\linewidth]{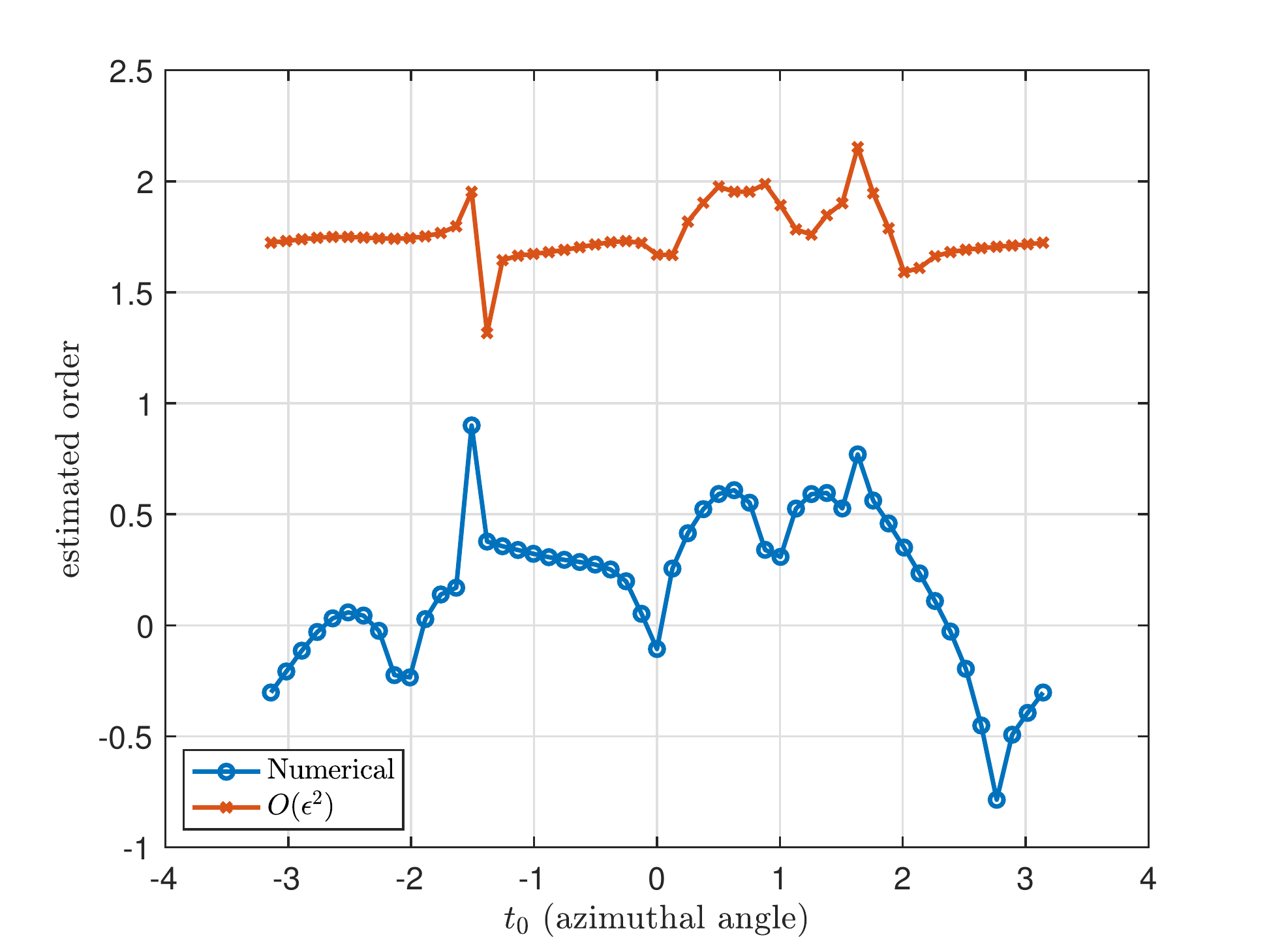}
  \caption{Estimated order of accuracy for the numerical
      approximation (blue $\circ$), and the $O(\epsilon^2)$ asymptotic
      approximation (red $\times$) on the $x_{1}x_{3}$-plane (left plot) and on
      the $x_{1}x_{2}$-plane (right plot). Results on the
      $x_{1}x_{3}$-plane are given in terms of the extended polar
      angle, $s_{0} \in [0,2\pi]$, which parameterizes the circle on
      the unit sphere lying on the $x_{1}x_{3}$-plane that starts and
      ends at the north pole. Results on the $x_{1}x_{2}$-plane are
      given in terms of the azimuthal angle, $t_{0} \in [0,2\pi]$.}
  \label{mushroom-order}
\end{figure}

We estimate the order of accuracy in Fig.~\ref{mushroom-order}.
The results for the esimated order of accuracy over
  the points intersecting the vertical $x_{1}x_{3}$-plane are shown in
  the left plot of Fig.~\ref{mushroom-order}. For those results we
  determine the estimated order of accuracy by determining the best
  fit line through the $\log-\log$ plot of the error versus $\epsilon$
  for several values of the extended polar angle,
  $s_{0} \in [0,2\pi]$. This extended polar angle parameterizes the
  circle on the unit sphere lying on the $x_{1}x_{3}$-plane that
  starts and ends at the north pole. The results for the esimated
  order of accuracy over the points intersecting the horizontal
  $x_{1}x_{2}$-plane are shown in the right plot of
  Fig.~\ref{mushroom-order}. For those results we determine the
  estimated order of accuracy by determining the best fit line through
  the $\log-\log$ plot of the error versus $\epsilon$ for several
  values of the azimuthal angle, $t_{0} \in [0,2\pi]$.
Because of the resolution limitation in the Galerkin method, we are
not able to see that the order of accuracy for the numerical method is
two. In fact, the error is nearly uniform with respect to $\epsilon$
because it is the truncation error of \eqref{eq:mu-galerkin} that is
dominating. Despite the resolution limitation in the Galerkin method,
we find that the $O(\epsilon^{2})$ asymptotic approximation has an
order accuracy of nearly two.

\subsubsection*{Summary of the results}
For three-dimensional problems, the subtraction method is more
effective when computed in an appropriate rotated coordinate system
than for two-dimensional problems. The subtraction method is more
effective in three dimensions because in this rotated coordinate
system, integration with respect to the azimuthal angle is a natural
averaging operation that regularizes the integral thereby allowing for
the use of a high-order quadrature rules for integration with respect
to the polar angle.  Provided that the density is sufficiently
resolved, the subtraction method has been shown to decay quadratically
with the distance away from the boundary~\cite{carvalho2018_3D}.  The
$O(\epsilon^{2})$ asymptotic approximation also decays
quadratically. However, it is not as sensitive to the accuracy of
the density. For this reason, we find that the asymptotic
approximation is still useful for three-dimensional problems,
especially when the density is not highly resolved. 

\section{Extension to forward-peaked scattering in radiative transfer}
\label{sec:RTE}
The asymptotic approximations developed here for the close
evaluation of the double-layer potential can be extended to other problems. 
In particular, there is an interesting connection between the close evaluation
problem in potential theory and forward-peaked
scattering in radiative transfer.  We first establish this connection
and then, we apply the asymptotic analysis used above to this problem.

Radiative transfer describes the multiple scattering of
light~\cite{chandrasekhar1960, ishimaru}. This theory has several
applications for light propagation and scattering in geophysical
media~\cite{marshak20053d,thomas2002radiative}, biological
tissues~\cite{wang2012biomedical}, and computer
graphics~\cite{jensen2001realistic}, among others. Let
$\psi: S^{2} \times \mathbb{R}^{3} \times (0,T) \to \mathbb{R}^{+} \ge 0$ denote the
specific intensity. The specific intensity gives the flow of power in
direction $\Omega \in S^{2}$, at position $r \in \mathbb{R}^{3}$, at time
$t \in [0,T]$. The radiative transfer equation,
\begin{equation}
  c^{-1} \psi_{t} + \Omega \cdot \nabla \psi + \mu_{a} \psi - \mu_{s}
  L \psi = Q,
  \label{eq:RTE}
\end{equation}
governs $\psi$ in a medium that absorbs, scatters, and emits
light. Here, $c$ denotes the speed of light in the background medium, $\mu_{a} \ge 0$, $\mu_{s} \ge 0$ denote absorption and
scattering coefficients, respectively, and $Q$ denotes a source. The
scattering operator, $L$, is defined by
\begin{equation}
  L \psi(\Omega,r,t) = \displaystyle \frac{1}{4 \pi} \int_{S^{2}} p(\Omega \cdot \Omega') \left[
    \psi(\Omega',r,t) - \psi(\Omega,r,t) \right]
  \mathrm{d}\sigma_{\Omega'},
  \label{eq:scattering-operator}
\end{equation}
where $p$ denotes the scattering phase function which gives the
fraction of light incident in direction $\Omega'$ that is scattered in
direction $\Omega$.

In forward-peaked scattering media, $p$ is sharply peaked about
$\Omega = \Omega'$ so that most of the scattering occurs in a small
angular cone about the incident direction. This problem is important
for several applications, especially light propagation in biological
tissues~\cite{Kim:03light}.  Mathematically, forward-peaked scattering
corresponds to the case in which \eqref{eq:scattering-operator} is a
nearly singular integral.  

There have been several studies on the asymptotic behavior of
\eqref{eq:RTE} with forward-peaked scattering leading to the
Fokker-Planck approximation~\cite{pomraning1992fokker} and its
generalizations~\cite{prinja2001generalized,
  leakeas2001generalized}. Solving the radiative transfer equation~\eqref{eq:RTE} with these approximate scattering operators has led to useful physical
insight~\cite{Kim:04backscattering, Kim:04beam,
  gonzalez2009comparison}. However, for most scattering models used
for applications, the Fokker-Planck approximation is
inaccurate. The error made by the Fokker-Planck
approximation is due to neglecting the tail away from the sharply peaked behavior of the kernel $p$, corresponding to large-angle scattering. That tail
typically decays too slowly to make a local approximation
appropriate. In fact, it has been shown that for several specific
scattering kernels, the leading order behavior is nonlocal and given
by a pseudo-differential operator~\cite{pomraning1996higher,
  larsen1999linear}.

We study forward-peaked scattering using the same asymptotic analysis
used to study the close evaluation of layer potentials discussed
above. For this study, we consider the specific choice of the
Henyey-Greenstein scattering phase function~\cite{henyey1941diffuse},
\begin{equation}
  p_{\text{HG}}(\Omega \cdot \Omega') =\frac{1 -
    g^{2}}{( 1 + g^{2} - 2 g \Omega \cdot \Omega')^{3/2}}, \quad -1
  \le g \le 1,
  \label{eq:HG}
\end{equation}
as the kernel.
The Henyey-Greenstein scattering phase function is extensively used in
applications since it provides a simple model with enough
sophistication to study multiple scattering of light. The anisotropy
factor, $g$, is the mean cosine of the scattering angle. It sets the
amount of scattering that is forward peaked. When $g = 0$, scattering is
istropic, and when $g = \pm 1$, scattering is restricted to only the
forward/backward direction. Forward-peaked scattering corresponds to
the asymptotic limit as $g \to 1$.

Henyey-Greenstein scattering is directly related to Poisson's formula
for the unit sphere,
\begin{equation}
  u(x) = \frac{1}{4\pi} \int_{S^{2}}
  \frac{1 - |x|^{2}}{| x - y |^{3}}
  f(y) \mathrm{d}\sigma_{y}, \quad |x| < 1,
  \label{eq:PoissonSphere}
\end{equation}
which gives the solution of the boundary value problem for
Laplace's equation:
\begin{subequations}
  \begin{gather}
    \varDelta u = 0 \quad \text{in $|x| < 1$},\\    
    u = f \quad \text{on $|x| = 1$}.  
  \end{gather}
  \label{eq:LaplaceSphere}
\end{subequations}
For the close evaluation point $x = (1 - \epsilon) y^{\star}$ with
$|y^{\star}| = 1$ and $0 < \epsilon \ll 1$, we find that
\begin{equation}
  u((1-\epsilon) y^{\star}) = \frac{1}{4\pi} \int_{S^{2}}
  \frac{1 - (1-\epsilon)^{2}}{\left[ 1 + (1-\epsilon)^{2} - 2 (1 -
      \epsilon) y^{\star} \cdot y \right]^{3/2}}
  f(y) \mathrm{d}\sigma_{y}.
  \label{eq:Poisson-CloseEval}
\end{equation}
Notice that the kernel in \eqref{eq:Poisson-CloseEval} is the same as
\eqref{eq:HG} with $g = 1 - \epsilon$.  Consequently, the asymptotic
limit of forward-peaked Henyey-Greenstein scattering is closely
related to the close evaluation of Poisson's formula for the unit sphere.

To evaluate $L$ given in \eqref{eq:scattering-operator} using
\eqref{eq:HG} as the scattering phase function, we use a spherical
coordinate system with $\Omega$ defining its north pole. In this coordinate system, $\theta$ denotes the polar angle and 
$\phi$ denotes the azimuthal angle.  For that
case, we have
\begin{equation}
  L \psi(\Omega) = \frac{1}{2}\int_{0}^{\pi}
  p_{\text{HG}}(\cos\theta;\epsilon) \tilde{\psi}(\theta) \sin \theta
  \mathrm{d}\theta
\end{equation}
with $ \Omega \cdot \Omega' = \cos\theta$, $g = 1 - \epsilon$, 
\begin{equation}
  p_{\text{HG}}(\cos\theta;\epsilon) =\frac{1 - (1 -
    \epsilon)^{2}}{[ 1 + (1 - \epsilon)^{2} - 2 (1 - \epsilon)
    \cos\theta ]^{3/2}},
  \label{eq:fHG}
\end{equation}
and
\begin{equation}
  \tilde{\psi}(\theta) = \frac{1}{2\pi} \int_{-\pi}^{\pi}
  \left[ \psi(\theta,\varphi) - \psi(0,\cdot) \right]
  \mathrm{d}\varphi.
\end{equation}
Due to the regularity of the solution at the north pole, we have
$\psi(\theta,\varphi + \pi) = \psi(-\theta,\varphi)$. As a result, we
write
\begin{equation}
  \tilde{\psi}(\theta) = \frac{1}{2\pi} \int_{0}^{\pi}
  \left[ \psi(\theta,\varphi) + \psi(-\theta,\varphi) - 2
    \psi(0,\cdot) \right] \mathrm{d}\varphi.
  \label{eq:u-tilde}
\end{equation}
Note here $\psi(0,\cdot) = \psi(\Omega)$. To study
forward-peaked scattering
we study the asymptotic limit corresponding to $\epsilon \to 0^{+}$.

We let
\begin{equation}
  L \psi(\Omega) = 
  L^{\text{in}} \psi(\Omega) +
  L^{\text{out}} \psi(\Omega),
\end{equation}
where
\begin{equation}
  L^{\text{in}} \psi(\Omega) =\frac{1}{2}\int_{0}^{\delta}
  p_{\text{HG}}(\cos\theta;\epsilon) \tilde{\psi}(\theta) \sin\theta
  \mathrm{d}\theta,
  \label{eq:Lpeak}
\end{equation}
is the inner expansion, and
\begin{equation}
  L^{\text{out}} \psi(\Omega) = \frac{1}{2} \int_{\delta}^{\pi}
  p_{\text{HG}}(\cos\theta;\epsilon) \tilde{\psi}(\theta) \sin\theta
  \mathrm{d}\theta,
  \label{eq:Ltail}
\end{equation}
is the outer expansion. Just as we have done for the double-layer
potential, we consider the asymptotic limit in which
$0 < \epsilon \ll \delta \ll 1$. Details of the calculations below can
be found in a developed \textit{Mathematica} notebook available on
GitHub~\cite{CKK-2018Codes}.

To find the leading-order behavior for $L^{\text{in}}$, we substitute
\eqref{eq:fHG} into \eqref{eq:Lpeak} and make the substitution
$\theta = \epsilon \Theta$, 
\begin{align}
  \begin{split}
    L^{\text{in}} \psi(\Omega) &= \int_{0}^{\delta/\epsilon}
    \frac{1}{2} \frac{1 - (1 - \epsilon)^{2}} {[1 + (1 - \epsilon)^{2}
      - 2 (1 - \epsilon) \cos\epsilon \Theta ]^{3/2}}
    \tilde{\psi}(\epsilon \Theta) \sin \epsilon\Theta
    \epsilon \mathrm{d}\Theta\\
    &= \int_{0}^{\delta/\epsilon} \left[ \frac{\Theta}{(1 +
        \Theta^{2})^{3/2}} - \frac{\epsilon}{2} \frac{\Theta - 2
        \Theta^{3}}{(1 + \Theta^{2})^{5/2}} + O(\epsilon^{2}) \right]
    \tilde{\psi}(\epsilon \Theta) \mathrm{d}\Theta.
  \end{split}
  \label{eq:5.14}
\end{align}
Next, we compute the expansion
\begin{align}
  \begin{split}
    \tilde{\psi}(\epsilon \Theta) &= \frac{1}{2\pi} \int_{0}^{\pi}
    \left[ \psi(\epsilon \Theta,\varphi) + \psi(-\epsilon
      \Theta,\varphi) - 2 \psi(0,\cdot) \right] \mathrm{d}\varphi\\
    &= \frac{\epsilon^{2} \Theta^{2}}{2\pi} \int_{0}^{\pi}
    \psi_{\theta\theta}(0,\cdot) \mathrm{d}\varphi + O(\epsilon^{4}).
  \end{split}
\end{align}
Substituting \eqref{eq:spherical-laplacian} (see Appendix
\ref{sec:spherical-laplacian}) into this result, we find that
\begin{equation}
  \tilde{\psi}(\epsilon \Theta) = \frac{\epsilon^{2}
    \Theta^{2}}{2 \pi} \int_{0}^{\pi} 
  \psi_{\theta\theta}(0,\cdot) \mathrm{d}\varphi + O(\epsilon^{4}) =
  \frac{\epsilon^{2} \Theta^{2}}{4}
  \varDelta_{S^{2}} \psi(\Omega) + O(\epsilon^{4}),
  \label{eq:sphericallaplacian}
\end{equation}
Substituting \eqref{eq:sphericallaplacian} into \eqref{eq:5.14}, we
determine by integrating and then expanding about $\epsilon = 0$ that
\begin{align}
  \begin{split}
    L^{\text{in}} \psi(\Omega) &= \varDelta_{S^{2}} \psi(\Omega)
    \int_{0}^{\delta/\epsilon} \left[ \frac{\epsilon^{2} \Theta^{3}}{4
        (1 + \Theta^{2})^{3/2}} - \frac{\epsilon^{3}}{8}
      \frac{\Theta^{3} - 2 \Theta^{5}}{(1 + \Theta^{2})^{5/2}} +
      O(\epsilon^{4}) \right]
    \mathrm{d}\Theta,\\
    &= \left[ \frac{\epsilon \delta}{4} - \frac{\epsilon^{2}}{2} +
      \frac{\epsilon^{2} \delta}{4} \right] \varDelta_{S^{2}} \psi(\Omega)
    + O(\epsilon^{3}).
  \end{split}
  \label{eq:L-inner}
\end{align}

To compute the outer expansion, we expand \eqref{eq:fHG} about
$\epsilon = 0$ then substitute the result into \eqref{eq:Ltail}, and
obtain
\begin{equation}
  L^{\text{out}} \psi(\Omega) = \int_{\delta}^{\pi} \frac{\epsilon +
    \epsilon^{2}}{2 \sqrt{2} (1 -
    \cos\theta)^{3/2}}  \tilde{\psi}(\theta) \sin\theta
  \mathrm{d}\theta + O(\epsilon^{3}).
\end{equation}
To remove $\delta$ from the lower limit of integration, we write
\begin{multline}
  L^{\text{out}} \psi(\Omega) = \int_{0}^{\pi} \frac{\epsilon +
    \epsilon^{2}}{2 \sqrt{2} (1 - \cos\theta)^{3/2}}
  \tilde{\psi}(\theta) \sin\theta \mathrm{d}\theta   - \int_{0}^{\delta} \frac{\epsilon + \epsilon^{2}}{2 \sqrt{2} (1 -
    \cos\theta)^{3/2}} \tilde{\psi}(\theta) \sin\theta
  \mathrm{d}\theta + O(\epsilon^{3}).
  \label{eq:Ltail-intermediate}
\end{multline}
For the second term in \eqref{eq:Ltail-intermediate}, we substitute
$\theta = \epsilon \Theta$ and find that
\begin{align}
  \begin{split}
    \int_{0}^{\delta/\epsilon} \frac{\epsilon + \epsilon^{2}}{2
      \sqrt{2} (1 - \cos\epsilon \Theta)^{3/2}} &
    \tilde{\psi}(\epsilon
    \Theta) \sin\epsilon\Theta \epsilon \mathrm{d}\Theta\\
    &= \int_{0}^{\delta/\epsilon} \left[ \frac{1 +
        \epsilon}{\Theta^{2}} + O(\epsilon^{2}) \right]
    \tilde{\psi}(\epsilon \Theta)
    \mathrm{d}\Theta\\
    &= \varDelta_{S^{2}} \psi(\Omega) \int_{0}^{\delta/\epsilon}
    \left[ \frac{\epsilon^{2} + \epsilon^{3}}{4} + O(\epsilon^{4})
    \right] \mathrm{d}\Theta\\
    &= \left[ \frac{\epsilon \delta}{4} + \frac{\epsilon^{2}
        \delta}{4} \right] \varDelta_{S^{2}} \psi(\Omega) +
    O(\epsilon^{3}),
  \end{split}
  \label{eq:5.17}
\end{align}
where we have made use of \eqref{eq:sphericallaplacian}
 Thus, the leading-order asymptotic behavior
for $L^{\text{out}} \psi(\Omega)$ is given by
\begin{equation}
  L^{\text{out}} \psi(\Omega) = \frac{\epsilon + \epsilon^{2}}{2
    \sqrt{2}} \int_{0}^{\pi} \frac{\tilde{\psi}(\theta)}{(1 -
    \cos\theta)^{3/2}} \sin\theta \mathrm{d}\theta - \left[
    \frac{\epsilon\delta}{4} + \frac{\epsilon^{2} \delta}{4} \right]
  \varDelta_{S^{2}} \psi(\Omega) + O(\epsilon^{3}).
  \label{eq:L-outer}
\end{equation}

By summing \eqref{eq:L-inner} and \eqref{eq:L-outer}, and
substituting \eqref{eq:u-tilde}, we find that
%
%
\begin{equation}
  L \psi(\Omega) = \left(\epsilon + \epsilon^{2}
  \right) L_{3/2} \psi(\Omega) - \frac{\epsilon^{2}}{2}
  \varDelta_{S^{2}} \psi(\Omega) + O(\epsilon^{3}).
  \label{eq:HG-asympt}
\end{equation}
where
\begin{equation}
  L_{3/2}\psi(\Omega)   = \frac{1}{4 \sqrt{2} \pi} \int_{0}^{\pi} \int_{0}^{\pi}
  \frac{\psi(\theta,\varphi) + \psi(-\theta,\varphi) - 2
    \psi(0,\cdot)}{(1 - \cos\theta)^{3/2}} \sin \theta
  \mathrm{d}\theta \mathrm{d}\varphi.
  \label{eq:L_3/2}
\end{equation}
%
The integral in \eqref{eq:L_3/2} appears to be singular, but since
\begin{equation}
  \frac{1}{4 \sqrt{2} \pi} \int_{0}^{\pi} \frac{\psi(\theta,\varphi) +
    \psi(-\theta,\varphi) - 2 \psi(0,\cdot)}{(1 - \cos\theta)^{3/2}}
  \sin\theta \mathrm{d}\varphi = \frac{1}{4} \varDelta_{S^{2}} \psi(\Omega)
  + O(\theta^{2}),
\end{equation}
$L_{3/2} \psi(\Omega)$ is well defined.

The leading-order asymptotic behavior of $L$ given in
\eqref{eq:HG-asympt} is equivalent to a result by
Larsen~\cite[Eq. 31]{larsen1999linear} who derived an asymptotic approximation
for $L$ using a spectral analysis.  In contrast to that asymptotic
approximation, the asymptotic analysis given above directly
addresses the balance between the forward peak given by the inner
expansion, and the long tail given by the outer expansion of the
scattering operator with the Henyey-Greenstein scattering
kernel. 

\section{Conclusion}
\label{sec:conclusion}

%
We have computed the leading-order asymptotic behavior for the close
evaluation of the double-layer potential in two and three
dimensions. By developing numerical methods to evaluate these
asymptotic approximations, we obtain effective methods for computing
double-layer potentials at close evaluation points. Our numerical
examples demonstrate the effectiveness of these asymptotic
approximations and corresponding numerical methods.

The key to this asymptotic analysis is the insight it provides. The
leading-order asymptotic behavior of the close evaluation of the
double-layer potential is given by its local Dirichlet data plus a
correction that is nonlocal. It is this nonlocal term that makes the
close evaluation problem challenging to address using only numerical
methods. It is consistent with the fact that solutions of boundary
value problems for elliptic partial differential equations have a
global dependence on their boundary data. By explicitly computing this
correction using asymptotic analysis, we have been able to develop an
effective numerical method for it. Moreover, the asymptotic error
estimates provide guidance on where to apply these approximations,
namely, for evaluation points closer to the boundary than the boundary
mesh spacing. The result of this work is an accurate and efficient
method for computing the close evaluation of the double-layer
potential.

The asymptotic approximation method can be extended to other layer potentials and other boundary value problems.  
In this paper, we show how these methods discover valuable insight for forward-peaked scattering in radiative transfer theory. 
Future work will entail of extending these methods to applications of Stokes flow and plasmonics.

\appendix

\section{Rotations on the sphere}
\label{sec:rotation}
We give the explicit rotation formulas over the sphere used in the
numerical method for the asymptotic approximation in three dimensions.
Consider $y, y^{\star} \in S^{2}$. We introduce the parameters
$\theta \in [0,\pi]$ and $\varphi \in [-\pi,\pi]$ and write
\begin{equation}
  y = y(\theta,\varphi) = \sin \theta \cos \varphi \, \ihat + \sin
  \theta \sin \varphi \, \jhat + \cos \theta \,
  \khat.
  \label{eq:y-xyz}
\end{equation}
The parameter values, $\theta^{\star}$ and $\varphi^{\star}$, are
set such that that $y^{\star} = y(\theta^{\star},\varphi^{\star})$.

We would like to work in the rotated, $\rm uvw$-coordinate system in which
\begin{equation}
  \begin{aligned}
    \uhat &= \cos \theta^{\star} \cos \varphi^{\star} \, \ihat +
    \cos \theta^{\star} \sin \varphi^{\star} \, \jhat - \sin
    \theta^{\star} \, \khat,\\
    \vhat &= - \sin \varphi^{\star} \, \ihat + \cos
    \varphi^{\star} \, \jhat,\\
    \what &= \sin \theta^{\star} \cos \varphi^{\star} \, \ihat +
    \sin \theta^{\star} \sin \varphi^{\star} \, \jhat + \cos
    \theta^{\star} \, \khat.
  \end{aligned}
  \label{eq:transformations}
\end{equation}
Notice that $\hat{\rm w} = y^{\star}$. For this rotated coordinate system,
we introduce the parameters $s \in [0,\pi]$ and $t \in [-\pi,\pi]$
such that
\begin{equation}
  y = y(s,t) = \sin s \cos t \, \uhat + \sin s
  \sin t \, \vhat + \cos s \, \what.
  \label{eq:y-uvw}
\end{equation}
It follows that $y^{\star} = y(0,\cdot)$.  By equating
\eqref{eq:y-xyz} and \eqref{eq:y-uvw} and substituting
\eqref{eq:transformations} into that result, we obtain
\begin{equation}
  \begin{bmatrix} \sin \theta \cos \varphi\\ \sin \theta \sin \varphi \\
    \cos \theta \end{bmatrix}
  = \begin{bmatrix} \cos \theta^{\star} \cos \varphi^{\star} & - \sin
    \varphi^{\star} & \sin \theta^{\star} \cos \varphi^{\star}\\
    \cos \theta^{\star} \sin \varphi^{\star} & \cos \varphi^{\star} & \sin
    \theta^{\star} \sin \varphi^{\star} \\
    - \sin \theta^{\star} & 0 & \cos \theta^{\star}
  \end{bmatrix}
  \begin{bmatrix} \sin s \cos t\\ \sin s \sin t\\
    \cos s \end{bmatrix}.
  \label{eq:y-identity}
\end{equation}
We rewrite \eqref{eq:y-identity} compactly as
$\hat{y}(\theta,\varphi) = R(\theta^{\star},\varphi^{\star})
\hat{y}(s,t)$
with $R(\theta^{\star},\varphi^{\star})$ denoting the $3 \times 3$ orthogonal 
rotation matrix.

We now seek to write $\theta = \theta(s,t)$ and
$\varphi = \varphi(s,t)$. To do so, we introduce
\begin{align}
  \xi(s,t;\theta^{\star},\varphi^{\star}) 
  &= \cos \theta^{\star} \cos \varphi^{\star} \sin s \cos t -
    \sin \varphi^{\star} \sin s \sin t + \sin \theta^{\star} \cos
    \varphi^{\star} \cos s,\\
  \eta(s,t,\theta^{\star},\varphi^{\star})
  &= \cos \theta^{\star} \sin \varphi^{\star} \sin s \cos t + \cos
    \varphi^{\star} \sin s \sin t + \sin \theta^{\star} \sin
    \varphi^{\star} \cos s,\\
  \zeta(s,t,\theta^{\star},\varphi^{\star})
  &= - \sin \theta^{\star} \sin s \cos t + \cos \theta^{\star} \cos s.
    \label{eq:angle-rotations}
\end{align}
From \eqref{eq:y-identity}, we find that 
\begin{equation}
  \theta = \arctan \left( \frac{\sqrt{\xi^{2} + \eta^{2}}}{\zeta}
  \right),
  \label{eq:theta-st}
\end{equation}
and
\begin{equation}
  \varphi = \arctan\left( \frac{\eta}{\xi} \right).
  \label{eq:varphi-st}
\end{equation}
With these formulas, we can write $\theta = \theta(s,t)$ and
$\varphi = \varphi(s,t)$.

\section{Spherical Laplacian}
\label{sec:spherical-laplacian}
In this Appendix, we establish the result given in
\eqref{eq:spherical-laplacian}.  We first seek an expression for
$\partial_{s}^{2}[ \cdot ]|_{s = 0}$ in terms of $\theta$ and
$\varphi$. By the chain rule, we find that
\begin{equation}
  \frac{\partial^{2}}{\partial s^{2}} [ \cdot ] \bigg|_{s = 0} =
  \left[ \left( \frac{\partial \theta}{\partial s} \right)^{2}
    \frac{\partial^{2}}{\partial \theta^{2}} + \left( \frac{\partial
        \varphi}{\partial s} \right)^{2} \frac{\partial^{2}}{\partial
      \varphi^{2}} + 2 \frac{\partial \theta}{\partial
      s} \frac{\partial \varphi}{\partial s}
    \frac{\partial^{2}}{\partial \theta \partial \varphi} +
    \frac{\partial^{2} \theta}{\partial s^{2}}
    \frac{\partial}{\partial \theta} + \frac{\partial^{2}
      \varphi}{\partial s^{2}} \frac{\partial}{\partial \varphi}
  \right] \bigg|_{s = 0}.
  \label{eq:chain-rule}
\end{equation}
Using $\theta$ defined in \eqref{eq:theta-st} and $\varphi$ defined in
\eqref{eq:varphi-st}, we find that
\begin{align}
  \frac{\partial \theta(s,t)}{\partial s} \bigg|_{s = 0} 
  &= \cos t, \label{eq:B.2}\\
  \frac{\partial^{2} \theta(s,t)}{\partial s^{2}} \bigg|_{s = 0} 
  &= \frac{\cos \theta^{\star}}{\sin\theta^{\star}} \sin^{2}
    t, \label{eq:B.3}\\
  \frac{\partial \varphi(s,t)}{\partial s} \bigg|_{s = 0} 
  &= \frac{\sin t}{\sin \theta^{\star}}, \label{eq:B.4}\\
  \frac{\partial^{2} \varphi(s,t)}{\partial s^{2}} \bigg|_{s = 0} 
  &= - \frac{\cos \theta^{\star}}{\sin^{2}\theta^{\star}} \sin 2 t.
    \label{eq:B.5}
\end{align}
Note that at $s = 0$, we have $\theta^{\star} = \theta$.

Substituting \eqref{eq:B.2} -- \eqref{eq:B.5} into
\eqref{eq:chain-rule} and replacing $\theta^{\star}$ by $\theta$, we
obtain
\begin{multline}
  \frac{\partial^{2}}{\partial s^{2}} [ \cdot ] \bigg|_{s = 0} =
  \cos^{2} t \frac{\partial^{2}}{\partial \theta^{2}} + \sin^{2} t
  \frac{1}{\sin^{2} \theta} \frac{\partial^{2}}{\partial \varphi^{2}}
  + 2 \cos t \sin t \frac{1}{\sin \theta} \frac{\partial^{2}}{\partial
    \theta \partial \varphi}\\
  + \sin^{2} t \frac{\cos \theta}{\sin \theta}
  \frac{\partial}{\partial \theta} - \sin 2 t \frac{\cos
    \theta}{\sin^{2} \theta} \frac{\partial}{\partial \varphi},
\end{multline}
from which it follows that
\begin{equation}
  \frac{1}{\pi} \int_{0}^{\pi} \frac{\partial^{2}}{\partial s^{2}} [
  \cdot ] \bigg|_{s = 0} \mathrm{d}t
  = \frac{1}{2} \left[ \frac{\partial^{2}}{\partial \theta^{2}} +
    \frac{\cos \theta}{\sin \theta} \frac{\partial}{\partial \theta} +
    \frac{1}{\sin^{2} \theta} \frac{\partial^{2}}{\partial
      \varphi^{2}} \right]
  = \frac{1}{2} \varDelta_{S^{2}}.
\end{equation}



\begin{thebibliography}{10}

\bibitem{af2016fast}
{\sc L.~af~Klinteberg and A.-K. Tornberg}, {\em A fast integral equation method
  for solid particles in viscous flow using quadrature by expansion}, J.
  Comput. Phys., 326 (2016), pp.~420--445.

\bibitem{af2017error}
{\sc L.~af~Klinteberg and A.-K. Tornberg}, {\em Error estimation for quadrature
  by expansion in layer potential evaluation}, Adv. Comput. Math., 43 (2017),
  pp.~195--234.

\bibitem{akselrod2014probing}
{\sc G.~M. Akselrod, C.~Argyropoulos, T.~B. Hoang, C.~Cirac{\`\i}, C.~Fang,
  J.~Huang, D.~R. Smith, and M.~H. Mikkelsen}, {\em Probing the mechanisms of
  large {P}urcell enhancement in plasmonic nanoantennas}, Nat. Photonics, 8
  (2014), pp.~835--840.

\bibitem{atkinson1982laplace}
{\sc K.~E. Atkinson}, {\em The numerical solution {L}aplace's equation in three
  dimensions}, SIAM J. Numer. Anal., 19 (1982), pp.~263--274.

\bibitem{atkinson1985algorithm}
{\sc K.~E. Atkinson}, {\em Algorithm 629: An integral equation program for
  {L}aplace's equation in three dimensions}, ACM Trans. Math. Softw., 11
  (1985), pp.~85--96.

\bibitem{atkinson1990survey}
{\sc K.~E. Atkinson}, {\em A survey of boundary integral equation methods for
  the numerical solution of {L}aplace's equation in three dimensions}, in
  Numerical {S}olution of {I}ntegral {E}quations, Springer, 1990, pp.~1--34.

\bibitem{atkinson1997numerical}
{\sc K.~E. Atkinson}, {\em The Numerical Solution of Integral Equations of the
  Second Kind}, Cambridge University Press, 1997.

\bibitem{barnett2015spectrally}
{\sc A.~Barnett, B.~Wu, and S.~Veerapaneni}, {\em Spectrally accurate
  quadratures for evaluation of layer potentials close to the boundary for the
  2d {S}tokes and {L}aplace equations}, SIAM J. Sci. Comput., 37 (2015),
  pp.~B519--B542.

\bibitem{barnett2014evaluation}
{\sc A.~H. Barnett}, {\em Evaluation of layer potentials close to the boundary
  for {L}aplace and {H}elmholtz problems on analytic planar domains}, SIAM J.
  Sci. Comput., 36 (2014), pp.~A427--A451.

\bibitem{beale2001method}
{\sc J.~T. Beale and M.-C. Lai}, {\em A method for computing nearly singular
  integrals}, SIAM J. Numer. Anal., 38 (2001), pp.~1902--1925.

\bibitem{beale2016simple}
{\sc J.~T. Beale, W.~Ying, and J.~R. Wilson}, {\em A simple method for
  computing singular or nearly singular integrals on closed surfaces}, Commun.
  Comput. Phys., 20 (2016), pp.~733--753.

\bibitem{carvalho2018asymptotic}
{\sc C.~Carvalho, S.~Khatri, and A.~D. Kim}, {\em Asymptotic analysis for close
  evaluation of layer potentials}, J. Comput. Phys., 355 (2018), pp.~327--341.

\bibitem{carvalho2018_3D}
{\sc C.~Carvalho, S.~Khatri, and A.~D. Kim}, {\em Close evaluation of layer
  potentials in three dimensions}, arXiv:1807.02474,  (2018).

\bibitem{chandrasekhar1960}
{\sc S.~Chandrasekhar}, {\em Radiative Transfer}, Dover Publications, 1960.

\bibitem{epstein2013convergence}
{\sc C.~L. Epstein, L.~Greengard, and A.~Kl{\"o}ckner}, {\em On the convergence
  of local expansions of layer potentials}, SIAM J. Numer. Anal., 51 (2013),
  pp.~2660--2679.

\bibitem{gonzalez2009comparison}
{\sc P.~Gonz{\'a}lez-Rodr{\'\i}guez and A.~D. Kim}, {\em Comparison of light
  scattering models for diffuse optical tomography}, Opt. Express, 17 (2009),
  pp.~8756--8774.

\bibitem{guenther1996partial}
{\sc R.~B. Guenther and J.~W. Lee}, {\em Partial Differential Equations of
  Mathematical Physics and Integral Equations}, Dover Publications, 1996.

\bibitem{helsing2008evaluation}
{\sc J.~Helsing and R.~Ojala}, {\em On the evaluation of layer potentials close
  to their sources}, J. Comput. Phys., 227 (2008), pp.~2899--2921.

\bibitem{henyey1941diffuse}
{\sc L.~G. Henyey and J.~L. Greenstein}, {\em Diffuse radiation in the galaxy},
  Astrophys. J., 93 (1941), pp.~70--83.

\bibitem{ishimaru}
{\sc A.~Ishimaru}, {\em Wave Propagation and Scattering in Random Media}, IEEE
  Press, Piscataway, NJ, 1997.

\bibitem{jensen2001realistic}
{\sc H.~W. Jensen}, {\em Realistic Image Synthesis Using Photon Mapping}, AK
  Peters, Ltd., 2001.

\bibitem{keaveny2011applying}
{\sc E.~E. Keaveny and M.~J. Shelley}, {\em Applying a second-kind boundary
  integral equation for surface tractions in {S}tokes flow}, J. Comput. Phys.,
  230 (2011), pp.~2141--2159.

\bibitem{CKK-2018Codes}
{\sc A.~D. Kim}, {\em Asymptotic-{DLP}}.
\texttt{https://github.com/arnolddkim/Asymptotic-DLP}, 2018.

\bibitem{Kim:03light}
{\sc A.~D. Kim and J.~B. Keller}, {\em Light propagation in biological tissue},
  J. Opt. Soc. Am. A, 20 (2003), pp.~92--98.

\bibitem{Kim:04backscattering}
{\sc A.~D. Kim and M.~Moscoso}, {\em Backscattering of beams by forward-peaked
  scattering media}, Opt. Lett., 29 (2004), pp.~74--76.

\bibitem{Kim:04beam}
{\sc A.~D. Kim and M.~Moscoso}, {\em Beam propagation in sharply peaked forward
  scattering media}, J. Opt. Soc. Am. A, 21 (2004), pp.~797--803.

\bibitem{klockner2013quadrature}
{\sc A.~Kl{\"o}ckner, A.~Barnett, L.~Greengard, and M.~O'Neil}, {\em Quadrature
  by expansion: A new method for the evaluation of layer potentials}, J.
  Comput. Phys., 252 (2013), pp.~332--349.

\bibitem{larsen1999linear}
{\sc E.~W. Larsen}, {\em The linear {B}oltzmann equation in optically thick
  systems with forward-peaked scattering}, Prog. Nucl. Energy, 34 (1999),
  pp.~413--423.

\bibitem{leakeas2001generalized}
{\sc C.~L. Leakeas and E.~W. Larsen}, {\em Generalized {F}okker-{P}lanck
  approximations of particle transport with highly forward-peaked scattering},
  Nucl. Sci. Eng., 137 (2001), pp.~236--250.

\bibitem{Maier07}
{\sc S.~A. Maier}, {\em Plasmonics: {F}undamentals and {A}pplications},
  Springer, 2007.

\bibitem{marple2016fast}
{\sc G.~R. Marple, A.~Barnett, A.~Gillman, and S.~Veerapaneni}, {\em A fast
  algorithm for simulating multiphase flows through periodic geometries of
  arbitrary shape}, SIAM J. Sci. Comput., 38 (2016), pp.~B740--B772.

\bibitem{marshak20053d}
{\sc A.~Marshak and A.~Davis}, {\em 3D {R}adiative {T}ransfer in {C}loudy
  {A}tmospheres}, Springer Science \& Business Media, 2005.

\bibitem{mayer2008label}
{\sc K.~M. Mayer, S.~Lee, H.~Liao, B.~C. Rostro, A.~Fuentes, P.~T. Scully,
  C.~L. Nehl, and J.~H. Hafner}, {\em A label-free immunoassay based upon
  localized surface plasmon resonance of gold nanorods}, ACS Nano, 2 (2008),
  pp.~687--692.

\bibitem{novotny2011antennas}
{\sc L.~Novotny and N.~Van~Hulst}, {\em Antennas for light}, Nat. Photonics, 5
  (2011), pp.~83--90.

\bibitem{pomraning1992fokker}
{\sc G.~C. Pomraning}, {\em The {F}okker-{P}lanck operator as an asymptotic
  limit}, Math. Models Methods Appl. Sci., 2 (1992), pp.~21--36.

\bibitem{pomraning1996higher}
{\sc G.~C. Pomraning}, {\em Higher order {F}okker-{P}lanck operators}, Nucl.
  Sci. Eng., 124 (1996), pp.~390--397.

\bibitem{prinja2001generalized}
{\sc A.~K. Prinja and G.~C. Pomraning}, {\em A generalized {F}okker-{P}lanck
  model for transport of collimated beams}, Nucl. Sci. Eng., 137 (2001),
  pp.~227--235.

\bibitem{rachh2017fast}
{\sc M.~Rachh, A.~Kl{\"o}ckner, and M.~O'Neil}, {\em Fast algorithms for
  quadrature by expansion i: Globally valid expansions}, J. Comput. Phys., 345
  (2017), pp.~706--731.

\bibitem{sannomiya2008situ}
{\sc T.~Sannomiya, C.~Hafner, and J.~Voros}, {\em In situ sensing of single
  binding events by localized surface plasmon resonance}, Nano Lett., 8 (2008),
  pp.~3450--3455.

\bibitem{schwab1999extraction}
{\sc C.~Schwab and W.~Wendland}, {\em On the extraction technique in boundary
  integral equations}, Math. Comput., 68 (1999), pp.~91--122.

\bibitem{smith2009boundary}
{\sc D.~J. Smith}, {\em A boundary element regularized {S}tokeslet method
  applied to cilia-and flagella-driven flow}, Proc. R. Soc. Lond. A, 465
  (2009), pp.~3605--3626.

\bibitem{thomas2002radiative}
{\sc G.~E. Thomas and K.~Stamnes}, {\em Radiative {T}ransfer in the
  {A}tmosphere and {O}cean}, Cambridge University Press, 2002.

\bibitem{wala20183DQBX}
{\sc M.~Wala and A.~Kl{\"o}ckner}, {\em A fast algorithm for {Q}uadrature by
  {E}xpansion in three dimensions}, arXiv:1805.06106v1,  (2018).

\bibitem{wang2012biomedical}
{\sc L.~V. Wang and H.-I. Wu}, {\em Biomedical Optics: Principles and Imaging},
  John Wiley \& Sons, 2012.

\end{thebibliography}
\end{document}